\documentclass[a4paper,twoside,10pt]{article}

\usepackage[a4paper,left=3cm,right=3cm, top=3cm, bottom=3cm]{geometry}
\usepackage[latin1]{inputenc}
\usepackage{mathrsfs}
\usepackage{graphicx}
\usepackage{epstopdf}
\usepackage{amsmath}
\usepackage{amsthm}
\usepackage{amssymb}
\usepackage{hyperref}
\usepackage{stmaryrd}
\usepackage{float}
\usepackage{bigints}
\usepackage{cite}
\usepackage{color}
\usepackage[abs]{overpic}
\usepackage[font=footnotesize,labelfont=bf]{caption}
\usepackage{cases}
\usepackage{tikz}
\usepackage{algorithmic}
\usepackage{rotating}
\usepackage{blkarray}
\usetikzlibrary{matrix,calc,arrows}
\usepackage{pdfcomment}
\usepackage{soul,xcolor}
\usepackage{enumitem}  
\usepackage{verbatim}
\usepackage{graphicx}
\usepackage{subcaption}
\usepackage{mwe}
\usepackage{algorithm}
\usepackage{pifont}
\usepackage{multirow}
\usepackage{datetime}

\restylefloat{table}
\theoremstyle{plain}
\newtheorem{thm}{Theorem}[section]
\newtheorem{cor}[thm]{Corollary}
\newtheorem{lem}[thm]{Lemma}

\theoremstyle{definition}

\theoremstyle{remark}
\newtheorem{remark}{Remark}

\newcommand{\ap}[1]{{\color{red}{#1}}} %colors
\newcommand{\N}{\mathbb{N}}
\newcommand{\R}{\mathbb{R}}

\renewcommand{\div}{\operatorname{\rm div}}

\newcommand{\rot}{\operatorname{\rm rot}}

\newcommand{\diam}{\mathrm{diam}}

\newcommand {\bu}{\boldsymbol{u}}
\newcommand {\buh}{\boldsymbol{u}_h}
\newcommand {\bv}{\boldsymbol{v}}
\newcommand {\bvh}{\boldsymbol{v}_h}
\newcommand {\bx}{\boldsymbol{x}}
\newcommand {\by}{\boldsymbol{y}}

\newcommand {\bn}{\boldsymbol{n}}

\newcommand {\bsigma}{\boldsymbol{\sigma}}

\newcommand {\bmu}{\boldsymbol{\mu}}
\newcommand {\bgamma}{\boldsymbol{\gamma}}

\newcommand {\bzetah}{\boldsymbol{\zeta}_h}
\newcommand {\bmuh}{\boldsymbol{\mu}_h}
\newcommand {\ch}{c_h}
\newcommand {\zh}{z_h}
\newcommand {\ph}{p_h}
\newcommand {\qh}{q_h}
\newcommand {\Thetah}{\Theta_h}
\newcommand {\czh}{c_{0,h}}
\newcommand {\Gammach}{\Gamma_{c,h}}
\newcommand {\Gammac}{\Gamma_c}

\newcommand {\qplus}{q^+}
\newcommand {\qminus}{q^-}

\newcommand {\tn}{t_n}
\newcommand {\tnmo}{t_{n-1}}

\newcommand {\bun}{\boldsymbol{u}^n}

\newcommand {\buno}{\boldsymbol{u}^{n-1}}
\newcommand {\bunI}{\boldsymbol{u}^n_I}

\newcommand {\cn}{c^n}
\newcommand {\cno}{c^{n-1}}

\newcommand {\pn}{p^n}
\newcommand {\bUn}{\boldsymbol{U}^n}
\newcommand {\bUno}{\boldsymbol{U}^{n-1}}
\newcommand {\Cn}{C^n}
\newcommand {\Cno}{C^{n-1}}

\newcommand {\Pn}{P^n}
\newcommand {\pnI}{p^n_I}
\newcommand {\chat}{\widehat{c}}
\newcommand {\chatn}{\widehat{c}^{\, n}}
\newcommand {\Gn}{G^n}
\newcommand {\qplusn}{(q^+)^n}
\newcommand {\qminusn}{(q^-)^n}

\newcommand {\bdeltan}{\boldsymbol{\delta}^n}
\newcommand {\Phatn}{\widehat{P}^n}
\newcommand {\Whatn}{\widehat{\boldsymbol{W}}^n}
\newcommand {\Pc}{\mathcal{P}_c }
\newcommand {\varthetao}{\vartheta^0}
\newcommand {\varthetan}{\vartheta^n}
\newcommand {\varthetano}{\vartheta^{n-1}}
\newcommand {\rhon}{\rho^n}
\newcommand {\rhono}{\rho^{n-1}}

\newcommand{\E}{K} % generic polygon
\newcommand{\taun}{\mathcal T_h} % polygonal decomposition
\newcommand{\e}{e} % generic edge
\newcommand{\hE}{h_{\E}} % mesh size function of \E
\newcommand{\he}{h_e} % length of edge e
\newcommand{\En}{\mathcal E_h} % edges
\newcommand{\EE}{\mathcal E^{\E}} % edges of element K

\newcommand{\dof}{\textup{dof}} % degrees of freedom
\newcommand{\dofVh}{\dof^{\boldsymbol{V}_h}} % degrees of freedom
\newcommand{\dofVhE}{\dof^{\boldsymbol{V}_h(K)}} % degrees of freedom
\newcommand{\dofQh}{\dof^{Q_h}} % degrees of freedom
\newcommand{\dofQhE}{\dof^{Q_h(K)}} % degrees of freedom
\newcommand{\dofZhtildeE}{\dof^{\widetilde{Z_h}(K)}} % degrees of freedom
\newcommand{\dofZh}{\dof^{Z_h}} % degrees of freedom
\newcommand{\dofZhE}{\dof^{Z_h(K)}} % degrees of freedom
\newcommand{\NVh}{{\textrm{dim}{V}_h}} % number of dofs
\newcommand{\NVhE}{{\textrm{dim}{V}_h(K)}} % number of dofs
\newcommand{\NQh}{\textrm{dim}{Q_h}} % number of dofs
\newcommand{\NQhE}{\textrm{dim}{Q_h(K)}} % number of dofs
\newcommand{\NZhtildeE}{\textrm{dim}{\widetilde{Z_h}(K)}} % number of dofs
\newcommand{\NZh}{\textrm{dim}{Z_h}} % number of dofs
\newcommand{\NZhE}{\textrm{dim}{Z_h(K)}} % number of dofs

\newcommand{\Ne}{n_\E} % number of edges of element K
\newcommand{\philE}{\varphi_j^\E} % local canonical basis fcts
\newcommand{\psilE}{\boldsymbol{\psi}_j^\E} % local canonical basis fcts

\newcommand{\VhE}{\boldsymbol{V}_h(\E)} 
\newcommand{\Vh}{\boldsymbol{V}_h} 
\newcommand{\QhE}{Q_h(\E)} 
\newcommand{\Qh}{Q_h} 
\newcommand{\ZhE}{Z_h(\E)} 
\newcommand{\ZhEtilde}{\widetilde{Z_h}(\E)} 
\newcommand{\Zh}{Z_h} 
\newcommand{\HdivE}{H(\div;\E)}
\newcommand{\HrotE}{H(\rot;\E)}
\newcommand{\HdivOmega}{H(\div;\Omega)}

\newcommand{\Acal}{\mathcal{A}}

\newcommand{\Acalh}{\mathcal{A}_h}
\newcommand{\AcalhE}{\mathcal{A}_h^\E}
\newcommand{\Dcal}{\mathcal{D}}
\newcommand{\Dcalh}{\mathcal{D}_h}
\newcommand{\DcalhE}{\mathcal{D}_h^\E}
\newcommand{\Khcal}{\mathcal{K}_h}
\newcommand{\Kcal}{\mathcal{K}}
\newcommand{\Mcal}{\mathcal{M}}
\newcommand{\Mcalh}{\mathcal{M}_h}
\newcommand{\McalhE}{\mathcal{M}_h^\E}

\newcommand{\ds}{\textup{d}s}
\newcommand{\dx}{\textup{d}x}

\newcommand{\SEA}{S^\E_{\Acal}} % local stabilization
\newcommand{\nuAE}{\nu_{\Acal}^\E} % local stabilization
\newcommand{\nuAplus}{\nu_{\Acal}^+} % stabilization constant
\newcommand{\nuAminus}{\nu_{\Acal}^-} % stabilization constant
\newcommand{\SEM}{S^\E_{\Mcal}} % local stabilization
\newcommand{\nuME}{\nu_{\Mcal}^\E} % local stabilization
\newcommand{\nuMplus}{\nu_{\Mcal}^+} % stabilization constant
\newcommand{\nuMminus}{\nu_{\Mcal}^-} % stabilization constant
\newcommand{\SED}{S^\E_{D}} % local stabilization
\newcommand{\nuDE}{\nu_{D}^\E} % local stabilization
\newcommand{\nuDplus}{\nu_{\Dcal}^+} % stabilization constant
\newcommand{\nuDminus}{\nu_{\Dcal}^-} % stabilization constant

\newcommand{\Pinabla}{\Pi^{\nabla}_{k+1}}
\newcommand{\PinablaE}{\Pi^{\nabla,\E}_{k+1}}
\newcommand{\Pizbkminusone}{\boldsymbol{\Pi^{0,\E}_{k-1}}}

\newcommand{\PizbkE}{\boldsymbol{\Pi^{0,\E}_{k}}}
\newcommand{\Pizbk}{\boldsymbol{\Pi^{0}_{k}}}
\newcommand{\PizkE}{\Pi^{0,\E}_{k}}

\newcommand{\PizkoE}{\Pi^{0,\E}_{k+1}}
\newcommand{\Pizko}{\Pi^{0}_{k+1}}

\newcommand{\PizbzE}{\boldsymbol{\Pi^{0,\E}_{0}}}
\newcommand{\PizzE}{\Pi^{0,\E}_{0}}

\newcommand{\bchi}{\boldsymbol{\chi}}
\newcommand{\bpsi}{\boldsymbol{\psi}}
\newcommand{\bkappa}{\boldsymbol{\kappa}}
\newcommand{\PizrE}{\Pi^{0,\E}_{r}}
\newcommand{\Pizr}{\Pi^{0}_{r}}
\newcommand{\PizbsE}{\boldsymbol{\Pi^{0,\E}_{s}}}
\newcommand{\Pizbs}{\boldsymbol{\Pi^{0}_{s}}}
\newcommand{\PizbtE}{\boldsymbol{\Pi^{0,\E}_{t}}}
\newcommand{\Pizbt}{\boldsymbol{\Pi^{0}_{t}}}
\newcommand{\sigmahat}{\widehat{\sigma}}
\newcommand{\mr}{m_r}
\newcommand{\ms}{m_s}
\newcommand{\mt}{m_t}

\begin{tiny}
	\author{
		\normalsize{
	}}
\end{tiny}

\title{A virtual element method for the miscible displacement of incompressible fluids in porous media}
\date{}
\author{L. Beir\~{a}o da Veiga\thanks{Dipartimento di Matematica e Applicazioni, Universit\`{a} di Milano-Bicocca, 20125 Milano-Bicocca, Italy (lourenco.beirao@unimib.it, giuseppe.vacca@unimib.it)},\ A. Pichler\thanks{Faculty of Mathematics, University of Vienna, 1090 Vienna, Austria (alex.pichler@univie.ac.at)},\ G.Vacca\footnotemark[1]}

\begin{document}
%%%%%%%%%%%%%%%%%%%%%%%%%%%%%%%%%%%%%
\maketitle

% \medskip\noindent
% \textbf{AMS subject classification}: fff
		
% \medskip\noindent
% \textbf{Keywords}: fff

\begin{abstract}
In the present contribution, we construct a virtual element (VE) discretization for the problem of miscible displacement of one incompressible fluid by another, described by a time-dependent coupled system of nonlinear partial differential equations. Our work represents a first study to investigate the premises of virtual element methods (VEM) for complex fluid flow problems. We combine the VEM discretization with a time stepping scheme and develop a complete theoretical analysis of the method under the assumption of a regular solution. The scheme is then tested both on a regular and on a more realistic test case.%, both inspired from the literature.

\medskip\noindent
\textbf{AMS subject classification}: 65M12, 65M60, 76S05
% 65N12, 65N30, 76R99

\medskip\noindent
\textbf{Keywords}: virtual element methods, miscible fluid flow, porous media, polygonal meshes
\end{abstract}

\section{Introduction}

The virtual element method (VEM) was introduced in \cite{VEMvolley,hitchhikersguideVEM} (see also \cite{VEM3Dbasic}) as a generalization of the finite element method (FEM) that allows to use general polygonal and polyhedral meshes. Since its recent birth in 2013, VEM enjoyed a rapid growth in the mathematics and engineering communities. Among the large numbers of papers in the literature, we here cite only 
\cite{Brezzi-Marini:2012,
absv_VEM_cahnhilliard,
cangianimanzinisutton_VEMconformingandnonconforming,
VEM_DD_basic,
VEMnavierstokes,
VEMnewtonianflows,
VEMmagnetostatic3D,
VEM_nc_SUPG,
VEMtransmissioneigen,
VEMcurved,
TVEM_Helmholtz_num,
VEMcrank,
VEMtopopt}
as representatives, see also the references therein. 

In the realm of diffusion problems, virtual elements have been developed for linear model diffusion-convection-reaction equations in primal and mixed form 
\cite{
VEMvolley,
BBMR_generalsecondorder,
BBMR_generalsecondordermixed,
VEM3Dbasic,
sukumar,
vacca2015virtual,
VEM_SUPG
}. 
It was soon recognized that the flexibility of VEM in terms of meshing could lead to appealing advantages in the presence of complex geometries, such as for discrete fracture network simulation 
\cite{
VEM_dfn1,
VEM_dfn2,
VEM_dfn3} 
and, more in general, in the presence of fractures in porous 3D media
\cite{
VEM_dfn4,
VEM_fumagalli}. 
Nevertheless, although in other frameworks (such as solid mechanics) VEM have indeed proven themselves also on tough nonlinear problems, to the best knowledge of the authors, virtual elements have never been developed and tested for more complex diffusion models. Since VEM have indeed been hardly tested (with very promising outcomes) on linear diffusion problems with complex geometries, as often encountered in geophysical flows, developing a VEM also for more complex (and realistic) flow models becomes a key step towards a competitive methodology for applications. 
%The present paper is a first investigation in such direction.

In the present contribution, we consider the miscible displacement of one incompressible fluid by another in
a reservoir, described by a time-dependent coupled system of nonlinear partial differential equations, that is a basic (but meaningful) model instrumental to applications such as oil recovery and environmental pollution 
\cite{
russell1983finite,
peaceman2000fundamentals,
ewing1983approximation,
chainais2007convergence,
droniou2}.
One must note that, on this and similar models, there already exists a large literature with many competitive schemes, adopting for instance finite elements \cite{
fem1,
fem2,
ewing1983approximation}, 
discontinuous Galerkin methods
\cite{
dg1,
dg2,
dg3}, and
finite volumes 
\cite{chainais2007convergence,fv1,fv2}. 
%and also other more recent polytopal techs \cite{HHO?}. 
The aim of the present paper is to make a first study on the premises of VEM in this framework (by proposing a numerical scheme, giving a theoretical backbone to it, and finally testing it numerically). We believe that the VEM for this kind of problems could have a value due to its strong robustness with respect to the mesh features and, more heuristically, its potential in terms of flexibility in general terms.

From the mathematical viewpoint, the above model yields a nonlinear time-dependent coupled problem for concentration, velocity and pressure, also with potential issues of stability (at the discrete level) due to possible convection-dominated regimes. 
We propose a continuous ($H^1$ conforming) approximation for the concentration variable, thus leading to nodal virtual elements \cite{VEMvolley,BBMR_generalsecondorder}, and an $H({\rm div})$ conforming approximation of the Darcy velocity, thus leading to face virtual elements \cite{Brezzi-Falk-Marini,BBMR_generalsecondordermixed}. For the pressure, we adopt a standard piecewise discontinuous polynomial space. Due to the presence of the non-linear coefficients coupling the two set of equations, we make use of projection operators to approximate such terms and of stabilization factors that are suitably chosen.
We combine the VEM discretization in space with a simple discretization procedure in time, that is a backward Euler approximation that is explicit in the coefficient terms. As a consequence, the system to be solved at each time step is linear and decoupled, leading to a cheap procedure. Extending the proposed scheme to different time discretization procedures would be, on the basis of the work presented here, quite trivial.

After proposing the method, we develop an error analysis under the assumption of a regular solution. Although such regularity conditions are unrealistic in most cases of interest, we believe that the derived results are still critical in order to give a theoretical backbone to the method. They serve the purpose of showing that the method indeed delivers a solution with the potential to yield the correct approximation order whenever this is feasible (given the approximability of the target solution by the discrete space). No time step size condition is needed in the analysis.
Finally, we test the proposed scheme in two different ways. We firstly consider a problem with known regular solution inspired from~\cite{miscibledispl2018hu}, in order to validate the convergence properties of the method also in practice and to test some other practical aspect such as the possibility of having different time step sizes for the two different equations. Then, we consider a more realistic test, taken from~\cite{wang2000approximation}, in order to have a qualitative comparison with the expected benchmark solution from the literature. In this second test, there is also the risk of overshoots and undershoots in the discrete solution due to strong convection. We here deal with this aspect by introducing in our scheme a modification borrowed from~\cite{FCT}, that is recognized in~\cite{smalldiffusion} to be one of the  best choices in practice. 
From the present first theoretical and numerical studies, we believe the VEM is promising and has the possibility to become, after further developments, a competitive scheme for complex flow problems. 

The structure of the paper is the following. In Section \ref{sec:2}, we introduce the continuous problem, in strong and weak form. In Section \ref{sec:3}, we describe the proposed virtual element discretization, in space and time. In Section \ref{sec:4}, we develop the theoretical convergence analysis of the scheme. Finally, in Section \ref{sec:5}, we show the numerical results.

%%%%%%%%%%%%%%%%%%%%%%%%%%%%%%%%%%%%%%%%%%%%%%%
\section{Problem description}\label{sec:2}

We consider the miscible displacement of one incompressible fluid by another in a porous medium. This problem can be formulated in terms of a system of partial differential equations, where a parabolic diffusion-convection-reaction type equation is nonlinearly coupled with an elliptic system, see also~\cite{
peaceman2000fundamentals,
ewing1983approximation,
chainais2007convergence,
droniou2}.

We need to introduce some notation and conventions to be adopted throughout the paper. We denote by $\N$ and $\N_0$ the sets of all natural numbers without and including zero, respectively. Moreover, we employ the standard notation for Sobolev spaces, norms, and seminorms. More precisely, for a given bounded Lipschitz domain $D \subset \R^2$, $k \in \N \cup \{\infty\}$, and $p \in \N$, we define by $W^{k,p}(D)$ the space of all $L^p$ integrable functions over $D$ whose weak derivatives up to order~$k$ are again $L^p$ integrable. Sobolev spaces with fractional order can be defined for instance via interpolation theory~\cite{Triebel}. For $p=2$, we write $H^k(D):=W^{k,2}(D)$, and we use $(\cdot,\cdot)_{k,D}$, $|\cdot|_{k,D}$, and $\lVert \cdot \rVert_{k,D}$, to denote the corresponding inner product, seminorm, and norm, respectively. The standard $L^2$ inner product over $D$ is written as $(\cdot,\cdot)_{0,D}$ with corresponding norm $\lVert \cdot \rVert_{0,\Omega}$. Further, $\mathbb{P}_k(D)$ is the space of polynomials up to order $k$, and $[\mathbb{P}_k(D)]^2$ the corresponding vector valued space. Additionally, $|\cdot|$ is the standard Euclidean norm for scalars and vectors. 
Finally, throughout the paper, $\eta$ denotes a generic constant, possibly varying from one occurrence to the other, but independent of the mesh size and, apart from Theorem~\ref{thm:Cn_cn}, also independent of the variables.

\subsection{Continuous Problem}
Let $\Omega \subset \R^2$ be a polygonal bounded, convex domain, describing a reservoir of unit thickness. Given a time interval $J:=[0,T]$, for $T>0$, we are interested in finding $\bu=\bu(\bx,t)$, representing the Darcy velocity (volume of fluid flowing cross a unit across-section per unit time), the pressure $p=p(\bx,t)$ in the fluid mixture, and the concentration $c=c(\bx,t)$ of one of the fluids (amount of the fluid per unit volume in the fluid mixture), with $(\bx,t) \in \Omega_T:=\Omega \times J$, such that
\begin{equation} \label{eq:model_problem}
\left\{
\begin{alignedat}{2}
\phi \, \frac{\partial c}{\partial t} + \bu \cdot \nabla c - \div (D(\bu) \nabla c) &= \qplus(\chat-c) \\
\div \, \bu &= G \\
\bu &= -a(c) (\nabla p - \bgamma(c)),
\end{alignedat}
\right.
\end{equation}
where $\phi=\phi(\bx)$ is the porosity of the medium, $\qplus=\qplus(\bx,t)$ and $\qminus=\qminus(\bx,t)$ are the (non negative) injection and production source terms, respectively, $\chat=\chat(\bx,t)$ is the concentration of the injected fluid, and
\begin{equation} \label{eq:def_G}
G:=\qplus-\qminus.
\end{equation} 
Moreover, $D(\bu) \in \R^{2\times 2}$ is the diffusion tensor given by 
\begin{equation} \label{eq:def_D}
D(\bu):=\phi \left[ d_m I + |\bu|(d_\ell E(\bu) + d_t E^\perp(\bu)) \right],
\end{equation}
with matrices 
\begin{equation*} \label{eq:def_E(u)}
E(\bu):=\left( \frac{\bu_i \bu_j}{|\bu|^2} \right)_{i,j=1,2}=\frac{\bu \bu^T}{|\bu|^2}, \quad E^\perp(\bu):=I-E(\bu),
\end{equation*}
and molecular diffusion coefficient $d_m$, longitudinal dispersion coefficient $d_\ell$, and transversal dispersion coefficient $d_t$. Further, $\boldsymbol{\gamma}(c)$ in~\eqref{eq:model_problem} describes the force density due to gravity (typically written as $\boldsymbol{\gamma}(c)=\gamma_0(c)\boldsymbol{\rho}$ with $\gamma_0(c)$ being the density of the fluid and $\boldsymbol{\rho}$ the gravitational acceleration vector), and $a(c)=a(c,\bx)$ is the scalar valued function given by
\begin{equation*}
a(c):=\frac{k}{\mu(c)},
\end{equation*}
where $k=k(\bx)$ represents the permeability of the porous rock, and $\mu(c)$ is the viscosity of the fluid mixture, which can be modeled by
\begin{equation*}
\mu(c)=\mu(0) \left( 1+\left(M^{\frac{1}{4}}-1\right)c \right)^{-4}, \quad \text{ in } [0,1],
\end{equation*}
with mobility ratio $M:=\frac{\mu(0)}{\mu(1)}$. Note that $\mu$ can be set to $\mu(0)$ for $c<0$, and to $\mu(1)$ for $c>1$. We also highlight that, in the literature, $k$ is sometimes assumed to be a tensor. The following analysis can be straightforwardly generalized to that case.

Assuming impermeability of $\partial \Omega$, the system~\eqref{eq:model_problem} is closed by requiring \textit{no-flow boundary conditions} of the form

\begin{equation} \label{eq:bdry_cond}
\left\{
\begin{alignedat}{2}
\bu \cdot \bn&=0 \quad \text{ on } \partial \Omega \times J \\
D(\bu) \nabla c \cdot \bn&=0 \quad \text{ on } \partial \Omega \times J,
\end{alignedat}
\right.
\end{equation}
and initial condition 
\begin{equation} \label{eq:initial_cond}
c(\bx,0) = c_0(\bx) \quad \text{ in } \Omega,
\end{equation}
where $0\le c_0(\bx) \le 1$ is an initial concentration.
  
By use of the divergence theorem, the boundary conditions~\eqref{eq:bdry_cond} directly imply the following compatibility condition for $\qplus$ and $\qminus$:
\begin{equation*}
\int_{\Omega} \qplus(\bx,t) \, \dx = \int_{\Omega} \qminus(\bx,t) \, \dx,
\end{equation*}
for every $t \in J$.

We highlight that, in the forthcoming theoretical analysis, we will always assume sufficient regularity of the exact solution and the involved functions, such as $\qplus$, $\qminus$, $\chat$, \textit{et cetera}, as better motivated in the corresponding section. Moreover, we will make use of the following assumptions.

First of all, we suppose that the functions $a$ and $\phi$ are positive and uniformly bounded from below and above, i.e. there exist positive constants $a_\ast$, $a^\ast$, $\phi_\ast$, and $\phi^\ast$, such that
\begin{equation} \label{ass:a_phi}
a_\ast \le a(z,\bx) \le a^\ast, \qquad 
\phi_\ast \le \phi(\bx) \le \phi^\ast,
\end{equation}
for all $\bx \in \Omega$ and $z=z(t)$. For the sake of readability, we define
\begin{equation*}
A(z)(\bx):=a^{-1}(z,\bx).
\end{equation*}
Additionally, we will make use of the following relation of the diffusion and dispersion coefficients, which was observed in laboratory experiments:
\begin{equation} \label{ass:dl_dt_dm}
0 < d_m \le d_t \le d_\ell.
\end{equation}

Finally, we recall that the source terms $\qplus$ and $\qminus$ are, as usual, assumed to be non-negative functions.

Existence of weak solutions to this model problem was shown in~\cite{feng1995existence} for $\bgamma(c)=0$. An extension of this result to 3D spatial domains, including the presence of $\bgamma(c)$ and various boundary conditions was discussed in~\cite{chen1999mathematical}.

\subsection{Weak formulation of the continuous problem}

Here, we fix the basic notation and the functional setting.

To this purpose, given $\Omega$ as above, we first introduce the Sobolev space 
\begin{equation*}
\HdivOmega:=\{ \bv \in [L^2(\Omega)]^2: \, \div \bv \in L^2(\Omega) \}.
\end{equation*}
Then, we define the velocity space $\boldsymbol{V}$, the pressure space $Q$, and the concentration space $Z$ by
\begin{equation} \label{eq:spaces_V_Q_Z}
\begin{split}
\boldsymbol{V}&:=\{ \bv \in \HdivOmega: \, \bv \cdot \bn = 0 \, \text{ on } \partial \Omega \} \\
Q&:=L^2_0(\Omega):=\{ \varphi \in L^2(\Omega): \, (\varphi,1)_{0,\Omega}=0 \} \\
Z&:=H^1(\Omega),
\end{split}
\end{equation}
respectively. These spaces are endowed, respectively, with the following norms: 
\begin{equation*}
\begin{split}
\lVert \bu \rVert_{\boldsymbol{V}}^2:=\lVert \bu \rVert_{0,\Omega}^2 + \lVert \div \bu \rVert_{0,\Omega}^2, \qquad
\lVert q \rVert_Q^2:=\lVert q \rVert_{0,\Omega}^2, \qquad
\lVert z \rVert_Z^2:=\lVert z \rVert_{1,\Omega}^2:=\lVert \nabla z \rVert^2_{0,\Omega} + \lVert z \rVert^2_{0,\Omega}.
\end{split}
\end{equation*}
Note that $\div \boldsymbol{V}=Q$.

As usual in the framework of parabolic problems, we use the notation
\begin{equation} \label{ansatz:cont_fcts}
\bu(t)(x):=\bu(x,t), \qquad p(t)(x):=p(x,t), \qquad c(t)(x):=c(x,t).
\end{equation}
For $0 \le a \le b$, we further introduce
\begin{equation*}
\lVert \bv \rVert_{L^2(a,b;\boldsymbol{V})}:=\left(\int_{a}^{b} \lVert \bv(t) \rVert^2_{\boldsymbol{V}} \, \dx\right)^{\frac{1}{2}}, \quad
\lVert \bv \rVert_{L^\infty(a,b;\boldsymbol{V})}:=\underset{t\in [a,b]}{\text{ess sup}} \lVert \bv(t) \rVert_{\boldsymbol{V}};
\end{equation*}
analogously for $p$ and $c$.

Having this, the continuous problem reads as follows: find 
$c \in L^2(0,T;Z) \cap C^0([0, T]; L^2(\Omega))$, 
%with $\frac{\partial c}{\partial t} \in L^2(0,T;Z')$, 
$\bu \in L^2(0,T;\boldsymbol{V})$, and $p \in L^2(0,T;Q)$, such that
\begin{comment}
\begin{equation*}
\begin{split}
\{\bu,p,c\}: J \to \boldsymbol{V} \times Q \times Z, \quad t \mapsto \{\bu(t),p(t),c(t)\},
\end{split}
\end{equation*}
such that
\end{comment}
\begin{equation} \label{eq:cont_form}
\left\{
\begin{alignedat}{2}
\Mcal \left(\frac{\partial c(t)}{\partial t},z \right) + \left( \bu(t) \cdot \nabla c(t),z \right)_{0,\Omega}  + \Dcal(\bu(t);c(t),z) &= 
\left( \qplus(\chat-c)(t),z \right)_{0,\Omega} \\
\Acal(c(t);\bu(t),\bv) + B(\bv,p(t)) &= (\bgamma(c(t)),\bv)_{0,\Omega}  \\
B(\bu(t),q) &= - \left( G(t), q \right)_{0,\Omega}
\end{alignedat}
\right.
\end{equation}
for all $\bv \in \boldsymbol{V}$, $q \in Q$, and $z \in Z$, for almost all $t \in J$ and with initial condition $c(0)=c_0$, where
\begin{alignat}{2} 
\Mcal(c,z)&:=\left(\phi \, c,z\right)_{0,\Omega}, &\qquad  \Dcal(\bu;c,z)&:=\left( D(\bu) \nabla c,\nabla z \right)_{0,\Omega}, \label{eq:def_A_B_D} \\
\Acal(c;\bu,\bv)&:=\left( A(c) \bu,\bv \right)_{0,\Omega} , & B(\bv,q)&:=-\left( \div \bv,q \right)_{0,\Omega}. \nonumber
\end{alignat}
Note that $c \in L^2(0,T;Z) \cap C^0([0, T]; L^2(\Omega))$ implies $\frac{\partial c}{\partial t} \in L^2(0,T;Z')$, see e.g.~\cite[Thm. 11.1.1]{quarteroni2008numerical}.
%{\gv{It holds that $c \in L^2(0,T;Z) \cap C^0([0, T]; L^2(\Omega))$ that implies $\frac{\partial c}{\partial t} \in L^2(0,T;Z')$. For a reference see the book Quarteroni-Valli, chapter 11, thm 11.1.1 }}

For the sake of readability, we suppressed $(t)$ in~\eqref{eq:def_A_B_D}. From now on, we will use the convention that by writing $\bu$, we mean in fact $\bu(t)$; similarly for the other functions depending on space and time. In general it will be clear from the context whether $\bu$ represents $\bu(t)$ for a fixed $t \in J$, i.e. as a function of space only, or for varying $\bx$ and $t$, as a function of both space and time.

Moreover, we will use the following alternative form for the concentration equation:
\begin{equation} \label{eq:cont_form_altern}
\begin{split}
\Mcal \left(\frac{\partial c}{\partial t},z \right)  +\Theta(\bu,c;z)+ \Dcal(\bu;c,z) 
= \left( \qplus \, \chat,z \right)_{0,\Omega},
\end{split}
\end{equation}
where
\begin{equation*}
\Theta(\bu,c;z):=\frac{1}{2} \bigg[\left( \bu \cdot \nabla c,z \right)_{0,\Omega}  + ((\qplus+\qminus) \, c, z )_{0,\Omega}  - \left( \bu, c \, \nabla z \right)_{0,\Omega} \bigg].
\end{equation*}
This version is obtained from the original one in~\eqref{eq:cont_form} by rewriting the convective term as
\begin{equation*}
\left( \bu \cdot \nabla c,z \right)_{0,\Omega} 
=\frac{1}{2} \big[\left( \bu \cdot \nabla c,z \right)_{0,\Omega}  - (G, c \, z)_{0,\Omega}  - \left( \bu, c \, \nabla z \right)_{0,\Omega}  \big],
\end{equation*}
where we first integrated by parts, then employed the fact that $\nabla \cdot \bu=G$, together with the definition of $G$ in~\eqref{eq:def_G}, and afterwards combined this term with $(\qplus\, c,z)_{0,\Omega}$ from the right hand side of~\eqref{eq:cont_form}. This representation was inspired by the theory of VEM for general elliptic problems~\cite{cangianimanzinisutton_VEMconformingandnonconforming} and helps to ensure that properties of the continuous bilinear will be preserved after discretization. 

In the rest of this section, we summarize some properties of the forms $\Mcal(\cdot,\cdot)$, $\Acal(\cdot,\cdot,\cdot)$ and $\Dcal(\cdot;\cdot,\cdot)$, all defined in~\eqref{eq:def_A_B_D}, which will be needed later on. %Let $t\in J$ be fixed for the moment.

To start with, for $\Mcal(\cdot,\cdot)$, it directly holds with the Cauchy-Schwarz inequality and~\eqref{ass:a_phi}
\begin{equation*}
\Mcal(c,z) \le \phi^\ast \lVert c \rVert_{0,\Omega} \lVert z \rVert_{0,\Omega}, \qquad 
\Mcal(z,z) \ge \phi_\ast \lVert z \rVert_{0,\Omega}^2,
\end{equation*}
for all $c,z\in Z$.

Concerning $\Acal(\cdot;\cdot,\cdot)$, again employing~\eqref{ass:a_phi}, for all $c\in L^{\infty}(\Omega)$ and $\bu,\bv \in [L^2(\Omega)]^2$, we have
\begin{equation*}
\Acal(c;\bu,\bv) \le \frac{1}{a_\ast} \lVert \bu \rVert_{0,\Omega} \lVert \bv \rVert_{0,\Omega}.
\end{equation*}
Further, if $c \in L^2(\Omega)$, $\bu \in [L^\infty(\Omega)]^2$ and $\bv \in [L^2(\Omega)]^2$, it holds true that
\begin{equation*}
\Acal(c;\bu,\bv) \le \lVert A(c) \rVert_{0,\Omega} \lVert \bu \rVert_{\infty,\Omega} \lVert \bv \rVert_{0,\Omega}.
\end{equation*}
We also have the coercivity bound
\begin{equation*} \label{eq:coercivity_A}
\Acal(c;\bv,\bv) \ge \frac{1}{a^\ast} \lVert \bv \rVert_{0,\Omega}^2
\end{equation*}
for all $c \in L^\infty(\Omega)$ and $\bv \in [L^2(\Omega)]^2$, from which, after defining the kernel
\begin{equation} \label{eq:cont_kernel}
\Kcal:=\{ \bv \in \boldsymbol{V}: \, B(\bv,q)=0 \quad \forall q \in Q \},
\end{equation}
coercivity of $\Acal(c;\cdot,\cdot)$ on $\Kcal$ in the norm $\lVert \cdot \rVert_{V}$ follows.% for all $c \in \ap{L^\infty(\Omega)}$.

Regarding $\Dcal(\cdot;\cdot,\cdot)$, the following continuity properties can be shown.
Firstly, for all $\bu \in [L^\infty(\Omega)]^2$ and $c,z \in H^1(\Omega)$, we have
\begin{equation} \label{eq:cont_D_1}
\begin{split}
\Dcal(\bu;c,z) \le \phi^\ast \left[ d_m + \lVert \bu \rVert_{\infty,\Omega} (d_\ell+d_t)\right] \lVert \nabla c \rVert_{0,\Omega} \lVert \nabla z \rVert_{0,\Omega},
\end{split}
\end{equation}
which follows directly from the Cauchy-Schwarz inequality, the definition of $D(\bu)$ in~\eqref{eq:def_D}, and the fact that $|E(\bu) \bv| \le |\bv|$ and $|E^\perp(\bu) \bv| \le |\bv|$ for all $\bv \in \R^2$. Moreover, for all $\bu \in [L^2(\Omega)]^2$ and $c,z \in H^1(\Omega)$ with $\nabla c \in L^\infty (\Omega)$, we have the bound 
\begin{equation} \label{eq:cont_D_2}
\Dcal(\bu;c,z) \le \lVert D(\bu) \rVert_{0,\Omega} \lVert \nabla c \rVert_{\infty,\Omega} \rVert \nabla z \rVert_{0,\Omega} 
\le \eta_{\mathcal{D}} (1+\lVert \bu \rVert_{0,\Omega}) \lVert \nabla c \rVert_{\infty,\Omega} \rVert \nabla z \rVert_{0,\Omega} ,
\end{equation}
with matrix norm $\lVert D(\bu) \rVert_{0,\Omega}:=\left(\sum_{i,j=1}^{2} \lVert D_{i,j}(\bu) \rVert^2_{0,\Omega} \right)^{\frac{1}{2}}$, and some positive constant $\eta_{\mathcal{D}}$ depending only on $d_m$, $d_{\ell}$, and $d_t$. In addition, coercivity of $\Dcal(\bu;\cdot,\cdot)$ for all $\bu \in [L^\infty(\Omega)]^2$, with respect to $\lVert \cdot \rVert_{0,\Omega}$, follows from
\begin{equation} \label{eq:Du_mu_mu}
\begin{split}
(D(\bu) \, \bmu ,\bmu)_{0,\Omega} &= (\phi \, d_m \, \bmu,\bmu)_{0,\Omega} + (\phi \, |\bu| \, (d_\ell E(\bu) + d_t E^\perp (\bu)) \, \bmu,\bmu)_{0,\Omega} \\
&\ge \phi_\ast \, d_m \, \lVert \bmu \rVert_{0,\Omega}^2 + (\phi \,  |\bu| (d_\ell-d_t) E(\bu) \bmu,\bmu )_{0,\Omega} + (\phi \, |\bu| d_t \,  \bmu,\bmu )_{0,\Omega} \\
&\ge \phi_\ast \left(d_m \, \lVert \bmu \rVert_{0,\Omega}^2 + d_t \, \lVert |\bu|^{\frac{1}{2}} \bmu \rVert_{0,\Omega}^2 \right)
\end{split}
\end{equation}
for all $\bmu \in [L^2(\Omega)]^2$, where we also employed~\eqref{ass:a_phi} and \eqref{ass:dl_dt_dm}.

%\todo{ALEXANDER: Remark about well-posedness? YES, IF YOU FIND A DIRECT REFERENCE IN THE LITERATURE (NOT IF WE NEED TO PROVE IT). CAN YOU CHECK?}

\section{The virtual element method}\label{sec:3}

In this section, we derive a virtual element formulation for the model problem~\eqref{eq:cont_form}. To this purpose, we firstly fix the concept of polygonal decompositions of $\Omega$ in Section~\ref{subsec:polygonal_decomp}, and then, we introduce a set of discrete spaces, discrete bilinear forms, and projectors in Section~\ref{subsec:discr_spaces_proj}. Having these ingredients, we state a semidiscrete formulation which is continuous in time and discrete in space in Section~\ref{subsec:semidiscr_form}. The fully discrete formulation is the subject of Section~\ref{subsec:fully_discr}.

\subsection{Polygonal decompositions} \label{subsec:polygonal_decomp}

Let $\taun$ be a discretization of $\Omega$ into polygons $\E$. We denote by $\En$ the set of all edges of $\taun$, and, for a given element $\E \in \taun$, by $\EE$ the set of edges belonging to $\E$. Furthermore, $\Ne$ is the number of edges of $\E$, $\hE$ is the diameter of $\E$, and $h:=\max_{\E \in \taun} \hE$. For a given edge $\e \in \En$, we write $\he$ for its length. Having this, we make the following assumptions on $\taun$: there exists $\rho_0>0$ such that, for all $h > 0$ and for all $\E \in\taun$,

\begin{itemize}
\item[(\textbf{D1})] $\E$ is star-shaped with respect to a ball of radius $\rho \ge \rho_0 \hE$;
\item[(\textbf{D2})] $\he \ge \rho_0 \hE$ for all $\e \in \EE$.
\end{itemize}
Note that these two assumptions imply that the number of edges of each element is uniformly bounded. Additionally, we will require quasi-uniformity: 
\begin{itemize}
\item[(\textbf{D3})] for all $h>0$ and for all $\E \in \taun$, it holds $h_K \ge \rho_1 h$, for some positive uniform constant $\rho_1$.
\end{itemize}

Given $\taun$, we define, for all $s>0$, the broken Sobolev spaces on $\taun$ as
\begin{equation*} \label{broken_Sobolev_space}
H^{s}(\taun) := \{v \in L^2(\Omega) \mid v_{|_\E} \in H^s(\E) \ \forall \E \in \taun\},
\end{equation*}
together with the corresponding broken seminorms and norms
\begin{equation} \label{broken_Sobolev_norm}
\vert v \vert^2_{s,\taun} := \sum_{\E \in \taun} \vert v \vert_{s,\E}^2, \quad \quad \quad \Vert v \Vert^2_{s,\taun} := \sum_{\E \in \taun} \Vert v \Vert_{s,\E}^2.
\end{equation}

\begin{remark}
Both assumptions ({\bf D1}) and ({\bf D2}) are standard in the virtual element literature. While condition ({\bf D1}) is quite critical in the following analysis, assumption ({\bf D2}) could be possiby avoided by following steps similar to \cite{beiraolovadinarusso_stabilityVEM,BrennerGuanSung_someestimatesVEM}, at the expense of making the proofs even more lenghty and technical. Finally, assumption ({\bf D3}) (that can be found also in many FEM papers on the same subject) is only needed to prove bound \eqref{eq:bound_Pi_nabla_Pc} below.
\end{remark}

\subsection{Discrete spaces and projectors} \label{subsec:discr_spaces_proj}

Here, we introduce the local discrete VE spaces corresponding to $\boldsymbol{V}$, $Q$ and $Z$ in~\eqref{eq:spaces_V_Q_Z}, a set of local projectors mapping from these VE spaces into spaces made of polynomials, and finally, the related global counterparts.

\subsubsection{Local discrete spaces}

Let $\E \in \taun$ and let $k \in \N_0$ be a given \textit{degree of accuracy}. Then, the local velocity and pressure VE spaces are defined by
\begin{equation} \label{eq:def_local_spaces}
\begin{split}
\VhE&:=\{ \bv \in \HdivE \cap \HrotE: \, {\bv \cdot \bn}_{|_e} \in \mathbb{P}_k(e) \, \forall e \in \EE, \\
&\qquad \div \bv \in \mathbb{P}_{k}(\E), \, \rot \bv \in \mathbb{P}_{k-1}(\E) \}\\
\QhE&:=\{ q \in L^2(\E): \, q \in \mathbb{P}_{k}(\E) \}.
\end{split}
\end{equation}
These spaces are coupled with the {\it preliminary} local concentration space
\begin{equation} \label{eq:def_local_spaces_conc}
\ZhEtilde:=\{ z \in H^1(\E): \, z_{|_{\partial \E}} \in C^0(\partial \E), \, z_{|_e} \in \mathbb{P}_{k+1}(e) \, \forall e \in \EE, \, \Delta z \in \mathbb{P}_{k-1}(\E) \}.
\end{equation}
%Their local dimensions are given by 
%\begin{equation*}
%\begin{split}
%\NVh&:=\dim \VhE
%=\Ne (k+1) + (k+1)^2 -1 \\
%\NQh&:=\dim \QhE=\frac{(k+1)(k+2)}{2}\\
%\NZhtilde&:=\dim \ZhEtilde = \Ne (k+1)+\frac{k(k+1)}{2}.
%\end{split}
%\end{equation*}
Moreover, it is important to observe that $[\mathbb{P}_k(\E)]^2 \subset \VhE$ and $\mathbb{P}_{k+1}(\E) \subseteq \ZhEtilde$. Associated sets of local degrees of freedom are given as follows:
\begin{itemize}
\item for $\VhE$, a set of degrees of freedom $\{\dofVhE_j\}_{j=1}^{\NVhE}$ is defined by
\begin{equation} \label{eq:dofs_VhK}
\begin{split}
1.& \quad \frac{1}{|e|} \int_\e \bv \cdot \bn \, p_k \, \ds \qquad \forall p_k \in \mathbb{P}_k(e) \quad \forall e \in \EE \\
2.& \quad \frac{1}{|\E|^{\frac{1}{2}}} \int_\E (\div \bv) \, p_{k} \, \dx \qquad \forall p_{k} \in \mathbb{P}_{k}(\E)/\R \\
3.& \quad \frac{1}{|\E|} \int_\E \bv \cdot \boldsymbol{x}^\perp \, p_{k-1} \, \dx \qquad \forall p_{k-1} \in \mathbb{P}_{k-1}(\E),
\end{split}
\end{equation}
with $\bx^\perp:=(\bx_2,-\bx_1)^T$, where we assume the coordinates to be centered at the barycenter of the element;
\item for $\QhE$, we consider $\{\dofQhE_j\}_{j=1}^{\NQhE}$ with 
\begin{equation} \label{eq:dofs_QhK}
\frac{1}{|\E|} \int_\E q \, p_{k} \, \dx \qquad \forall p_{k} \in \mathbb{P}_{k}(\E);
\end{equation}
\item for $\ZhEtilde$, we take $\{\dofZhtildeE_j\}_{j=1}^{\NZhtildeE}$ with
\begin{equation} \label{eq:dofs_ZhK}
\begin{split}
1.& \quad \text{ pointwise values at the vertices: } \bv(z) \\
2.& \quad \text{ on each edge $e\in \EE$, the values of $z$ at the $k$ internal Gau\ss{}-Lobatto points} \\
3.& \quad \frac{1}{|\E|} \int_\E  z \, q_{k-1} \, \dx \qquad \forall q_{k-1} \in \mathbb{P}_{k-1}(\E).
\end{split}
\end{equation}
\end{itemize}
In all three cases, unisolvency is provided. More precisely, for $\VhE$, this was proven in e.g.~\cite{HdivHcurlVEM}, for $\QhE$ it is immediate, and for $\ZhEtilde$, see e.g.~\cite{VEMvolley}. 

We also highlight that $\VhE$ endowed with~\eqref{eq:dofs_VhK} mimics the Raviart-Thomas element, but in fact those two elements only coincide in the special case of triangles and $k=0$. An analogous result is true for $\ZhEtilde$, when compared to finite elements.

\begin{remark}
We note that, for $k=0$, one obtains the lowest order local VE spaces. More precisely, in this case, the velocity space $\VhE$ consists of all rotation free vector fields with constant divergence and edgewise constant normal traces, the pressure space $\QhE$ only contains the constant functions, and the concentration space $\ZhEtilde$ is made of all harmonic functions that are linear on each edge. This motivates the choice of the present polynomial degrees for the spaces. However, in general, it is also possible to  choose a degree of accuracy $k_1$ for $\VhE$ and $\QhE$, and another strictly positive one $k_2$ for $\ZhEtilde$; see e.g.~\cite{ewing1983approximation} for FEM. The following analysis can be extended easily to such more general case just by keeping track of the different polynomial degrees.
\end{remark}

\begin{remark}
In order to really have a set of degrees of freedom in the computer code, one clearly needs to choose a basis for the polynomial test spaces appearing in \eqref{eq:dofs_VhK} and \eqref{eq:dofs_ZhK}. We here assume to take the classical choice, that is any monomial basis $\{ m_1,m_2,..,m_\ell \}$ of the polynomial space satisfying $\| m_i \|_{L^\infty} \simeq 1$, $i=1,2,..,\ell$, where the $L^\infty$ norm has to be taken over the corresponding edge or bulk.
\end{remark}

\subsubsection{Local projections} \label{subsec:proj}

For the construction of the method, we will need some tools to deal with VE functions due to the lack of their explicit knowledge in closed form. These tools will be provided in the form of local operators mapping VE functions onto polynomials. To this purpose, following~\cite{VEMvolley,hitchhikersguideVEM}, we introduce the subsequent projectors.

\medskip\noindent
The projector $\PizbkE:\, [L^2(\E)]^2 \to [\mathbb{P}_{k}(\E)]^2$ is defined as the $L^2$ projector onto vector valued polynomials of degree at most $k$ in each component: Given $\boldsymbol{f} \in [L^2(\Omega)]^2$,
\begin{equation} \label{eq:L2_proj}
(\PizbkE \boldsymbol{f}, \boldsymbol{p}_{k})_{0,\E}=( \boldsymbol{f}, \boldsymbol{p}_{k})_{0,\E} \quad \forall \boldsymbol{p}_{k} \in [\mathbb{P}_{k}(\E)]^2.
\end{equation}
It can be shown, see~\cite{BBMR_generalsecondorder}, that this operator is computable for functions in $\VhE$ only by knowing their values at the degrees of freedom~\eqref{eq:dofs_VhK}. 
Moreover, one has computability also for functions of the form $\nabla \zh$ with $\zh \in \ZhEtilde$. This can be seen by using integration by parts:
\begin{equation*}
\int_\E (\PizbkE \nabla \zh) \cdot \boldsymbol{p}_{k} \, \ds = \int_\E \nabla \zh \cdot \boldsymbol{p}_{k} \, \ds
= -\int_\E \zh \, \underbrace{\div \boldsymbol{p}_{k}}_{\in \mathbb{P}_{k-1}(\E)} \, \ds +  \int_{\partial \E} \zh \, \boldsymbol{p}_{k} \cdot \bn \, \ds,
\end{equation*}
for all $\boldsymbol{p}_{k} \in [\mathbb{P}_{k}(\E)]^2$, where the right hand side is computable by means of~\eqref{eq:dofs_ZhK}.

\medskip\noindent
The projector $\PinablaE: \, H^1(\E) \to \mathbb{P}_{k+1}(\E)$ is given, for every $z \in H^1(\E)$, by
\begin{equation*}
\left\{
\begin{split}
(\nabla \PinablaE z, \nabla p_k)_{0,\E} &= (\nabla z, \nabla p_k)_{0,\E} \quad \forall p_{k+1} \in \mathbb{P}_{k+1}(\E) \\
\frac{1}{|\partial \E|} \int_{\partial \E} \PinablaE z \, \ds & = \frac{1}{|\partial \E|} \int_{\partial \E} z \, \ds,
% P^{\partial \E}(\PinablaE \zh) &= P^{\partial \E}(\zh),
\end{split}
\right.
\end{equation*}
where the second identity 
% containing $P^{\partial \E}(\zh):=\frac{1}{|\partial \E|} \int_{\partial \E} \zh \, \ds$
is needed to fix the constants. Computability of this mapping for functions in $\ZhEtilde$ was shown in~\cite{VEMvolley,hitchhikersguideVEM}. 

% --------------------------------------------------------
\subsubsection{Discrete space for concentrations} \label{subsec:discr_conc_space}
% -------------------------------------------------------- 
The space introduced in \eqref{eq:def_local_spaces_conc} was a preliminary space, useful to introduce the main idea of the construction. Nevertheless, we will here make use of a more advanced space for the discrete concentration variable.
Indeed, one can use the operator $\PinablaE$ to pinpoint the local enhanced space
\begin{equation*}
\begin{split}
\ZhE:=\{z \in H^1(\E)&: \, z_{|_{\partial \E}} \in C^0(\partial \E), \, z_{|_e} \in \mathbb{P}_{k+1}(e) \, \forall e \in \EE, \, \Delta z \in \mathbb{P}_{k+1}(\E), \\
&\int_\E z \, p_k \, \dx = \int_\E (\PinablaE z) \, p_k \, \dx \quad \forall p_k \in \mathbb{P}_{k+1}/\mathbb{P}_{k-1}(K)  \},
\end{split}
\end{equation*}
where $\mathbb{P}_{k+1}/\mathbb{P}_{k-1}(K)$ is the space of polynomials in $\mathbb{P}_{k+1}(\E)$ which are $L^2(\E)$ orthogonal to $\mathbb{P}_{k-1}(\E)$. It can be shown that the space $\ZhE$ has the same dimension and the same degrees of freedom~\eqref{eq:dofs_ZhK} as $\ZhEtilde$, see~\cite{equivalentprojectorsforVEM, bbmr_VEM_generalsecondorderelliptic}.
The advantage of the space $\ZhE$, when compared to $\ZhEtilde$, is that \emph{also} the $L^2$ projector $\PizkoE:\, L^2(\E) \to \mathbb{P}_{k+1}(\E)$ onto polynomials of degree at most $k+1$, defined analogously to~\eqref{eq:L2_proj}, is computable \cite{hitchhikersguideVEM}

Finally, we state the following approximation result for the three projectors above~\cite[Lemma 5.1]{BBMR_generalsecondorder}:
\begin{lem} \label{lem:approx_properties}
Given $\E \in \taun$, let $\psi$ and $\boldsymbol{\psi}$ be sufficiently smooth scalar and vector valued functions, respectively. Then, it holds, for all $k \in \N_0$,
\begin{equation*} \label{eq:approx_properties}
\begin{split}
\lVert \psi - \PizkE \psi \rVert_{\ell,\E} &\le \zeta \, \hE^{s-\ell} \, |\psi|_{s,\E}, \quad 0 \le \ell \le s \le k+1 \\
\lVert \boldsymbol{\psi} - \PizbkE \boldsymbol{\psi} \rVert_{\ell,\E} &\le \zeta \, \hE^{s-\ell} \, |\boldsymbol{\psi}|_{s,\E}, \quad 0 \le \ell \le s \le k+1 \\
\lVert \psi - \Pi^{\nabla,\E}_k \psi \rVert_{\ell,\E} &\le \zeta \, \hE^{s-\ell} \, |\psi|_{s,\E}, \quad 0 \le \ell \le s \le k+1, \, s \ge 1,
\end{split}
\end{equation*}
where $\zeta>0$ only depends on the shape-regularity parameter $\rho_0$ in assumption (\textbf{D1}), and $k$.
\end{lem}

\subsubsection{Global discrete spaces and projectors}

The global discrete spaces are defined via their local counterparts:
\begin{equation*}
\begin{split}
\Vh&:=\{ \bv \in \boldsymbol{V}: \, {\bv}_{|_\E} \in \VhE \, \forall \E \in \taun \} \\
\Qh&:=\{ q \in Q: \, q_{|_\E} \in \QhE \, \forall \E \in \taun \} \\
\Zh&:=\{ z \in Z:  \, z_{|_\E} \in \ZhE \, \forall \E \in \taun \} 
\end{split}
\end{equation*}
with the obvious sets of global degrees of freedom.

In addition to the broken Sobolev norm~\eqref{broken_Sobolev_norm}, we introduce, for all $\buh \in \Vh$,
\begin{equation*}
\lVert \buh \rVert_{\Vh}^2
:=\sum_{\E \in \taun} \lVert \buh \rVert_{V,\E}^2
:=\sum_{\E \in \taun} \left[\lVert \buh \rVert_{0,\E}^2 + \lVert \div \buh \rVert_{0,\E}^2 \right].
\end{equation*}

Moreover, we will denote by $\Pizbk$, $\Pinabla$ and $\Pizko$, the global projectors which are defined elementwise as the corresponding local ones in Section~\ref{subsec:proj} and~\ref{subsec:discr_conc_space}.

The sets of global degrees of freedom $\{\dofVh_j\}_{j=1}^{\NVh}$, $\{\dofQh_j\}_{j=1}^{\NQh}$, and $\{\dofZh_j\}_{j=1}^{\NZh}$ are obtained by coupling the local counterparts given in~\eqref{eq:dofs_VhK},~\eqref{eq:dofs_QhK}, and~\eqref{eq:dofs_ZhK}, respectively.

\subsection{Semidiscrete formulation} \label{subsec:semidiscr_form}

Our aim in this section is to find a semidiscrete formulation for~\eqref{eq:cont_form} which is continuous in time and discrete in space. To this purpose, we employ the same notation for the numerical approximants $\buh$, $\ph$, and $\ch$, as in~\eqref{ansatz:cont_fcts} for $\bu$, $p$, and $c$, namely
\begin{equation*}
\buh(t)(x):=\buh(x,t), \qquad p_h(t)(x):=p_h(x,t), \qquad \ch(t)(x):=\ch(x,t),
\end{equation*}
where the dependence on $(t)$ will be again suppressed in the sequel.

A semidiscrete variational formulation for~\eqref{eq:cont_form} can then be written in an abstract way as follows: 
for almost every $t \in J$, find $\buh \in \Vh$, $\ph \in \Qh$, and $\ch \in \Zh$, such that
\begin{equation} \label{eq:semidiscr_var_form}
\left\{
\begin{alignedat}{2}
\Mcalh \left(\frac{\partial \ch}{\partial t},\zh \right) + \Thetah(\buh,\ch;\zh) + \Dcalh(\buh;\ch,\zh) &= \left( \qplus \,\chat,\zh \right)_h \\
\Acalh(\ch;\buh,\bvh) + B(\bvh,p_h) &= (\bgamma(\ch),\bvh)_h \\
B(\buh,\qh) &= - \left( G, \qh \right)_{0,\Omega}
\end{alignedat}
\right.
\end{equation}
%(in principle: $\buh(t)$, $p_h(t)$, $\ch(t)$; for readability: '(t)' suppressed) 
for all $\bvh \in \Vh$, $\qh \in \Qh$, and $\zh \in \Zh$, and the initial condition
\begin{equation*}
\ch(0)=c_{0,h} := I_h c_0
\end{equation*}
is satisfied, where $I_h c_0$ is the VEM interpolant of $c_0$ in $Z_h$, and where the involved forms and terms in~\eqref{eq:semidiscr_var_form} are specified in the forthcoming lines. 

Starting from the continuous problem~\eqref{eq:cont_form}, by simply replacing the continuous functions by their discrete counterparts, most of the resulting terms cannot be computed any more, owing to the fact that VE functions are not known explicitly in closed form. Thus, these terms need to be substituted by computable versions in the spirit of the VEM philosophy. To this purpose, the following replacements were made:
\begin{itemize}
\item The term $\Mcal\left(\frac{\partial \ch}{\partial t},\zh \right)$ in the concentration equation was replaced by
\begin{equation}\label{eq:def_Mh}
\Mcalh \left(\frac{\partial \ch}{\partial t},\zh \right):=\sum_{\E \in \taun} \McalhE\left(\frac{\partial \ch}{\partial t},\zh\right),
\end{equation}
where the local contributions are given as
\begin{equation}  \label{eq:def_Mh_loc}
\begin{split}
\McalhE\left(\ch,\zh \right)&:=\int_\E \phi \, (\PizkoE \ch) \, (\PizkoE \zh) \, \dx \\
&\qquad+ \nuME(\phi) \SEM\left((I-\PizkoE) \ch, (I-\PizkoE) \zh \right),
\end{split}
\end{equation}
with $\SEM(\cdot,\cdot)$ denoting a stabilization term with certain properties and a constant $\nuME(\phi)$, both described in Section~\ref{subsec:stabilizations} below.
\item Next, the term $\Theta(\buh,\ch;\zh)$ was substituted by 
\begin{equation} \label{eq:def_Theta_h}
\Thetah(\buh,\ch;\zh):=\frac{1}{2} \bigg[ (\buh \cdot \nabla \ch,\zh)_h + ((\qplus+\qminus)\, \ch,\zh)_h - (\buh \, \ch,\nabla \zh)_h \bigg],
\end{equation}
where
\begin{equation*} 
\begin{split}
(\buh \cdot \nabla \ch,\zh)_h&:=\sum_{\E \in \taun} \int_\E \PizbkE \buh \cdot \PizbkE(\nabla \ch) \, \PizkoE \zh \, \dx \\
((\qplus+\qminus)\, \ch,\zh)_h&:=\sum_{\E \in \taun} \int_\E (\qplus+\qminus) \, \PizkoE \ch \, \PizkoE \zh \, \dx \\
(\buh \, \ch,\nabla \zh)_h&:=\sum_{\E \in \taun} \int_\E \PizbkE \buh \, \PizkoE \ch \, \cdot \PizbkE(\nabla \zh) \, \dx.
\end{split}
\end{equation*}
\item Moreover, the term $\Dcal(\buh;\ch,\zh)$ was replaced by
\begin{equation} \label{eq:def_Dcalh}
\Dcalh \left(\buh;\ch,\zh \right):=\sum_{\E \in \taun} \DcalhE \left(\buh;\ch,\zh \right)
\end{equation}
with local contributions
\begin{equation} \label{eq:def_Dcalh_loc}
\begin{split}
\DcalhE \left(\buh;\ch,\zh \right)&:=\int_\E D(\PizbkE \buh) \, \PizbkE(\nabla \ch) \cdot  \PizbkE(\nabla \zh) \, \dx \\
&\qquad + \nuDE(\buh) \, \SED\left((I-\PinablaE) \ch,(I-\PinablaE) \zh) \right),
\end{split}
\end{equation}
where $\SED(\cdot,\cdot)$ is a stabilization term with certain properties and a constant $\nuDE(\buh)$, both described in Section~\ref{subsec:stabilizations} below.
\item Concerning $\left(\qplus\,\chat ,\zh \right)_{0,\Omega}$, this term was approximated by
\begin{equation*}
\left( \qplus\,\chat,\zh \right)_h:=\sum_{\E \in \taun} \left[\int_\E  \qplus\,\chat \, \PizkoE \zh \, \dx \right].
\end{equation*} 
\item Regarding the mixed problem, the term $\Acal(\ch;\buh,\bvh)$ was substituted by
\begin{equation} \label{eq:def_Acalh}
\Acalh(\ch;\buh,\bvh):=\sum_{\E \in \taun} \AcalhE(\ch;\buh,\bvh)
\end{equation}
with local forms
\begin{equation} \label{eq:def_Acalh_loc}
\begin{split}
\AcalhE(\ch;\buh,\bvh)&:=\int_\E A(\PizkoE \ch) \PizbkE \buh \cdot \PizbkE \bvh \, \dx \\
&\qquad + \nuAE(\ch) \, \SEA((I-\PizbkE)\buh,(I-\PizbkE)\bvh),
\end{split}
\end{equation}
where, similarly as before, $\SEA(\cdot,\cdot)$ is a stabilization term and $\nuAE(\ch)$ a constant, both described in Section~\ref{subsec:stabilizations} below.
\item Finally, the term $(\bgamma(\ch),\bvh)_{0,\Omega}$ was replaced by
\begin{equation*}
(\bgamma(\ch),\bvh)_h:=\sum_{\E \in \taun} \left[\int_\E \bgamma(\PizkoE \ch) \cdot \PizbkE \bvh \, \dx \right].
\end{equation*} 
\end{itemize}

At this point, we highlight that the bilinear form $B(\cdot,\cdot)$ needs not to be substituted since it is computable for VE functions due to the choice of degrees of freedom~\eqref{eq:dofs_VhK}. Furthermore, the right hand side term $\left( G, \qh \right)_{0,\Omega}$ remains unchanged.

\begin{remark}%\label{rem:h-idea}
Note that we here use the convention that terms which are written in caligraphic letters, such as $\Mcalh$, $\Dcalh$ and $\Acalh$, include a stabilization term, whereas those in non-caligraphic fashion and those of the form $(\cdot,\cdot)_h$ with subscript $h$ do not. In general, the terms of the type $(\cdot,\cdot)_h$ are approximations of the corresponding (possibly weighted) $L^2$ scalar products $(\cdot,\cdot)_{0,\Omega}$, obtained by introducing projections onto polynomials for all virtual functions, but not for the data terms that are known exactly. 
\end{remark}

\subsubsection{Construction of the stabilizations} \label{subsec:stabilizations}

Here, we specify the assumptions on the stabilizations $\SEM(\cdot,\cdot): \, \Zh \times \Zh \to \R$, $\SED(\cdot,\cdot):\, \Zh \times \Zh \to \R$, and $\SEA(\cdot,\cdot):\, \Vh \times \Vh \to \R$, in~\eqref{eq:def_Mh},~\eqref{eq:def_Dcalh} and~\eqref{eq:def_Acalh}, respectively. 

We require that these terms represent computable, symmetric, and positive definite bilinear forms that satisfy, for all $\E \in \taun$, the following property: there exist positive constants $M_0^{\Mcal}$, $M_1^{\Mcal}$, $M_0^{\Dcal}$, $M_1^{\Dcal}$, $M_0^{\Acal}$, $M_1^{\Acal}$, which are independent of $h$ and $\E$, such that
\begin{equation} \label{eq:stab_assumpt}
\begin{split}
M_0^{\Mcal} \lVert \zh \rVert_{0,\E}^2 \le \SEM(\zh,\zh) &\le M_1^{\Mcal} \lVert \zh \rVert_{0,\E}^2 
\qquad \forall \zh \in \Zh\cap\ker(\PizkoE) \\
M_0^{\Dcal} \lVert \nabla \zh \rVert_{0,\E}^2 \le \SED(\zh,\zh) 
& \le M_1^{\Dcal} \lVert \nabla \zh \rVert_{0,\E}^2 \qquad 
\forall \zh \in \Zh\cap\ker(\PinablaE) \\
M_0^{\Acal} \lVert \bvh \rVert_{0,\E}^2 \le \SEA(\bvh,\bvh) &\le M_1^{\Acal} \lVert \bvh \rVert_{0,\E}^2 \qquad 
\forall \bvh \in \Vh \cap \ker(\PizbkE).
\end{split}
\end{equation}
Note that continuity follows immediately from the properties:
\begin{equation*}
\SEM(\zh,\widetilde{\zh}) \le \left( \SEM(\zh,\zh) \right)^{\frac{1}{2}} \left( \SEM(\widetilde{\zh},\widetilde{\zh}) \right)^{\frac{1}{2}} \le M_1^{\Mcal} \lVert \zh \rVert_{0,\E} \lVert \widetilde{\zh} \rVert_{0,\E}
\end{equation*}
for all $\zh, \widetilde{\zh} \in \Zh\cap\ker(\PizkoE)$;
analogously for the other forms.
In practice, under mesh assumptions (\textbf{D1})-(\textbf{D2}), one can take the following scaled stabilizations corresponding to the degrees of freedom:
\begin{equation} \label{eq:stabs}
\begin{split}
\SEM(\ch,\zh)&=|\E| \sum_{j=1}^{\NZhE} \dofZhE_j(\ch) \, \dofZhE_j(\zh) \\
\SED(\ch,\zh)&=\sum_{j=1}^{\NZhE} \dofZhE_j(\ch) \, \dofZhE_j(\zh) \\
\SEA(\buh,\bvh)&=|\E|\sum_{j=1}^{\NVhE} \dofVhE_j(\buh) \, \dofVhE_j(\bvh).
\end{split}
\end{equation}
Regarding the constants appearing in front of the stabilizations in~\eqref{eq:def_Mh},~\eqref{eq:def_Dcalh} and~\eqref{eq:def_Acalh}, respectively, we pick:
\begin{equation} \label{eq:stab_constants}
\nuME(\phi)=\left|\PizzE \phi \right|, \quad \nuDE(\buh)=\nuME(\phi) (d_m+d_t |\PizbzE \buh|), \quad
\nuAE(\ch)=|A(\PizzE(\ch))|,
\end{equation}
where $\PizzE:\, L^2(\E) \to \mathbb{P}_0(\E)$ and $\PizbzE:\,  [L^2(\E)]^2 \to [\mathbb{P}_0(\E)]^2$ are the $L^2$ projectors onto scalar and vector valued constants, respectively.

\subsubsection{Well-posedness of the semidiscrete problem} %\label{subsec:well_posedness}

We first define the constants
\begin{equation*}
\nuMminus:=\min_{\E \in \taun} \nuME, \qquad \nuMplus:=\max_{\E \in \taun} \nuME.
\end{equation*}
Analogously, we introduce $\nuDminus$, $\nuDplus$, $\nuAminus$ and $\nuAplus$.
Recalling \eqref{eq:def_D} and \eqref{ass:a_phi}, it is easy to check the following (mesh-uniform) bounds for the above constants:
$$
\begin{aligned}
& \phi_\ast \le \nuMminus \le \nuMplus \le \phi^\ast  \ , \quad
(a^\ast)^{-1} \le \nuAminus \le \nuAplus \le a_\ast^{-1}  \\
& \phi_\ast d_m \le \nuDminus \le \nuDplus \le \phi^\ast (d_m + (d_\ell+d_t) \| \buh \|_{\infty,\Omega}).
\end{aligned}
$$
Then, similarly as for their continuous counterparts, the following continuity and coercivity properties for $\Mcalh(\cdot,\cdot)$ $\Dcalh(\cdot;\cdot,\cdot)$, and $\Acalh(\cdot;\cdot,\cdot)$, defined in~\eqref{eq:def_Mh},~\eqref{eq:def_Dcalh} and~\eqref{eq:def_Acalh}, respectively, hold true.

\begin{lem} \label{lem:properties_D_A}
For $\Mcalh(\cdot,\cdot)$, it holds, for all $\ch,\zh \in \Zh$,
\begin{equation} \label{eq:cont_coerc_Mcalh}
\begin{split}
\Mcalh(\ch,\zh) &\le \max\{\phi^\ast,\nuMplus M_1^{\Mcal}\} \lVert \ch \rVert_{0,\Omega} \lVert \zh \rVert_{0,\Omega} \\
\Mcalh(\zh,\zh) &\ge \min\{\phi_\ast,\nuMminus M_0^{\Mcal} \} \lVert \zh \rVert_{0,\Omega}^2.
\end{split}
\end{equation}
Concerning $\Dcalh(\cdot;\cdot,\cdot)$, this form satisfies, for all $\buh \in \Vh$ and $\ch, \zh \in \Zh$,
\begin{equation} \label{eq:cont_coerc_Dcalh}
\begin{split}
\Dcalh(\buh;\ch,\zh) &\le  \left[ \phi^\ast \left( d_m + \eta \lVert \buh \rVert_{\infty,\Omega} (d_\ell+d_t)\right)
+ \nuDplus M_1^{\Dcal} \right] |\ch|_{1,\taun} |\zh|_{1,\taun} \\
\Dcalh(\buh;\zh,\zh) &\ge \min\{ \phi_\ast d_m, \nuDminus M_0^{\Dcal} \} |\zh|^2_{1,\taun}.
\end{split}
\end{equation}
Regarding $\Acalh(\cdot;\cdot,\cdot)$, for all $\ch \in \Zh$ and $\buh,\bvh \in \Vh$, it yields
\begin{equation} \label{eq:cont_coerc_Acalh}
\begin{split}
\Acalh(\ch;\buh,\bvh) &\le \max\left\{ \frac{1}{a_\ast},\nuAplus M_1^{\Acal} \right\} \lVert \buh \rVert_{0,\Omega} \lVert \bvh \rVert_{0,\Omega}\\
\Acalh(\ch;\bvh,\bvh) &\ge \min\left\{\frac{1}{a^\ast},\nuAminus M_0^{\Acal} \right\} \lVert \bvh \rVert^2_{0,\Omega}.
\end{split}
\end{equation}
Thus, $\Acalh(\ch;\cdot,\cdot)$ is coercive on the kernel 
\begin{equation} \label{eq:discr_kernel}
\Khcal:=\{ \bvh \in \Vh: \, B(\bvh,\qh)=0 \quad \forall \qh \in \Qh \} \subset \Kcal
\end{equation}
with respect to $\lVert \cdot \rVert_{\Vh}$, where $\Kcal$ is given in~\eqref{eq:cont_kernel}.
\end{lem}

\begin{proof}
The continuity bound in~\eqref{eq:cont_coerc_Mcalh} follows directly by using
\begin{equation} \label{eq:splitting_Mh}
\Mcalh(\ch,\zh) \le \Mcalh(\ch,\ch)^{\frac{1}{2}} \Mcalh(\zh,\zh)^{\frac{1}{2}},
\end{equation}
and then estimating
\begin{equation*}
\begin{split}
\Mcalh(\ch,\ch) &\le \phi^\ast \lVert \PizkoE \ch \rVert^2_{0,\E} + \nuMplus M_1^{\Mcal} \lVert (I-\PizkoE) \ch \rVert^2_{0,\E} \\
&\le \max\{\phi^\ast,\nuMplus M_1^{\Mcal}\} \left( \lVert \PizkoE \ch \rVert^2_{0,\E} +  \lVert (I-\PizkoE) \ch \rVert^2_{0,\E} \right) \\
&= \max\{\phi^\ast,\nuMplus M_1^{\Mcal}\} \lVert \ch \rVert^2_{0,\E},
\end{split}
\end{equation*}
where the Pythagorean theorem was applied in the last equality.
For the coercivity bound, one can use~\eqref{ass:a_phi},~\eqref{eq:stab_assumpt}, and the Pythagorean theorem.

Regarding the continuity estimate for $\Dcalh(\cdot;\cdot,\cdot)$, by using a splitting of the form~\eqref{eq:splitting_Mh}, together with an estimate as in~\eqref{eq:cont_D_1}, one can deduce at the local level
\begin{equation} \label{eq:estimate_DcalhE}
\begin{split}
\DcalhE(\buh;\ch,\ch) &\le 
\phi^\ast \left( d_m + \eta \lVert \PizbkE \buh \rVert_{\infty,\Omega} (d_\ell+d_t)\right)
\lVert \PizbkE(\nabla \ch) \rVert^2_{0,\E} \\
& + \left( \nuDplus M_1^{\Dcal} \right) \lVert \nabla (I-\PinablaE)\ch \rVert^2_{0,\E} \\
& \le  \left[ \phi^\ast \left( d_m + \eta \lVert \PizbkE\buh \rVert_{\infty,\Omega} (d_\ell+d_t)\right)
+ \nuDplus M_1^{\Dcal} \right] |\ch|_{1,\taun}^2 .
\end{split}
\end{equation}
By application of a polynomial inverse estimate~\cite[Lemma 4.5.3]{BrennerScott}, the continuity of the $L^2$ projector, and the H\"{o}lder inequality, we further estimate
\begin{equation} \label{eq:cont_Pi_infty}
\lVert \PizbkE \buh \rVert_{\infty,\E} \le \eta \, \hE^{-1} \lVert \PizbkE \buh \rVert_{0,\E} \le \eta \, \hE^{-1} \lVert \buh \rVert_{0,\E} \le \eta \lVert \buh \rVert_{\infty,\E}.
\end{equation}
After inserting~\eqref{eq:cont_Pi_infty} into~\eqref{eq:estimate_DcalhE}, taking the splitting into account, and summing over all elements, the stated bound follows.
Concerning the coercivity bound for $\Dcalh(\cdot,\cdot)$, one can proceed similarly as in~\eqref{eq:Du_mu_mu} for the consistency part, and employ~\eqref{eq:stab_assumpt} for the stabilization term, to obtain elementwise
\begin{equation*}
\DcalhE(\buh;\zh,\zh) \ge \min\{ \phi_\ast d_m,\nuDminus M_0^{\Dcal} \} \left[ \lVert \PizbkE \nabla \zh \rVert^2_{0,\E} + \lVert \nabla (I-\PinablaE) \zh \rVert^2_{0,\E} \right].
\end{equation*}
We now note that the definitions of $\PinablaE$ and $\PizbkE$ easily yield
\begin{equation} \label{eq:I-Pinabla_I-Piz}
\lVert \nabla (I-\PinablaE)\zh \rVert_{0,\E} \ge \lVert (I-\PizbkE)\nabla\zh \rVert_{0,\E} .
\end{equation}	
The estimate then follows with~\eqref{eq:I-Pinabla_I-Piz}, the Pythagorean theorem and summation over all elements.

The estimates for $\Acalh(\cdot;\cdot,\cdot)$ are derived in a similar fashion as those for $\Mcalh(\cdot,\cdot)$, using~\eqref{ass:a_phi}. The coercivity on $\Khcal$ follows from the fact 
\begin{equation*}
\Khcal \equiv \{ \bvh \in \Vh: \, \div \bvh = 0 \} \subset \Kcal,
\end{equation*}
owing to the definition of $\VhE$ in~\eqref{eq:def_local_spaces}.
\end{proof}

Well-posedness of problem~\eqref{eq:semidiscr_var_form} can be shown by combining the results in~\cite{vacca2015virtual} for parabolic problems with those in~\cite{Brezzi-Falk-Marini,BBMR_generalsecondorder} for mixed problems, using Lemma~\ref{lem:properties_D_A}. More precisely, in the spirit of the two-step strategy applied in~\cite{ewing1983approximation} for FEM, one can first show that for any given $\ch(t) \in L^\infty(\Omega)$, $t \in J$, the mixed problem
\begin{equation*}
\begin{split}
\Acalh(\ch;\buh,\bvh) + B(\bvh,p_h) &= (\bgamma(\ch),\bvh)_h \\
B(\buh,\qh) &= - \left( G, \qh \right)_{0,\Omega} 
\end{split}
\end{equation*} 
admits a unique solution by applying the techniques in~\cite{Brezzi-Falk-Marini,BBMR_generalsecondorder}, and then, by using the Gronwall lemma and Picard-Lindel\"{o}f (see e.g.~\cite[Ch.1.10]{braun1983differential}), that $\ch(t)$ is uniquely determined by the discrete concentration equation
\begin{equation*}
\begin{split}
\Mcalh\left(\frac{\partial \ch}{\partial t},\zh \right) + \Thetah(\buh,\ch;\zh) + \Dcalh(\buh;\ch,\zh) = \left( \qplus \, \chat,\zh \right)_h,
\end{split}
\end{equation*}
see also~\cite{vacca2015virtual}. We do not write here the details since we focus directly on the fully discrete case, see the next section.

\subsection{Fully discrete formulation} \label{subsec:fully_discr}

Here, our goal is to formulate a fully discrete version of~\eqref{eq:semidiscr_var_form}.

To start with, we introduce a sequence of time steps $\tn=n \tau$, $n=0,\dots,N$, with time step size $\tau$. Next, we define $\bun:=\bu(\tn)$, $\pn:=p(\tn)$, $\cn:=c(\tn)$, $\Gn:=G(\tn)$, $\qplusn:=\qplus(\tn)$, and $\chatn:=\chat(\tn)$ as the
evaluations of the corresponding functions at time~$\tn$, $n=0,\dots,N$. Moreover, we denote by $\bUn\approx \bu_h(\tn)$, $\Pn \approx p_h(\tn)$ and $\Cn \approx \ch(\tn)$, the approximations of the semidiscrete solutions at those times when using a time integrator method. Among many time discretization schemes, we here make a computationally cheap choice by choosing a backward Euler method that is explicit in the nonlinear terms. The fully discrete system consequently reads as follows:

\begin{itemize}
\item for $n=0$: Given $\czh \in \Zh$, solve 
\begin{equation} \label{eq:fully_discr_model_0}
\begin{split}
\Acalh(\czh;\bUn,\bvh) + B(\bvh,\Pn) &= (\bgamma(\czh),\bvh)_h \\
B(\bUn,\qh) &= - \left( \Gn, \qh \right)_{0,\Omega}
\end{split}
\end{equation}
for all $\bvh \in \Vh$ and $\qh \in \Qh$.
\item for $n=1,\dots,N$: Solve first the concentration equation for $\Cn$:
\begin{equation} \label{eq:fully_discr_model_1}
\begin{split}
\Mcalh\left(\frac{\Cn-\Cno}{\tau},\zh \right) + \Thetah(\bUno;\Cn,\zh)+ \Dcalh(\bUno;\Cn,\zh) = \left( \qplusn \,\chatn,\zh \right)_h
\end{split}
\end{equation}
for all $\zh \in \Zh$, where $C^0:=c_{0,h}$. Then, solve the mixed problem for $\bUn$ and $\Pn$:
\begin{equation} \label{eq:fully_discr_model_2}
\begin{split}
\Acalh(\Cn;\bUn,\bvh) + B(\bvh,\Pn) &= (\bgamma(\Cn),\bvh)_h \\
B(\bUn,\qh) &= - \left( \Gn, \qh \right)_{0,\Omega}
\end{split}
\end{equation}
for all $\bvh \in \Vh$ and $\qh \in \Qh$.
\end{itemize}

\begin{lem} \label{lem:well_posedness_semidiscr}
Given $\tau>0$, provided that $\Gn,\qplusn,\Pn,\Cn \in L^\infty(\Omega)$, $\bgamma(\Cn) \in [L^2(\Omega)]^2$, and $ \bUn \in [L^\infty(\Omega)]^2$, for all $n=0,\dots,N$, the formulation~\eqref{eq:fully_discr_model_0}-\eqref{eq:fully_discr_model_2} is uniquely solvable.
\end{lem}

\begin{proof}
Similarly as for the semidiscrete case, well-posedness of~\eqref{eq:fully_discr_model_0} and~\eqref{eq:fully_discr_model_2} follows by using the tools of~\cite{Brezzi-Falk-Marini,BBMR_generalsecondorder}. Regarding~\eqref{eq:fully_discr_model_1}, we first rewrite that equation as
\begin{equation} \label{eq:fully_discr_model_4}
\begin{split}
\Mcalh\left(\Cn,\zh \right) &+ \tau \left[ \Thetah(\bUno;\Cn,\zh)+ \Dcalh(\bUno;\Cn,\zh) \right] \\
&= \tau \left( \qplusn \,\chatn,\zh \right)_h + \Mcalh\left(\Cno,\zh \right).
\end{split}
\end{equation}
We observe that all of the term are continuous with respect to  the norm $\lVert \cdot \rVert_{1,\taun}$. More precisely, for $\Mcalh(\cdot,\cdot)$ and $\Dcalh(\bUno;\cdot,\cdot)$, continuity follows from Lemma~\ref{lem:properties_D_A} and the definition of the broken $H^1$ norm. Next, for the term involving $\qplusn$, we simply apply the Cauchy-Schwarz inequality and the stability of the $L^2$ projector. Finally, for the term with $\Thetah$, we estimate
\begin{equation*}
\begin{split}
\Thetah(\bUno;\Cn,\zh) &=
\frac{1}{2} \left[\left( \bUno \cdot \nabla \Cn,\zh \right)_{h}  + ((\qplus+\qminus)\, \Cn, \zh )_{h}  - \left( \bUno \Cn, \nabla \zh \right)_{h} \right] \\
&\le \eta \left[ \lVert \bUno \rVert_{\infty,\Omega} (|\Cn|_{1,\taun} + \lVert \Cn\rVert_{0,\Omega}) + \lVert \qplus+\qminus \rVert_{\infty,\Omega} \lVert \Cn \rVert_{0,\Omega} \right] \lVert \zh \rVert_{1,\taun},
\end{split}
\end{equation*}
where we also employed an inverse inequality as in~\eqref{eq:cont_Pi_infty}. Thus, by the Lax-Milgram lemma, it only remains to show that the left hand side of~\eqref{eq:fully_discr_model_4} is coercive with respect to $\lVert \cdot \rVert_{1,\taun}$. This is however a direct consequence of 
\begin{equation*}
\Thetah(\bUno;\zh,\zh) = \frac{1}{2} ((\qplus+\qminus)\, \zh,\zh)_h \ge 0,
\end{equation*}
owing to the fact that $\qplus$ and $\qminus$ are non-negative, and the coercivity bounds~\eqref{eq:cont_coerc_Mcalh} and~\eqref{eq:cont_coerc_Dcalh}.
\end{proof}

Note that both problems~\eqref{eq:fully_discr_model_1} and~\eqref{eq:fully_discr_model_2} represent linear systems of equations which are decoupled from each other in the sense that, first, given $\chatn$ and $\qplusn$, one can determine $\Cn$ with knowledge of $\bUno$ only, and then one can use $\Cn$ to compute $\bUn$ and $\Pn$. The quantity $\Pn$ does in fact not influence the calculation of $\Cn$ directly, but rather takes the role of a Lagrange multiplier and derived variable. This decoupling, combined with the fact that the systems to be solved at each time step are linear, makes the method quite cheap per iteration.

\section{Error analysis for the fully discrete problem}\label{sec:4}

The error analysis is performed in two steps: firstly, we estimate the discretization errors for the velocity and pressure, $\lVert \bun-\bUn \rVert_{0,\Omega}$ and $\lVert \pn-\Pn \rVert_{0,\Omega}$, respectively, and then, in the second step, the concentration error $\lVert \cn - \Cn \rVert_{0,\Omega}$. In the following analysis, we assume all the needed regularity of the exact solution. Although such high regularity will not often be available in practice, the purpose of the following analysis is to give a theoretical backbone to the proposed scheme and to investigate its potential accuracy in the most favorable scenario.

\subsection{An auxiliary result}

The subsequent technical lemma will serve as an auxiliary result in the derivation of the error estimates and will be used in several occasions.

\begin{lem} \label{lem:abstr_res_1}
Let $r,s,t \in \N_0$. Denote by $\Pizr$ and $\Pizbs$, the elementwise defined $L^2$ projectors onto scalar and vector valued polynomials of degree at most $r$ and $s$, respectively. Given a scalar function $\sigma \in H^{\mr}(\taun)$, $0 \le \mr \le r+1$, let $\kappa(\sigma)$ be a tensor valued piecewise Lipschitz continuous function with respect to $\sigma$. Further, let $\sigmahat \in L^2(\Omega)$, and let $\bchi$ and $\bpsi$ be vector valued functions. We assume that $\kappa(\sigma) \in [L^\infty(\Omega)]^{2 \times 2}$, $\bchi \in [H^{\ms}(\taun) \cap L^\infty(\Omega)]^2$, $\bpsi \in [L^2(\Omega)]^2$, and $\kappa(\sigma) \bchi \in [H^{\mt}(\taun)]^2$, for some $0 \le \ms \le s+1$ and $0 \le \mt \le t+1$. Then,
\begin{equation*}
\begin{split}
& (\kappa(\sigma)\bchi,\bpsi)_{0,\Omega} - (\kappa(\Pizr \sigmahat) \Pizbs \bchi,\Pizbt \bpsi)_{0,\Omega} \\
& \le \eta \big[h^{\mt} |\kappa(\sigma) \bchi|_{\mt,\taun} + h^{\ms} |\bchi|_{\ms,\taun} \lVert \kappa(\sigma) \rVert_{\infty,\Omega}
+ (h^{\mr} |\sigma|_{\mr,\taun} + \lVert \sigma - \sigmahat \rVert_{0,\Omega}) \lVert \bchi \rVert_{\infty,\Omega} \big] \lVert \bpsi \rVert_{0,\Omega}.
\end{split}
\end{equation*}
\end{lem}

\begin{proof}
We first write
\begin{equation} \label{eq:abstr_lem_proof1}
\begin{split}
& (\kappa(\sigma)\bchi,\bpsi)_{0,\Omega} - (\kappa(\Pizr \sigmahat) \Pizbs \bchi,\Pizbt \bpsi)_{0,\Omega} \\
& = [(\kappa(\sigma)\bchi,\bpsi)_{0,\Omega} - (\kappa(\Pizr \sigma) \Pizbs \bchi,\Pizbt \bpsi)_{0,\Omega} ]
+ ((\kappa(\Pizr \sigma) - \kappa(\Pizr \sigmahat)) \Pizbs \bchi,\Pizbt \bpsi)_{0,\Omega}.
\end{split}
\end{equation}
Then, for the first part on the right hand side of~\eqref{eq:abstr_lem_proof1}, we recall that $\Pizbt$ is an $L^2$ projection and derive, on each element $\E \in \taun$,
\begin{equation*}
\begin{split}
& \left(\kappa(\sigma)\bchi,\bpsi \right)_{0,\E} - (\kappa(\PizrE \sigma) \PizbsE \bchi,\PizbtE \bpsi)_{0,\E}\\
&=[(\kappa(\sigma)\bchi,\bpsi )_{0,\E} - (\PizbtE(\kappa(\sigma)\bchi),\bpsi )_{0,\E}] \\
&\quad +[(\PizbtE(\kappa(\sigma)\bchi),\bpsi )_{0,\E} - (\PizbtE(\kappa(\sigma)\PizbsE\bchi),\bpsi )_{0,\E}] \\
&\quad +[(\PizbtE(\kappa(\sigma)\PizbsE\bchi),\bpsi )_{0,\E} - (\PizbtE(\kappa(\PizrE\sigma)\PizbsE\bchi),\bpsi )_{0,\E}]\\
&\le \eta \big[ h^{\mt} |\kappa(\sigma)\bchi|_{\mt,\E} + h^{\ms} |\bchi|_{\ms,\E} \lVert \kappa(\sigma) \rVert_{\infty,\E} + h^{\mr} |\sigma|_{\mr,\E} \lVert \PizbsE \bchi \rVert_{\infty,\E} \big] \lVert \bpsi \rVert_{0,\E}, 
\end{split}
\end{equation*}
\begin{comment}
\begin{equation} 
\begin{split}
& \left(\kappa(\sigma)\bchi,\bpsi \right)_{0,\E} - (\kappa(\PizrE \sigma) \PizbsE \bchi,\PizbtE \bpsi)_{0,\E} \\
&=(\kappa(\sigma)\bchi,\bpsi-\PizbtE \bpsi)_{0,\E} 
+ (\bchi-\PizbsE \bchi,\kappa(\sigma)\PizbtE \bpsi)_{0,\E} \\
&\qquad + ((\kappa(\sigma)-\kappa(\PizrE \sigma)) \PizbsE \bchi,\PizbtE \bpsi)_{0,\E} \\
&=(\kappa(\sigma)\bchi - \PizbsE(\kappa(\sigma)\bchi),\bpsi-\PizbsE \bpsi)_{0,\E} 
+ (\bchi-\PizbsE \bchi,\kappa(\sigma)(\PizbsE \bpsi-\bpsi))_{0,\E} \\
&\qquad +(\bchi-\PizbsE \bchi,\kappa(\sigma)\bpsi-\PizbsE(\kappa(\sigma)\bpsi))_{0,\E} 
+ ((\kappa(\sigma)-\kappa(\PizrE \sigma)) \PizbsE \bchi,\PizbsE \bpsi)_{0,\E},\\
&\le \eta \big[ h^{r} |\kappa(\sigma)\bchi|_{r,\E} + h^{r} |\bpsi|_{r,\E} \lVert \kappa(\sigma) \rVert_{\infty,\E} + h^{m} |\sigma|_{m,\E} \lVert \PizbsE \bchi \rVert_{\infty,\E} \big] \lVert \bpsi \rVert_{0,\E}, 
\end{split}
\end{equation}
\end{comment}
where in the last step we used Lemma~\ref{lem:approx_properties} and the fact that $\kappa$ is Lipschitz continuous with respect to $\sigma$. The term $\lVert \PizbsE \bchi \rVert_{\infty,\E}$ is estimated as in~\eqref{eq:cont_Pi_infty}.
%follows. By application of a polynomial inverse estimate \cite[Lemma 4.5.3]{BrennerScott}, the continuity of the $L^2$ projector and the H\"{o}lder inequality, we obtain
%\begin{equation} \label{eq:cont_Pi_infty_2}
%\lVert \PizbsE \bchi \rVert_{\infty,\E} \le \eta \hE^{-1} \lVert \PizbsE \bchi \rVert_{0,\E} \le \eta \hE^{-1} \lVert \bchi \rVert_{0,\E} \le \eta \lVert \bchi \rVert_{\infty,\E}.
%\end{equation}
Concerning the second part on the right hand side of~\eqref{eq:abstr_lem_proof1}, we have, for each $\E \in \taun$,
\begin{equation*}
\begin{split}
((\kappa(\Pizr \sigma)-\kappa(\Pizr \sigmahat)) \Pizbs \bchi,\Pizbt \bpsi)_{0,\E} &\le \lVert (\kappa(\Pizr \sigma)-\kappa(\Pizr \sigmahat) \rVert_{0,\E} \lVert \Pizbs \bchi \rVert_{\infty,\E} \lVert \Pizbt \bpsi \rVert_{0,\E} \\
&\le \lVert \sigma - \sigmahat\rVert_{0,\E} \lVert \bchi \rVert_{\infty,\E} \lVert \bpsi \rVert_{0,\E},
\end{split}
\end{equation*}
where we used again the Lipschitz continuity of $\kappa$, the continuity properties of the $L^2$ projectors, and the bound~\eqref{eq:cont_Pi_infty}.
The assertion of the lemma follows after combining the estimates and summing over all elements.
\end{proof}
 
Note that the above lemma can be easily transferred to the cases where $\sigma$, $\kappa(\sigma)$, $\chi$, and $\psi$ are scalar, and to vector valued $\bsigma$, $\bchi$ and scalar $\kappa(\bsigma)$, $\psi$. 
 
In the special case of $\chi=1$ and vector valued $\bkappa$, an adaptation of Lemma~\ref{lem:abstr_res_1} gives
\begin{equation} \label{eq:cor_abstr_res_4}
\begin{split}
(\bkappa(\sigma),\bpsi)_{0,\Omega} - (\bkappa(\Pizr \sigmahat),\Pizbt \bpsi)_{0,\Omega} \le \eta \big[h^{\mt} |\bkappa(\sigma)|_{\mt,\taun} + h^{\mr} |\sigma|_{\mr,\taun} + \lVert \sigma - \sigmahat \rVert_{0,\Omega} \big] \lVert \bpsi \rVert_{0,\Omega}.
\end{split}
\end{equation}

\subsection{Bounds on $\lVert \bun-\bUn \rVert_{0,\Omega}$ and $\lVert \pn-\Pn \rVert_{0,\Omega}$} 
 
We consider the mixed problem
\begin{equation} \label{eq:mixed_discr_model}
\begin{split}
\Acalh(\Cn;\bUn,\bvh) + B(\bvh,\Pn) &= (\bgamma(\Cn),\bvh)_h \\
B(\bUn,\qh) &= - \left( \Gn, \qh \right)_{0,\Omega} ,
\end{split}
\end{equation}
where $\Cn \in \Zh$ is the numerical solution of the concentration equation~\eqref{eq:fully_discr_model_1} for $n=1,\dots,N$, and $C^0=\czh$. The goal is to derive an upper bound for $\lVert \bun-\bUn \rVert_{0,\Omega}$ and $\lVert \pn-\Pn \rVert_{0,\Omega}$ with respect to $\lVert \cn-\Cn \rVert_{0,\Omega}$. For the analysis, we basically follow the ideas of~\cite{Brezzi-Falk-Marini,BBMR_generalsecondorder} with the major differences that, here, $\Acalh(\Cn;\cdot,\cdot)$ is not consistent with respect to $[\mathbb{P}_k(\E)]^2$ due to presence of $\Cn$, and, additionally, the right hand side of~\eqref{eq:mixed_discr_model} is inhomogeneous.

\begin{thm} \label{thm:Un_un_Pn_pn}
Given $\Cn \in \Zh$, let $(\bUn,\Pn) \in \Vh \times \Qh$ be the solution to~\eqref{eq:mixed_discr_model}. Let us assume that for the exact solution $(\bun,\pn,\cn)$ to~\eqref{eq:cont_form} at time $\tn$, it holds $\bun \in [H^{k+1}(\taun)]^2$, $\pn \in H^{k+1}(\taun)$, and $\cn \in H^{k+1}(\taun)$. Furthermore, we suppose that $\bgamma(c)$ and $A(c)$ are piecewise Lipschitz continuous functions with respect to $c \in L^2(\Omega)$, and that $\bgamma(\cn), A(\cn)\bun \in [H^{k+1}(\taun)]^2$. Then, the following error estimates hold for all $k \in \N_0$:
\begin{equation*}
\begin{split}
\lVert \bUn-\bun \rVert_{0,\Omega} &\le \lVert \Cn-\cn \rVert_{0,\Omega} \, \zeta_1^n(\bun) + h^{k+1} \, \zeta_2^n(\bun,\cn,\bgamma(\cn),A(\cn)\bun)  \\
\lVert \Pn-\pn \rVert_{0,\Omega} &\le \lVert \Cn-\cn \rVert_{0,\Omega} \, \zeta_3^n(\bun) + h^{k+1} \, \zeta_4^n(\bun,\cn,\bgamma(\cn),A(\cn)\bun,\pn),
\end{split}
\end{equation*}
where $\zeta_1^n$-$\zeta_4^n$ are positive constants independent of $h$ and depending only on the specified functions.
\end{thm}

\begin{proof}
The estimate for $\lVert \bUn-\bun \rVert_{0,\Omega}$ can be obtained as follows. 

By using the second equality in~\eqref{eq:mixed_discr_model}, we have $\div \bUn=\Pi^0_k G^n$ (use that $\div \bUn \in \mathbb{P}_k(\E)$ for every $\E \in \taun$), where we recall that $(\Pi^0_k)_{|_\E}=\PizkE$. Define now the interpolant $\bunI \in \Vh$ via the degrees of freedom~\eqref{eq:dofs_VhK}:
\begin{equation*}
\dofVh_i(\bunI)=\dofVh_i(\bun), \quad i=1,\dots,\NVh.
\end{equation*}
Then, it holds~\cite[eq.(28)]{BBMR_generalsecondorder}
\begin{equation} \label{eq:interpolation_error}
\lVert \bun-\bunI \rVert_{0,\Omega} \le \eta \, h^{k+1} \lVert \bun \rVert_{k+1,\taun}.
\end{equation}
Moreover, one has $\div \bunI=\Pi^0_k \Gn$. Thus, setting $\bdeltan:=\bUn-\bunI$, it holds that $\bdeltan \in \Khcal \subset \Kcal$, where $\Khcal$ and $\Kcal$ were defined in~\eqref{eq:discr_kernel} and~\eqref{eq:cont_kernel}, respectively, and therefore, $\lVert \bdeltan \rVert_{\Vh}=\lVert \bdeltan \rVert_{0,\Omega}$. Owing to the assumptions on $a(\cdot)$ in~\eqref{ass:a_phi} together with~\eqref{eq:cont_coerc_Acalh}, we have, further using~\eqref{eq:mixed_discr_model} with $\bvh=\bdeltan \in \Khcal$ and~\eqref{eq:cont_form},
\begin{equation} \label{eq:error_estimate_Un_un}
\begin{split}
\alpha \lVert \bdeltan \rVert^2_{0,\Omega} &\le \Acalh(\Cn;\bdeltan,\bdeltan) = \Acalh(\Cn;\bUn,\bdeltan) - \Acalh(\Cn;\bunI,\bdeltan)  \\
&= (\bgamma(\Cn),\bdeltan)_h - \Acalh(\Cn;\bunI,\bdeltan) \\
&=\left[(\bgamma(\Cn),\bdeltan)_h - (\bgamma(\cn),\bdeltan)_{0,\Omega} \right]
+ \Acalh(\Cn;\bun-\bunI,\bdeltan)  \\
&\quad + \bigg[\Acal(\cn;\bun,\bdeltan) - \Acalh(\Cn;\bun,\bdeltan) \bigg] \\
&=: T_1 + T_2 + T_3.
\end{split}
\end{equation}

The terms $T_1$-$T_3$ are bounded as follows:
\begin{itemize}
\item term $T_1$: We use equation~\eqref{eq:cor_abstr_res_4} with $\bkappa=\bgamma$, $\sigma=\cn$, $\sigmahat=\Cn$, $\bpsi=\bdeltan$, $r=k+1$, $t=k$, and $\mr=\mt=k+1$, and obtain
\begin{equation*}
\begin{split}
|T_1| &= |(\bgamma(\cn),\bdeltan)_{0,\Omega} - (\bgamma(\Pizko\Cn),\Pizbk\bdeltan)_{0,\Omega}| \\
&\le \eta \big[h^{k+1} (|\bgamma(\cn)|_{k+1,\taun} + |\cn|_{k+1,\taun}) + \lVert \cn - \Cn \rVert_{0,\Omega} \big] \lVert \bdeltan \rVert_{0,\Omega}.
\end{split}
\end{equation*}
\item term $T_2$: Owing to the continuity properties~\eqref{eq:cont_coerc_Acalh} of $\Acalh(\cdot;\cdot,\cdot)$ and the interpolation error estimate~\eqref{eq:interpolation_error}, it holds
\begin{equation*}
\begin{split}
|T_2|=|\Acalh(\Cn;\bun-\bunI,\bdeltan)| 
\le \eta \lVert \bun-\bunI \rVert_{0,\Omega} \lVert \bdeltan \rVert_{0,\Omega}  
\le \eta \, h^{k+1} \lVert \bun \rVert_{k+1,\taun} \lVert \bdeltan \rVert_{0,\Omega}.
\end{split}
\end{equation*}

\item term $T_3$: We have
\begin{equation*}
\begin{split}
|T_3|&=|\Acal(\cn;\bun,\bdeltan) - \Acalh(\Cn;\bun,\bdeltan)| \\
&\le |(A(\cn) \bun, \bdeltan)_{0,\Omega} - (A(\Pizko \Cn) \Pizbk \bun, \Pizbk \bdeltan)_{0,\Omega}| \\
& + \left|\sum_{\E \in \taun} \nuAE(\Cn) \, \SEA((I-\PizbkE)\bun,(I-\PizbkE)\bdeltan) \right|\\
&=:T_3^A + T_3^B.
\end{split}
\end{equation*}
For the term $T_3^A$, we use Lemma~\ref{lem:abstr_res_1} with $\kappa=A$, $\sigma=\cn$, $\sigmahat=\Cn$, $\bchi=\bun$, $\bpsi=\bdeltan$, $r=k+1$, $s=t=k$, and $\mr=\ms=\mt=k+1$, to get
\begin{equation*}
\begin{split}
T_3^A &\le \eta \bigg[h^{k+1} \big( |A(\cn) \bun|_{k+1,\taun} + |\bun|_{k+1,\taun} \lVert A(\cn) \rVert_{\infty,\Omega}
+ |\cn|_{k+1,\taun} \lVert \bun \rVert_{\infty,\Omega} \big)  \\
&\qquad+ \lVert \cn - \Cn \rVert_{0,\Omega} \lVert \bun \rVert_{\infty,\Omega} \bigg] \lVert \bdeltan \rVert_{0,\Omega}.
\end{split}
\end{equation*}
On the other hand, the term $T_3^B$ can be bounded with~\eqref{eq:stab_assumpt},~\eqref{ass:a_phi}, and Lemma~\eqref{lem:approx_properties}:
\begin{equation*}
T_3^B \le \eta \, h^{k+1} |\bun|_{k+1,\taun} \lVert \bdeltan \rVert_{0,\Omega}.
\end{equation*}
\end{itemize}
After plugging the bounds obtained for $T_1$-$T_3$ into~\eqref{eq:error_estimate_Un_un}, dividing by $\lVert \bdeltan \rVert_{0,\Omega}$, using the triangle inequality in the form
\begin{equation*}
\lVert \bUn-\bun \rVert_{0,\Omega} \le \lVert \bdeltan \rVert_{0,\Omega} +\lVert \bun-\bunI \rVert_{0,\Omega},
\end{equation*}
and employing~\eqref{eq:interpolation_error}, the convergence result follows.

The error estimate for the term $\lVert \Pn-\pn \rVert_{0,\Omega}$ follows easily by combining the above ideas with the argument in~\cite[Theorem 6.1]{Brezzi-Falk-Marini} and is therefore not shown.
\end{proof}

% -----------------------------------------------------------
\subsection{Bounds on $\lVert \cn-\Cn \rVert_{0,\Omega}$} 
% -----------------------------------------------------------

For fixed $\bu(t) \in \boldsymbol{V}$ and $t\in J$, we define the projector $\Pc:\, Z \to \Zh$
(that to each $c \in Z$ associates $\Pc c \in \Zh$) by
\begin{equation} \label{eq:proj_Pc}
\begin{split}
\Gammach(\bu(t);\Pc c,\zh)=\Gammac(\bu(t);c,\zh),
\end{split}
\end{equation}
for all $\zh \in \Zh$, where 
\begin{equation} \label{def:Gammalambda_Gammalambdah}
\begin{split}
\Gammach(\bu;c,\zh)&:=\Dcalh(\bu;c,\zh) + \Thetah(\bu;c,\zh) + (c,\zh)_h \\
\Gammac(\bu;c,\zh)&:=\Dcal(\bu;c,\zh) + \Theta(\bu;c,\zh) + (c,\zh)_{0,\Omega},
\end{split}
\end{equation}
with 
\begin{equation*}\label{eq:L1}
(c,\zh)_h:=\sum_{\E \in \taun} \int_\E c \, (\PizkoE \zh) \, \dx.
\end{equation*}
\begin{comment}
and the function $\lambda$ chosen as follows. On each element $\E \in \taun$, $\lambda_{|_\E}$ is constant with value
\begin{equation} \label{def:lambdan}
\lambda_{|_\E}(\bx)=\max\left\{ 0,1-\frac{1}{2} \min_{\by \in \E} G(\by) \right\},
\end{equation}
where we assumed here that $G \in L^\infty(\Omega)$.
With this choice it follows that $\frac{1}{2} G + \lambda \ge 1$ in~$\Omega$. Further, we highlight that this particular choice of $\lambda$ is needed in the proof of Lemma~\ref{lem:well_posedness_Pc}.
\end{comment}

\begin{lem} \label{lem:well_posedness_Pc}
The projector $\Pc:\, Z \to \Zh$ given in~\eqref{eq:proj_Pc} is well-defined under the assumption that $\bu$, $\qplus$, and $\qminus$ are bounded in $L^\infty(\Omega)$ for all $t\in J$.
\end{lem}
\begin{proof}
By the Lax-Milgram lemma, we have to show that the left hand side of~\eqref{eq:proj_Pc} defines a continuous and coercive bilinear form and that the right hand side is a continuous functional with respect to $\lVert \cdot \rVert_{1,\taun}$. Continuity of the latter one is obtained by combining~\eqref{eq:cont_D_1} with
\begin{equation*}
\begin{split}
\Theta(\bu;c,\zh) &+ (c,\zh)_{0,\Omega}=
\frac{1}{2} \left[\left( \bu \cdot \nabla c,\zh \right)_{0,\Omega}  + ((\qplus+\qminus+2) c, \zh )_{0,\Omega}  - \left( \bu \, c, \nabla \zh \right)_{0,\Omega} \right] \\
&\le \frac{1}{2} \left[ \lVert \bu \rVert_{\infty,\Omega} (|c|_{1,\taun} + \lVert c\rVert_{0,\Omega}) + \lVert \qplus+\qminus+2 \rVert_{\infty,\Omega} \lVert c \rVert_{0,\Omega} \right] \lVert \zh \rVert_{1,\taun}.
\end{split}
\end{equation*}
By using~\eqref{eq:cont_coerc_Dcalh} and performing similar computations as in the proof of Lemma~\ref{lem:well_posedness_semidiscr}, continuity of $\Gammach$ follows:
\begin{equation} \label{eq:cont_Gammalambdah}
\Gammach(\bu;c,\zh) \le \eta \, \zeta(\bu,\qplus,\qminus) \lVert c \rVert_{1,\taun} \lVert \zh \rVert_{1,\taun},
\end{equation}
where $\zeta$ only depends on the specified functions.
Regarding the coercivity of $\Gammach$, we first estimate
\begin{equation*}
\begin{split}
\Thetah(\bu;\zh,\zh) + (\zh,\zh)_h 
=\sum_{\E \in \taun} \left( \left( \frac{1}{2} (\qplus+\qminus) + 1 \right) \PizkoE \zh, \PizkoE \zh \right)_{0,\E} \ge \lVert \Pizko \zh \rVert^2_{0,\Omega},
\end{split}
\end{equation*}
where we recall that $(\Pizko)_{|_\E}=\PizkoE$ for all $\E \in \taun$. Then, combining this result with~\eqref{eq:cont_coerc_Dcalh} yields
\begin{equation*} 
\Gammach(\bu;\zh,\zh) \ge \eta \left[\left| \zh \right|^2_{1,\taun} + \lVert \PizkoE \zh \rVert^2_{0,\Omega} \right]
\ge \eta \left[\left| \zh \right|^2_{1,\taun} + \lVert \overline{\zh} \rVert^2_{0,\Omega} \right],
\end{equation*}
with $\overline{\zh}$ denoting the $L^2(\Omega)$ projection of $\zh$ onto $\mathbb{P}_0(\Omega)$. 
\begin{comment}
This follows from
\begin{equation*}
\lVert \PizzE \zh \rVert^2_{0,\E} \le \lVert \PizkoE \zh \rVert^2_{0,\E}.
\end{equation*}
due to
\begin{equation*}
\begin{split}
\lVert \PizzE \zh \rVert^2_{0,\E} &= (\PizzE \zh,\PizzE \zh)_{0,\E}
= (\zh,\PizzE \zh)_{0,\E} = (\zh,\PizkoE(\PizzE \zh))_{0,\E} \\
&= (\PizkoE \zh,\PizzE \zh)_{0,\E} \le \lVert \PizkoE \zh \rVert_{0,\E} \lVert \PizzE \zh \rVert_{0,\E}.
\end{split}
\end{equation*}
\end{comment}
Since $\overline{\zh}$ coincides with the average of $\zh$, one can use a Poincar\'{e}-Friedrichs inequality, see e.g.~\cite{brenner2003poincare}, to deduce 
\begin{equation*}
\left| \zh \right|^2_{1,\taun} + \lVert \overline{\zh} \rVert^2_{0,\Omega} \ge C_p^{-1} \diam(\Omega)^{-1} \lVert \zh \rVert^2_{1,\taun},
\end{equation*}
and consequently the coercivity of $\Gammach$.
\end{proof}

\begin{comment}
the inverse triangle inequality and the Poincar\'e-Friedrichs inequality imply
\begin{equation*}
\lVert \zh \rVert_{0,\E} - \lVert \PizzE \zh \rVert_{0,\E}
\le \lVert \zh - \PizzE \zh \rVert_{0,\E} \le C_p |\zh|_{1,\E}.
\end{equation*}
This leads to 
\begin{equation*}
\begin{split}
\frac{\widetilde{M}}{2} \left| \zh \right|^2_{1,\E} + \lVert \PizzE \zh \rVert^2_{0,\E} 
&\ge \frac{\widetilde{M}}{2 C_p^2} \lVert \zh - \PizzE \zh \rVert^2_{0,\E} +  \lVert \PizzE \zh \rVert^2_{0,\E}\\
& \ge \min\left\{ \frac{\widetilde{M}}{2 C_p^2},1 \right\}  \left( \lVert \zh - \PizzE \zh \rVert^2_{0,\E} +  \lVert \PizzE \zh \rVert^2_{0,\E} \right) \\
& \ge \widetilde{C_p} \lVert \zh \rVert^2_{0,\E}.
\end{split}
\end{equation*}
with $\widetilde{C_p}:=\min\left\{ \frac{\widetilde{M}}{4 C_p^2},\frac{1}{2} \right\} $.
Thus,
\begin{equation*}
(D(\bu) \nabla \zh,\nabla \zh)_h +\Thetah(\bu;\zh,\zh) + (\lambda \zh,\zh)_h \ge \sum_{\E \in \taun} \min\left\{ \frac{\widetilde{M}}{4}, \frac{\widetilde{C_p}}{2} \right\} \lVert \zh \rVert^2_{1,\E} \ge M  \lVert \zh \rVert^2_{1,\Omega}
\end{equation*}
with $M:=\min_{\E \in \taun} \left\{ \min\left\{ \frac{\widetilde{M}}{4}, \frac{\widetilde{C_p}}{2} \right\} \right\}$.
\end{itemize}
\end{proof}
\end{comment}

\begin{lem} \label{lem:cn_Pcn}
We assume that $\bu \in [H^{k+1}(\taun) \cap L^\infty(\Omega)]^2$, $c \in H^{k+2}(\taun)\cap W^{1,\infty}(\taun)$, $\qplus,\qminus \in L^\infty(\Omega)$, $(\qplus+\qminus) c \in H^{k+1}(\taun)$, $\bu \, c \in [H^{k+1}(\taun)]^2$, $\bu \cdot \nabla c \in H^{k+1}(\taun)$, and $D(\bu) \nabla c \in [H^{k+1}(\taun)]^2$ for all $t \in J$. Then, the following error bounds for $c-\Pc c$, where $\Pc c$ is defined in~\eqref{eq:proj_Pc}, hold for all $k \in \N_0$: 
\begin{equation} \label{eq:bound_rho}
\begin{split}
\lVert c -\Pc c \rVert_{1,\taun} &\le  h^{k+1} \, \xi_1(c,\bu,\qplus,\qminus,D(\bu)\nabla c,\nabla c,(\qplus+\qminus) c, \bu \cdot \nabla c, \bu \, c), \\
\lVert c -\Pc c \rVert_{0,\Omega} &\le  h^{k+2} \, \xi_0(c,\bu,\qplus,\qminus,D(\bu)\nabla c,\nabla c,(\qplus+\qminus) c, \bu \cdot \nabla c, \bu \, c),
\end{split}
\end{equation}
where the constants $\xi_1,\xi_0>0$ only depend on the listed terms and are independent of $h$.
\end{lem}

\begin{proof}
We focus on the error estimate in the broken $H^1$ norm at a fixed time $t \in J$. First, we state the following result. Given $c \in H^{k+2}(\taun)$, there exists an interpolant $c_I \in \Zh$ such that the following bounds hold true (see for instance~\cite{cangianigeorgulispryersutton_VEMaposteriori,beiraolovadinarusso_stabilityVEM,BrennerGuanSung_someestimatesVEM}):
\begin{equation} \label{eq:interpolation_error_cn}
\lVert c-c_I \rVert_{0,\Omega} \le \eta \, h^{k+2} \lVert c \rVert_{k+2,\taun}, \quad \lVert c-c_I \rVert_{1,\taun} \le \eta \, h^{k+1} \lVert c \rVert_{k+2,\taun}.
\end{equation}
After denoting $\nu:=\Pc c-c_I$, one obtains with the coercivity of $\Gammach$, see the proof of Lemma~\ref{lem:well_posedness_Pc}, and the definition of $\Pc c$ in~\eqref{eq:proj_Pc},
\begin{equation} \label{eq:estimate_S1_S2_S3}
\begin{split}
M \lVert \nu \rVert^2_{1,\taun} &\le \Gammach(\bu,\nu,\nu)
= \Gammach(\bu,\Pc c,\nu) - \Gammach(\bu,c_I,\nu) \\
&= [\Gammac(\bu,c,\nu) - \Gammach(\bu,c,\nu)] + \Gammach(\bu,c-c_I,\nu) \\
&=: S_1 + S_2,
\end{split}
\end{equation}
for a constant $M >0$. By employing the definitions of $\Gammac$ and $\Gammach$ in~\eqref{def:Gammalambda_Gammalambdah}, the term $S_1$ is split as follows:
\begin{equation*}
\begin{split}
S_1 &= [\Dcal(\bu;c,\nu)-\Dcalh(\bu;c,\nu)] + [\Theta(\bu;c,\nu)-\Thetah(\bu;c,\nu)] + [(c,\nu)_{0,\Omega}-(c,\nu)_h] \\
&=:S_1^A + S_1^B + S_1^C.
\end{split}
\end{equation*}

For $S_1^A$, we have
\begin{equation*}
\begin{split}
S_1^A&= [(D(\bu)\nabla c,\nabla \nu)_{0,\Omega} - (D(\Pizbk \bu) \, \Pizbk(\nabla c) ,  \Pizbk(\nabla \nu))_{0,\Omega}]\\
&\quad + \sum_{\E \in \taun} \nuDE(\bu) \SED((I-\PinablaE)c,(I-\PinablaE) \nu) \\
&\le \eta \, h^{k+1} \left[ |D(\bu) \nabla c|_{k+1,\taun} + |\nabla c|_{k+1,\taun} (\lVert D(\bu) \rVert_{\infty,\Omega}+1) + |\bu|_{k+1,\taun} \lVert \nabla c \rVert_{\infty,\Omega} \right] |\nu|_{1,\taun},
\end{split}
\end{equation*}
where in the inequality we applied Lemma~\ref{lem:abstr_res_1} to estimate the first part on the right hand side of $S_1^A$, and made use of the continuity properties~\eqref{eq:stab_assumpt} of $\SED(\cdot,\cdot)$, the trivial continuity property of $\PinablaE$ in the $H^1$ seminorm and its approximation properties (stated in Lemma~\ref{lem:approx_properties}) to estimate the stabilization term.

Next, for $S_1^B$, we compute
\begin{equation*}
\begin{split}
S_1^B &= \frac{1}{2} \bigg\{  \left[(\bu \cdot \nabla c,\nu)_{0,\Omega} - (\Pizbk \bu \cdot \Pizbk(\nabla c) ,\Pizko \nu)_{0,\Omega} \right] \\
&\qquad + \left[ ((\qplus+\qminus)c,\nu)_{0,\Omega} - ((\qplus+\qminus) \Pizko c, \Pizko \nu)_{0,\Omega} \right] \\
&\qquad - \left[(\bu \, c,\nabla \nu)_{0,\Omega} - (\Pizbk \bu \, \Pizko c,\Pizbk(\nabla \nu))_{0,\Omega} \right] \bigg\} \\
&\le \eta \, h^{k+1} \big[ |\bu \cdot \nabla c|_{k+1,\taun} +  (|\nabla c|_{k+1,\taun} + |c|_{k+1,\taun}) \lVert \bu \rVert_{\infty,\Omega} + |c|_{k+1,\taun} \lVert \qplus + \qminus \rVert_{\infty,\Omega}  \\
& \qquad\qquad +|(\qplus+\qminus) c|_{k+1,\taun} + |\bu \, c|_{k+1,\taun} + |\bu|_{k+1,\taun} (\lVert c \rVert_{\infty,\Omega} + \lVert \nabla c \rVert_{\infty,\Omega}) \big] \lVert \nu \rVert_{1,\taun},
\end{split}
\end{equation*}
where in the last inequality we used Lemma~\ref{lem:abstr_res_1} with $\kappa=id$ and $\sigma=\bu$ for the first and third term inside the curly bracket, and $\kappa=\qplus+\qminus$ and $\sigma=1$ for the second one.

Finally, for $S_1^C$, it holds with the definition of the $L^2$ projector and Lemma~\ref{lem:approx_properties}
\begin{equation*}
\begin{split}
S_1^C &= ((I-\Pizko) c,\nu)_{0,\Omega} \le \eta \, h^{k+1} |c|_{k+1,\taun} \lVert \nu \rVert_{0,\Omega}.
\end{split}
\end{equation*}
On the other hand, for $S_2$, we use the continuity of $\Gammach$ in~\eqref{eq:cont_Gammalambdah}, together with the interpolation error estimate~\eqref{eq:interpolation_error_cn}, to derive
\begin{equation*}
\Gammach(\bu;c-c_I,\nu) \le \eta \, \zeta(\bu,\qplus,\qminus) \lVert c-c_I \rVert_{1,\taun} \lVert \nu \rVert_{1,\taun} 
\le \eta \,\zeta(\bu,\qplus,\qminus) h^{k+1} \lVert c \rVert_{k+2,\taun} \lVert \nu \rVert_{1,\taun}.
\end{equation*}
The error bound in the broken $H^1$ norm follows by plugging first the estimates for $S_1^A$, $S_1^B$, and $S_1^C$ into $S_1$, then those obtained for $S_1$ and $S_2$ into~\eqref{eq:estimate_S1_S2_S3}, using the definition of the $H^1$ norm, dividing by $\lVert \nu \rVert_{1,\taun}$, and using the triangle inequality in the form
\begin{equation*}
\lVert c -\Pc c \rVert_{1,\taun} \le \lVert c -c_I \rVert_{1,\taun} + \lVert \nu \rVert_{1,\taun},
\end{equation*}
together with the approximation properties~\eqref{eq:interpolation_error_cn} of the interpolant $c_I$.

The $L^2$ error bound can be derived by combining the above arguments with a standard duality argument as in~\cite{vacca2015virtual}, also recalling the convexity of $\Omega$; it is omitted here.
\end{proof}

%%%%
%\begin{remark} \label{rem:higher_reg_c}
%We highlight that in Lemma~\ref{lem:cn_Pcn} above, we assumed additional regularity on $c$, namely we required that $c \in H^{k+2}(\taun)$ instead of $c \in H^{k+1}(\taun)$. By this assumption, one recovers in fact the ``expected'' rate of convergence: for instance, for the lowest-order case $k=0$, the space $\ZhEtilde$ given in~\eqref{eq:def_local_spaces} is exactly the one introduced in~\cite{VEMvolley,hitchhikersguideVEM}, for which can show a linear and quadratic order of convergence for the $H^1$ and $L^2$ discretization errors, respectively, provided that the solution is at least in $H^2(\taun)$.
%\end{remark}
%%%%

By differentiation of~\eqref{eq:proj_Pc} in time and use of similar techniques as in the proof of Lemma~\ref{lem:cn_Pcn}, an analogous result can be obtained for $\frac{\partial}{\partial t}(c-\Pc c)$, summarized in the following corollary.

\begin{cor} \label{cor:bound_partial_t_c_Pc}
Provided that the continuous data and solution are sufficiently regular in space and time, it holds 
\begin{equation*}
\begin{split}
\left\lVert \frac{\partial}{\partial t} (c -\Pc c) \right\rVert_{1,\taun} \le  h^{k+1} \, \xi_{1,t}, \qquad \left\lVert \frac{\partial}{\partial t} (c -\Pc c) \right\rVert_{0,\Omega} \le  h^{k+2} \, \xi_{0,t},
\end{split}
\end{equation*}
where the constants $\xi_{1,t},\xi_{0,t}>0$ are independent of $h$.
\end{cor}

Moreover, we will later on need the two subsequent bounds. 
\begin{lem} \label{lem:partial_cn}
Under sufficient smoothness of the continuous data and solution, it holds
\begin{equation*}
\left\lVert \frac{\partial \cn}{\partial t} - \frac{\Pc \cn - \Pc \cno}{\tau} \right\rVert_{0,\Omega}
\le \tau^{\frac{1}{2}} \left\lVert \frac{\partial^2 c}{\partial s^2} \right\rVert_{L^2(\tnmo,\tn;L^2(\Omega))} 
+\tau^{-\frac{1}{2}} h^{k+2} \left(\int_{\tnmo}^{\tn} \xi_{0,t}^2 \, \ds \right)^{\frac{1}{2}},
\end{equation*}
where $\xi_{0,t}$ can be found in Corollary~\ref{cor:bound_partial_t_c_Pc}.
\end{lem}
\begin{proof}
We estimate
\begin{equation*} %\label{eq:estimate_partial_cn}
\begin{split}
\left\lVert \frac{\partial c}{\partial t} - \frac{\Pc \cn - \Pc \cno}{\tau} \right\rVert_{0,\Omega} 
&\le \left\lVert \frac{\partial \cn}{\partial t} - \frac{\cn-\cno}{\tau} \right\rVert_{0,\Omega} + \left\lVert \frac{\Pc\cn-\Pc\cno}{\tau} - \frac{\cn-\cno}{\tau} \right\rVert_{0,\Omega}\\
&=:(I) + (II).
\end{split} 
\end{equation*}
The term $(I)$ can be estimated exactly as for standard finite elements, see for instance~\cite{thomee1984galerkin}:
\begin{equation*} \label{eq:term_I}
\begin{split}
&(I)=\left\lVert \frac{\partial \cn}{\partial t} - \frac{\cn-\cno}{\tau} \right\rVert_{0,\Omega}
\le \int_{\tnmo}^{\tn} \left\lVert \frac{\partial^2 c}{\partial s^2}(s)\right\rVert_{0,\Omega} \, \ds 
\le \tau^{\frac{1}{2}} \left(\int_{\tnmo}^{\tn} \left\lVert \frac{\partial^2 c}{\partial s^2}(s)\right\rVert^2_{0,\Omega} \, \ds \right)^{\frac{1}{2}},
\end{split}
\end{equation*}
where we also applied the H\"{o}lder inequality in the last step. 
Concerning $(II)$, this term can be bounded as follows, using Corollary~\ref{cor:bound_partial_t_c_Pc}:
\begin{equation*}
\begin{split}
(II) &= \left\lVert \frac{\Pc\cn-\Pc\cno}{\tau} - \frac{\cn-\cno}{\tau} \right\rVert_{0,\Omega}
= \frac{1}{\tau}\left\lVert \int_{\tnmo}^{\tn} \frac{\partial}{\partial s} (\Pc c-c)(s) \, \ds \right\rVert_{0,\Omega} \\
&\le \tau^{-\frac{1}{2}} \left(\int_{\tnmo}^{\tn} \left\lVert \frac{\partial}{\partial s} (\Pc c-c)(s)\right\rVert^2_{0,\Omega} \, \ds \right)^{\frac{1}{2}} 
\le \tau^{-\frac{1}{2}} h^{k+2} \left(\int_{\tnmo}^{\tn} \xi_{0,t}^2 \, \ds \right)^{\frac{1}{2}}.
\end{split}
\end{equation*}
The statement of the lemma follows.
\end{proof}

\begin{lem} \label{lem:un_Uno}
Provided that the continuous data and solution are sufficiently regular in space and time, it holds
\begin{equation*} 
\lVert \bun - \bUno \rVert_{0,\Omega} \le \tau \left\lVert \frac{\partial \bu}{\partial t} \right\rVert_{L^{\infty}(\tnmo,\tn;L^2(\Omega))} + \lVert \Cno-\cno \rVert_{0,\Omega} \, \zeta_1^{n-1} + h^{k+1} \, \zeta_2^{n-1},
\end{equation*}
where $\zeta_1^{n-1}$ and $\zeta_2^{n-1}$ are the constants from Theorem~\ref{thm:Un_un_Pn_pn}.
\begin{proof}
By using the triangle inequality, one obtains
\begin{equation*} %\label{eq:estimate_un_Uno}
\begin{split}
\lVert \bun - \bUno \rVert_{0,\Omega} &\le \left\lVert \bun-\buno \right\rVert_{0,\Omega} + \lVert \buno-\bUno \rVert_{0,\Omega}.
\end{split}
\end{equation*}
The first term on the right hand side is estimated by
\begin{equation*}
\lVert \bun-\buno \rVert_{0,\Omega} = \left\lVert \int_{\tnmo}^{\tn} \frac{\partial \bu(s)}{\partial s} \, \ds \right\rVert_{0,\Omega} 
\le \tau \left\lVert \frac{\partial \bu}{\partial t} \right\rVert_{L^{\infty}(\tnmo,\tn;L^2(\Omega))},
\end{equation*}
and the second one term is bounded with Theorem~\ref{thm:Un_un_Pn_pn}.
\end{proof}
\end{lem}

Now, we have all the ingredients to bound $\lVert \cn-\Cn \rVert_{0,\Omega}$. 

\begin{thm} \label{thm:Cn_cn}
Let the mesh assumptions (\textbf{D1})-(\textbf{D3}) be satisfied. Then, provided that the continuous data and solutions are sufficiently regular, it yields
\begin{equation*}
\lVert \cn-\Cn \rVert_{0,\Omega} \le \eta \left[ \lVert \czh-c^0 \rVert_{0,\Omega} + h^{k+1} \,\varphi_1 + \tau \: \varphi_2 \right],
\end{equation*}
where the regularity terms  $\varphi_1, \varphi_2$ and the positive constant $\eta$ now depend on $\bu$, $c$, $\qplus$, $\qminus$, $\chat$, $\frac{\partial \bu}{\partial t}$, $\frac{\partial^2 \bu}{\partial t^2}$, $\frac{\partial c}{\partial t}$, and $\frac{\partial^2 c}{\partial t^2}$ (and products of these functions).
\end{thm}
%%
% \ap{\begin{remark}
% We highlight that no CFL condition is needed in Theorem~\ref{thm:Cn_cn}.
% \end{remark}}
%%
\begin{proof}
To start with, we write
\begin{equation*}
\Cn-\cn = (\Cn - \Pc \cn) + (\Pc \cn - \cn)=:\varthetan+\rhon.
\end{equation*}
Equation~\eqref{eq:bound_rho} gives a bound on $\rhon$. In order to deal with $\varthetan$, we use the continuous concentration equation~\eqref{eq:cont_form_altern} with $z=\varthetan$, the fully discretized version~\eqref{eq:fully_discr_model_1} with $\zh=\varthetan$, and the definition of the projector $\Pc \cn$ in~\eqref{eq:proj_Pc} with $\zh=\varthetan$: 
\begin{equation} \label{eq:rec_rel_theta}
\begin{split}
&\Mcalh\left( \frac{\varthetan-\varthetano}{\tau}, \varthetan \right) 
+ \Dcalh(\bUno;\varthetan,\varthetan)\\
&=\left[ \Mcal\left(\frac{\partial \cn}{\partial t},\varthetan \right)_{0,\Omega} - \Mcalh\left( \frac{\Pc \cn - \Pc \cno}{\tau},\varthetan \right) \right] \\
&\quad + \left[ \Thetah(\bun;\Pc \cn,\varthetan) - \Thetah(\bUno;\Cn,\varthetan) \right] 
+ [\Dcalh\left(\bun;\Pc \cn, \varthetan \right) - \Dcalh\left(\bUno;\Pc \cn, \varthetan \right)] \\
&\quad + \left[ (\Pc \cn,\varthetan)_h-(\cn,\varthetan)_{0,\Omega} \right] + \left[ (\qplusn\chatn,\varthetan)_h - (\qplusn\chatn,\varthetan)_{0,\Omega} \right]\\
&\quad =:R_1 + R_2 + R_3 + R_4 + R_5.
\end{split}
\end{equation}
Owing to the coercivity properties in~\eqref{eq:cont_coerc_Dcalh}, the second term on the left hand side of~\eqref{eq:rec_rel_theta} can be estimated by
\begin{equation} \label{eq:estimate:Dh}
\Dcalh(\bUno;\varthetan,\varthetan) \ge  D_\ast \left| \varthetan \right|^2_{1,\taun},
\end{equation}
with some constant $D_\ast>0$ independent of $h$ and $\bUno$. 

The terms $R_1$-$R_5$ on the right hand side of~\eqref{eq:rec_rel_theta} are estimated as follows:
\begin{itemize}
\item term $R_1$: Using the definition of $\Mcalh(\cdot,\cdot)$ in~\eqref{eq:def_Mh} together with~\eqref{eq:stab_assumpt} yields
\begin{equation*}% \label{eq:R_1}
\begin{split}
R_1 &= \Mcal\left(\frac{\partial \cn}{\partial t},\varthetan \right)_{0,\Omega} - \Mcalh\left( \frac{\Pc \cn - \Pc \cno}{\tau},\varthetan \right) \\
&= \bigg[ \left(\phi \, \frac{\partial \cn}{\partial t}, \varthetan \right)_{0,\Omega} - \left(\Pizko \left(\phi \, \Pizko \left( \frac{\Pc \cn - \Pc \cno}{\tau} \right)\right), \varthetan \right)_{0,\Omega} \\
& \qquad - \sum_{\E \in \taun} \nuME(\phi) \SEM\left((I-\PizkoE) \left( \frac{\Pc \cn - \Pc \cno}{\tau} \right), (I-\PizkoE) \varthetan \right) \bigg] \\
&\le \eta \bigg[\left\lVert \phi \, \frac{\partial \cn}{\partial t} - \Pizko \left(\phi \, \Pizko \left( \frac{\Pc \cn - \Pc \cno}{\tau} \right)\right) \right\rVert_{0,\Omega} \\
&\qquad+ \left\lVert (I-\Pizko) \left( \frac{\Pc \cn - \Pc \cno}{\tau} \right) \right\rVert_{0,\Omega} \bigg] \lVert \varthetan \rVert_{0,\Omega} \\
&=:\eta [R_1^A + R_1^B] \lVert \varthetan \rVert_{0,\Omega}.
\end{split}
\end{equation*}
The term $R_1^A$ is estimated by using the continuity of the $L^2$ projector, the assumption~\eqref{ass:a_phi} on $\phi$, and the approximation properties in Lemma~\eqref{lem:approx_properties}:
\begin{equation*}
\begin{split}
R_1^A &\le \left\lVert (I-\Pizko) \left(\phi \, \frac{\partial \cn}{\partial t}\right) \right\rVert_{0,\Omega}
+ \left\lVert \Pizko \left(\phi \, \frac{\partial \cn}{\partial t} - \phi \, \Pizko \left(\frac{\partial \cn}{\partial t}\right)  \right) \right\rVert_{0,\Omega} \\
& \qquad + \left\lVert \Pizko \left( \phi \, \Pizko \left(\frac{\partial \cn}{\partial t} - \frac{\Pc \cn - \Pc \cno}{\tau} \right)\right) \right\rVert_{0,\Omega} \\
&\le \eta \bigg[ h^{k+2} \left(\left| \phi \frac{\partial \cn}{\partial t} \right|_{k+2,\taun} + \left|\frac{\partial \cn}{\partial t} \right|_{k+2,\taun} \right) + \left\lVert \frac{\partial \cn}{\partial t} - \frac{\Pc \cn - \Pc \cno}{\tau} \right\rVert_{0,\Omega} \bigg].
\end{split}
\end{equation*}
Next, we bound $R_1^B$ with similar tools as for $R_1^A$:
\begin{equation*}
\begin{split}
R_1^B &\le \left\lVert (I-\Pizko) \left( \frac{\Pc \cn - \Pc \cno}{\tau} - \frac{\partial \cn}{\partial t} \right) \right\rVert_{0,\Omega} + \left\lVert (I-\Pizko) \frac{\partial \cn}{\partial t} \right\rVert_{0,\Omega} \\
&\quad \le \left\lVert \frac{\Pc \cn - \Pc \cno}{\tau} - \frac{\partial \cn}{\partial t} \right\rVert_{0,\Omega} + \eta \, h^{k+2} \left| \frac{\partial \cn}{\partial t} \right|_{k+2,\taun}.
\end{split}
\end{equation*}
Thus, we deduce with Lemma~\ref{lem:partial_cn}
\begin{equation} \label{eq:def_R1}
\begin{split}
R_1 &\le \eta \bigg[ h^{k+2} \left(\left| \phi \frac{\partial \cn}{\partial t} \right|_{k+2,\taun} + \left|\frac{\partial \cn}{\partial t} \right|_{k+2,\taun} \right) +\tau^{-\frac{1}{2}} h^{k+2} \left(\int_{\tnmo}^{\tn} \xi_{0,t}^2 \, \ds \right)^{\frac{1}{2}}  \\
&\qquad + \tau^{\frac{1}{2}} \left\lVert \frac{\partial^2 c}{\partial s^2} \right\rVert_{L^2(\tnmo,\tn;L^2(\Omega))} \bigg] \lVert \varthetan \rVert_{0,\Omega} \\
&=:\bigg[ h^{k+2} R_1^{n,1} + \tau^{-\frac{1}{2}} h^{k+2} R_1^{n,2} + \tau^{\frac{1}{2}} R_1^{n,3} \bigg] \lVert\varthetan \rVert_{0,\Omega},
\end{split}
\end{equation}
with the obvious definitions for the regularity terms $R_1^{n,1}$, $R_1^{n,2}$, and $R_1^{n,3}$.
\begin{comment}
\begin{equation} \label{eq:termsR1}
\begin{split}
R_1^{n,1}&:=\eta\left(\left| \phi \frac{\partial \cn}{\partial t} \right|_{k+1,\taun} + \left|\frac{\partial \cn}{\partial t} \right|_{k+1,\taun}\right), \quad 
R_1^{n,2}:=\eta\left(\int_{\tnmo}^{\tn} \xi_{0,t}^2 \, \ds \right)^{\frac{1}{2}}, \\
R_1^{n,3}&:=\eta\left(\int_{\tnmo}^{\tn} \left\lVert \frac{\partial^2 c}{\partial s^2}(s)\right\rVert^2_{0,\Omega} \, \ds \right)^{\frac{1}{2}}.
\end{split}
\end{equation}.
\end{comment}
\item term $R_2$: By the definition of $\Thetah(\cdot;\cdot,\cdot)$ in~\eqref{eq:def_Theta_h}, the identity $\varthetan=\Cn-\Pc \cn$, and the fact that $\qplusn$ and $\qminusn$ are non-negative, it holds
\begin{equation*}
\begin{split}
&\Thetah(\bun;\Pc \cn,\varthetan) - \Thetah(\bUno;\Cn,\varthetan)\\
& = \frac{1}{2} \left[ \left(\bun \cdot \nabla \Pc \cn,\varthetan \right)_h - \left(\bUno \cdot \nabla \Cn,\varthetan \right)_h \right]
-\frac{1}{2} \left( (\qplus+\qminus) \varthetan,\varthetan \right)_{0,\Omega }\\
& - \frac{1}{2} \left[ \left( \bun \Pc \cn,\nabla \varthetan \right)_h - \left( \bUno \Cn,\nabla \varthetan \right)_h \right] \\
& \le \frac{1}{2} \left[ \left(\bun \cdot \nabla \Pc \cn,\varthetan \right)_h - \left(\bUno \cdot \nabla \Cn,\varthetan \right)_h - \left( \bun \Pc \cn,\nabla \varthetan \right)_h + \left( \bUno \Cn,\nabla \varthetan \right)_h \right] .
\end{split}
\end{equation*}
The above equation, after adding zero in the form
\begin{equation*}
\begin{split}
0&=(\bUno \cdot \nabla \varthetan,\varthetan)_{h}
-(\bUno \cdot \nabla \varthetan,\varthetan)_{h}\\
&=(\bUno \cdot \nabla \Cn,\varthetan)_h
-(\bUno \cdot \nabla \Pc \cn,\varthetan)_h
-(\bUno \cdot \nabla \varthetan,\Cn)_h
+(\bUno \cdot \nabla \varthetan,\Pc \cn)_h
\end{split}
\end{equation*}
to the right hand side, can be equivalently expressed as
\begin{equation*}
\begin{split}
&\Thetah(\bun;\Pc \cn,\varthetan) - \Thetah(\bUno;\Cn,\varthetan)\\
& \le \frac{1}{2} \left[ \left((\bun-\bUno) \cdot \nabla \Pc \cn,\varthetan \right)_h - \left( (\bun-\bUno) \Pc \cn,\nabla \varthetan \right)_h\right] =:R_2^A + R_2^B.
\end{split}
\end{equation*}
For $R_2^A$, we estimate
\begin{equation*}
\begin{split}
R_2^A=\frac{1}{2} \left((\bun-\bUno) \nabla \Pc \cn,\varthetan \right)_h \le \frac{1}{2} \lVert \bun-\bUno \rVert_{0,\Omega} \lVert \Pizbk \nabla \Pc \cn \rVert_{\infty,\Omega} \lVert \varthetan \rVert_{0,\Omega}.
\end{split}
\end{equation*}
We now use an inverse estimate~\cite[Lemma 4.5.3]{BrennerScott}, the continuity of $\PizbkE$, a triangle inequality, the assumption that $\taun$ is quasi-regular, and Lemma~\ref{lem:cn_Pcn}, to deduce, for every $\E \in \taun$,
\begin{equation} \label{eq:bound_Pi_nabla_Pc}
\begin{split}
&\lVert \PizbkE \nabla \Pc \cn \rVert_{\infty,\E} \le \eta \, \hE^{-1} \lVert \PizbkE \nabla \Pc \cn \rVert_{0,\E} \le \eta \, \hE^{-1} \lVert \nabla \Pc \cn \rVert_{0,\E} \\
& \quad \le \eta \, \hE^{-1} \left( \lVert \nabla \Pc \cn - \nabla \cn \rVert_{0,\E} + \lVert \nabla \cn \rVert_{0,\E} \right) \\
&\quad \le \eta \, \left(h^{-1} \lVert \nabla \Pc \cn - \nabla \cn \rVert_{0,\taun} + \lVert \nabla \cn \rVert_{\infty,\E} \right) 
% \\ &\quad \le \eta \, \left( \xi_1^n + \lVert \nabla \cn \rVert_{\infty,\E} \right) 
\le \eta, 
\end{split}
\end{equation}
Recalling Lemma~\ref{lem:un_Uno}, the definitions of $\varthetano$ and $\rhono$, and Lemma~\ref{lem:cn_Pcn}, we get
\begin{equation} \label{eq:bound_bun_bUno}
\begin{split}
&\lVert \bun - \bUno \rVert_{0,\Omega} \\
&\quad \le \tau \left\lVert {\partial \bu}/{\partial t} \right\rVert_{L^{\infty}(\tnmo,\tn;L^2(\Omega))} + (\lVert \varthetano \rVert_{0,\Omega} + \lVert \rhono \rVert_{0,\Omega}) \zeta_1^{n-1} + h^{k+1} \, \zeta_2^{n-1} \\
&\quad \le \tau \left\lVert {\partial \bu}/{\partial t} \right\rVert_{L^{\infty}(\tnmo,\tn;L^2(\Omega))} + (\lVert \varthetano \rVert_{0,\Omega} + h^{k+2} \xi_0^{n-1}) \zeta_1^{n-1} + h^{k+1} \, \zeta_2^{n-1},
\end{split}
\end{equation}
thus implying
\begin{equation*}
R_2^A \le \eta \bigg[h^{k+1} R_2^{n,1} + \tau R_2^{n,2} + \lVert \varthetano \rVert_{0,\Omega} R_2^{n,3}\bigg] \lVert \varthetan \rVert_{0,\Omega},
\end{equation*}
with the obvious definitions for the regularity terms $R_2^{n,1}$, $R_2^{n,2}$, and $R_2^{n,3}$.
\begin{comment}
\begin{equation} \label{eq:termsR2}
R_2^{n,1}:=h \xi_0^{n-1} \zeta_1^{n-1} + \zeta_2^{n-1}, \quad 
R_2^{n,2}:=\left\lVert \frac{\partial \bu}{\partial t} \right\rVert_{L^{\infty}(\tnmo,\tn;L^2(\Omega))}, \quad
R_2^{n,3}:=\zeta_1^{n-1}.
\end{equation}
\end{comment}

The term $R_2^B$ can be bounded analogously to $R_2^A$, giving
\begin{equation*}
R_2^B=\frac{1}{2} \left( (\bUno-\bun) \Pc \cn,\nabla \varthetan \right)_h 
\le \eta \, \lVert \bun-\bUno \rVert_{0,\Omega} |\varthetan|_{1,\taun} .
\end{equation*}
Using again the bound~\eqref{eq:bound_bun_bUno}, one obtains
\begin{equation*}
R_2^B \le \eta \bigg[h^{k+1} R_2^{n,1} + \tau R_2^{n,2} + \lVert \varthetano \rVert_{0,\Omega} R_2^{n,3}\bigg] |\varthetan|_{1,\taun}.
\end{equation*}
Thus,
\begin{equation*}
R_2 \le \eta \bigg[h^{k+1} R_2^{n,1} + \tau R_2^{n,2} + \lVert \varthetano \rVert_{0,\Omega} R_2^{n,3}\bigg] \left( \lVert \varthetan \rVert_{0,\Omega} + |\varthetan|_{1,\taun} \right).
\end{equation*}
\item term $R_3$:  We use the definition of $\Dcalh(\cdot;\cdot,\cdot)$ in~\eqref{eq:def_Dcalh}, a standard H\"{o}lder inequality in the spirit of~\eqref{eq:cont_D_2}, the estimate~\eqref{eq:bound_Pi_nabla_Pc}, the scaling properties of the stabilization in~\eqref{eq:stab_assumpt}, and the Lipschitz continuity of $D(\cdot;\cdot,\cdot)$ and $\nuDE$ in~\eqref{eq:stab_constants}, to deduce 
\begin{equation*}
\begin{split}
R_3 &= \Dcalh\left(\bun;\Pc \cn, \varthetan \right) - \Dcalh\left(\bUno;\Pc \cn, \varthetan \right) \\
&= ( (D(\Pizbk \bun)-D(\Pizbk\bUno)) \, \Pizbk(\nabla \Pc \cn) \cdot  \Pizbk(\nabla \varthetan) )_{0,\Omega} \\
&\quad + \sum_{\E \in \taun} (\nuDE(\bun)-\nuDE(\bUno)) \SED\left((I-\PinablaE) \Pc \cn,(I-\PinablaE) \varthetan) \right)  \\
&\le \eta \lVert \bun-\bUno \rVert_{0,\Omega} |\varthetan|_{1,\taun}.
\end{split}
\end{equation*}
Hence, with~\eqref{eq:bound_bun_bUno} we have
\begin{equation*} 
R_3 \le \eta \bigg[h^{k+1} R_2^{n,1} + \tau R_2^{n,2} + \lVert \varthetano \rVert_{0,\Omega} R_2^{n,3} \bigg] |\varthetan|_{1,\taun}.
\end{equation*}
\item term $R_4$: The use of Lemma~\ref{lem:cn_Pcn} yields 
\begin{equation*}
\begin{split}
R_4 &= - [(\cn,\varthetan)_{0,\Omega} - (\Pc \cn,\varthetan)_h] 
= -[((I-\Pizko) \cn,\varthetan)_{0,\Omega} + (\Pizko (\cn-\Pc \cn), \varthetan)_{0,\Omega}] \\
&\le \eta \, h^{k+2} \left[|\cn|_{k+2,\taun} + \xi_0^n \right] \lVert \varthetan \rVert_{0,\Omega} =:\eta \, h^{k+2} R_4^{n,1} \lVert \varthetan \rVert_{0,\Omega},
\end{split}
\end{equation*}
with the obvious definition of $R_4^{n,1}$.
\item term $R_5$: The approximation properties in Lemma~\ref{lem:approx_properties} yield
\begin{equation*}
\begin{split}
R_5 = - \left((I-\Pizko) (\qplusn \chatn), \varthetan \right)_{0,\Omega}  \le \eta \, h^{k+2} \, |\qplusn\chatn|_{k+2,\taun} \lVert\varthetan \rVert_{0,\Omega} =:\eta \, h^{k+2} R_5^{n,1} \lVert\varthetan \rVert_{0,\Omega},
\end{split}
\end{equation*}
with the obvious definition of $R_5^{n,1}$.
\end{itemize}
We now insert~\eqref{eq:estimate:Dh} and the bounds on $R_1$-$R_5$ into~\eqref{eq:rec_rel_theta}. Afterwards, we observe that all regularity terms $\{R_J^{n,i}\}$ above only depend on the continuous solution and can be assumed to be bounded uniformly in $h$. We only keep track of the terms $R_1^{n,2}$ and $R_1^{n,3}$. This yields
\begin{equation} \label{eq:ev_eq_proof}
\begin{split}
&\frac{1}{\tau} \Mcalh\left(\varthetan-\varthetano,\varthetan \right) + D_\ast \left| \varthetan \right|^2_{1,\taun} \\
&\le \lVert \varthetano \rVert_{0,\Omega} \lVert \varthetan \rVert_{0,\Omega} \,\omega^n_1 + \lVert \varthetano \rVert_{0,\Omega} |\varthetan|_{1,\taun} \omega^n_2 + \lVert \varthetan \rVert_{0,\Omega} \omega^n_3 + |\varthetan|_{1,\taun} \omega^n_4 \\
&= \lVert \varthetan \rVert_{0,\Omega} \left[\omega^n_3 + \lVert \varthetano \rVert_{0,\Omega} \,\omega^n_1 \right] +
|\varthetan|_{1,\taun} \left[ \omega^n_4 + \lVert \varthetano \rVert_{0,\Omega} \,\omega^n_2 \right],
\end{split}
\end{equation}
with the positive scalars
\begin{equation} \label{eq:omega_i}
\begin{split}
\omega_i^n \le \eta, \quad i=1,2, \quad
\omega^n_3\le \eta \left( \tau + h^{k+1} + \tau^{-\frac{1}{2}} h^{k+2} R_1^{n,2} + \tau^{\frac{1}{2}} R_1^{n,3} \right), \quad 
\omega^n_4\le \eta \left(\tau + h^{k+1} \right).
\end{split}
\end{equation}
Next, we introduce, for all $w_h \in \Zh$, the discrete norm
\begin{equation} \label{eq:discr_norm}
\lVert w_h \rVert^2_{0,h}:=\Mcalh(w_h,w_h).
\end{equation}
Owing to Lemma~\ref{lem:properties_D_A}, there exist positive constants $c_\ast$ and $c^\ast$, such that, for all $w_h \in \Zh$, it holds
\begin{equation} \label{eq:discr_norm_equiv}
c_\ast \lVert w_h \rVert_{0,h} \le \lVert w_h \rVert_{0,\Omega} \le c^\ast \lVert w_h \rVert_{0,h}.
\end{equation}
Reshaping~\eqref{eq:ev_eq_proof}, and employing~\eqref{eq:discr_norm} and~\eqref{eq:discr_norm_equiv}, then gives
\begin{equation} \label{eq:ev_eq_proof_2}
\begin{split}
&\lVert \varthetan \rVert^2_{0,h} + \tau D_\ast \left| \varthetan \right|^2_{1,\taun} \\
&\le \Mcalh(\varthetano,\varthetan) 
+ \tau \lVert \varthetan \rVert_{0,h} \left[c^\ast \omega^n_3 + \lVert \varthetano \rVert_{0,h} (c^\ast)^2 \omega^n_1 \right] 
+ \tau |\varthetan|_{1,\taun} \left[ \omega^n_4 + \lVert \varthetano \rVert_{0,h} \, c^\ast \omega^n_2 \right]\\
&=:T_1+T_2+T_3
\end{split}
\end{equation}
The terms $T_1$ and $T_2$ are bounded as follows: 
\begin{equation} \label{eq:T1_T2}
\begin{split}
T_1+T_2 &\le
\lVert \varthetan \rVert_{0,h} \left[ (1+\tau \eta) \lVert \varthetano \rVert_{0,h} + \tau c^\ast \omega^n_3 \right] \\
&\le \frac{1}{2} \left( \lVert \varthetan \rVert^2_{0,h} + \left[ (1+\tau \eta) \lVert \varthetano \rVert_{0,h} + \tau c^\ast \omega^n_3 \right]^2 \right),
\end{split}
\end{equation}
where we used~\eqref{eq:splitting_Mh} and~\eqref{eq:discr_norm} in the first step.
%and an inequality of the form $ab \le \frac{1}{2}(a^2+b^2)$, for all $a,b\ge 0$, in the second. 
The term $T_3$ is bounded as follows
%by applying an inequality of the type $ab \le \varepsilon a^2 + \frac{1}{4\varepsilon} b^2$, for all $a,b\ge 0$, and suitable $\varepsilon>0$:
\begin{equation} \label{eq:T3}
\begin{split}
T_3 &\le \tau D_\ast |\varthetan|^2_{1,\taun} + \frac{\tau}{4 D_\ast} \left[ \omega^n_4 + \lVert \varthetano \rVert_{0,h} \, c^\ast \omega^n_2 \right]^2  \\
&\le \tau D_\ast |\varthetan|^2_{1,\taun} + \frac{\tau}{2} \eta \left[ (\omega^n_4)^2 + \lVert \varthetano \rVert^2_{0,h} \right]^2.
\end{split}
\end{equation}
Next, we plug~\eqref{eq:T1_T2} and~\eqref{eq:T3} into~\eqref{eq:ev_eq_proof_2}, cancel the terms $\tau D_\ast |\varthetan|^2_{1,\taun}$ and manipulate the resulting inequality, to obtain
\begin{equation*}
\begin{split}
\lVert \varthetan \rVert^2_{0,h} \le \left[ (1+\tau \eta) \lVert \varthetano \rVert_{0,h} + \tau c^\ast \omega^n_3 \right]^2 + \tau \eta \left[ (\omega^n_4)^2 + \lVert \varthetano \rVert^2_{0,h} \right]^2.
\end{split}
\end{equation*}
Moreover, we estimate
%, using again $ab\le \frac{1}{2}(a^2+b^2)$, for all $a,b\ge 0$,
\begin{equation*}
\begin{split}
&\left[ (1+\tau \eta) \lVert \varthetano \rVert_{0,h} + \tau c^\ast \omega^n_3 \right]^2 \\
&\quad=(1+\tau \eta)^2 \lVert \varthetano \rVert^2_{0,h} 
+2 \tau^{\frac{1}{2}} \lVert \varthetano \rVert_{0,h} \tau^{\frac{1}{2}} (1+\tau \eta) c^\ast \omega^n_3 + \tau^2 (c^\ast)^2 (\omega^n_3)^2 \\
&\quad \le \left[ (1+\tau \eta)^2 + \tau \right] \lVert \varthetano \rVert^2_{0,h} + \left[ \tau(1+\tau \eta)^2 + \tau^2 \right] (c^\ast)^2 (\omega^n_3)^2 \\
&\quad \le \left( 1+\tau \eta \right) \lVert \varthetano \rVert^2_{0,h} + \tau \eta (\omega^n_3)^2.
\end{split}
\end{equation*}
Hence,
\begin{equation*}
\lVert \varthetan \rVert^2_{0,h} \le (1+\tau \eta) \lVert \varthetano \rVert^2_{0,h} + \tau \eta \left[ (\omega^n_3)^2 + (\omega^n_4)^2 \right].
\end{equation*}
Defining
\begin{equation*}
\gamma^n:=(\omega^n_3)^2 + (\omega^n_4)^2
\end{equation*}
and solving the recursion then leads to
\begin{equation*}
\begin{split}
\lVert \varthetan \rVert^2_{0,h} \le (1+\tau \eta)^n \lVert \varthetao \rVert^2_{0,h} + \tau \eta \sum_{j=1}^n \gamma^j
\le \eta \lVert \varthetao \rVert^2_{0,h} + \tau \eta \sum_{j=1}^n \gamma^j,
\end{split}
\end{equation*}
where we recall that $n \le T/\tau$ with $T$ the final time instant. With~\eqref{eq:discr_norm_equiv} the estimate in the $L^2$ norm is a direct consequence:
\begin{equation} \label{eq:thetan_est1}
\begin{split}
\lVert \varthetan \rVert^2_{0,\Omega} \le \eta \lVert \varthetao \rVert^2_{0,\Omega} + \tau \eta \sum_{j=1}^n \gamma^j.
\end{split}
\end{equation}
The initial term $\lVert \varthetao \rVert_{0,\Omega}^2$ is estimated by
\begin{equation}\label{eq:L3}
\lVert \varthetao \rVert_{0,\Omega} = \lVert \czh - \Pc c^0 \rVert_{0,\Omega} \le \lVert \czh - c^0 \rVert_{0,\Omega} + \lVert c^0 - \Pc c^0 \rVert_{0,\Omega} \le \lVert \czh - c^0 \rVert_{0,\Omega} + h^{k+2} \, \xi_0^0,
\end{equation}
where we applied Lemma~\ref{lem:cn_Pcn}. Moreover, using~\eqref{eq:omega_i}, the fact that $\sum_{j=1}^n \tau \le T$, and the definitions of $R_1^{j,2}$ and $R_1^{j,3}$ in~\eqref{eq:def_R1}, after some simple manipulations, we obtain
\begin{equation} \label{eq:thetan_est2}
\begin{split}
\tau \eta \sum_{j=1}^n \gamma_j &\le \eta \left( \sum_{j=1}^n \tau (\omega^j_3)^2 + \sum_{j=1}^n \tau (\omega_4^j)^2 \right) \\
&\le \eta \left[ \sum_{j=1}^n \tau (\tau+h^{k+1})^2 + (h^{k+2})^2 \sum_{j=1}^n (R_1^{j,2})^2 + \tau^2 \sum_{j=1}^{n} (R_1^{j,3})^2 \right] \\
&\le \eta \left[ (\tau+h^{k+1})^2 + (h^{k+2})^2 \sum_{j=1}^n (R_1^{j,2})^2 + \tau^2 \sum_{j=1}^{n} (R_1^{j,3})^2 \right] \\
&\le \eta \left[ (\tau+h^{k+1})^2 + (h^{k+2})^2 \int_{0}^{\tn} \xi_{0,t}^2 \, \ds + \tau^2 \int_{0}^{\tn} \left\lVert \frac{\partial^2 c}{\partial s^2}(s)\right\rVert^2_{0,\Omega} \, \ds \right].
\end{split}
\end{equation}
The assertion of the theorem follows by combining~\eqref{eq:thetan_est1} with~\eqref{eq:thetan_est2} and \eqref{eq:L3}.
\end{proof}

\begin{comment}
\begin{remark}
We underline that, in contrast to Lemma~\ref{lem:cn_Pcn}, the higher regularity assumption on $\cn$, namely $\cn \in H^{k+2}(\taun)$ instead of $\cn \in H^{k+1}(\taun)$ (see also Remark~\ref{rem:higher_reg_c}) needed in the proof of Theorem~\ref{thm:Cn_cn}, does not lead to a convergence rate of order $k+2$ in the $L^2$ norm for the dicretization error of the concentration. Although, at a first glance, this could seem to be a bit odd, this can in fact be seen as a consequence of the nonlinear coupling of the equations.
\end{remark}
\end{comment}

\section{Numerical experiments}\label{sec:5}

In this section, we demonstrate the performance of the method on the basis of numerical experiments, focusing on the lowest order case $k=0$. To this purpose, we first consider an ideal test case (\textit{Example 1}), and then a more realistic one (\textit{Example 2}). The aim of the first test is to validate (also numerically) the convergence of the method on a problem with regular known solution, whereas those of the second test is to check the method's performance on a well-known benchmark that mimics a more realistic situation. 

\begin{figure}[h]
\begin{minipage}[t]{0.485\textwidth}
\centering
\includegraphics[width=\textwidth]{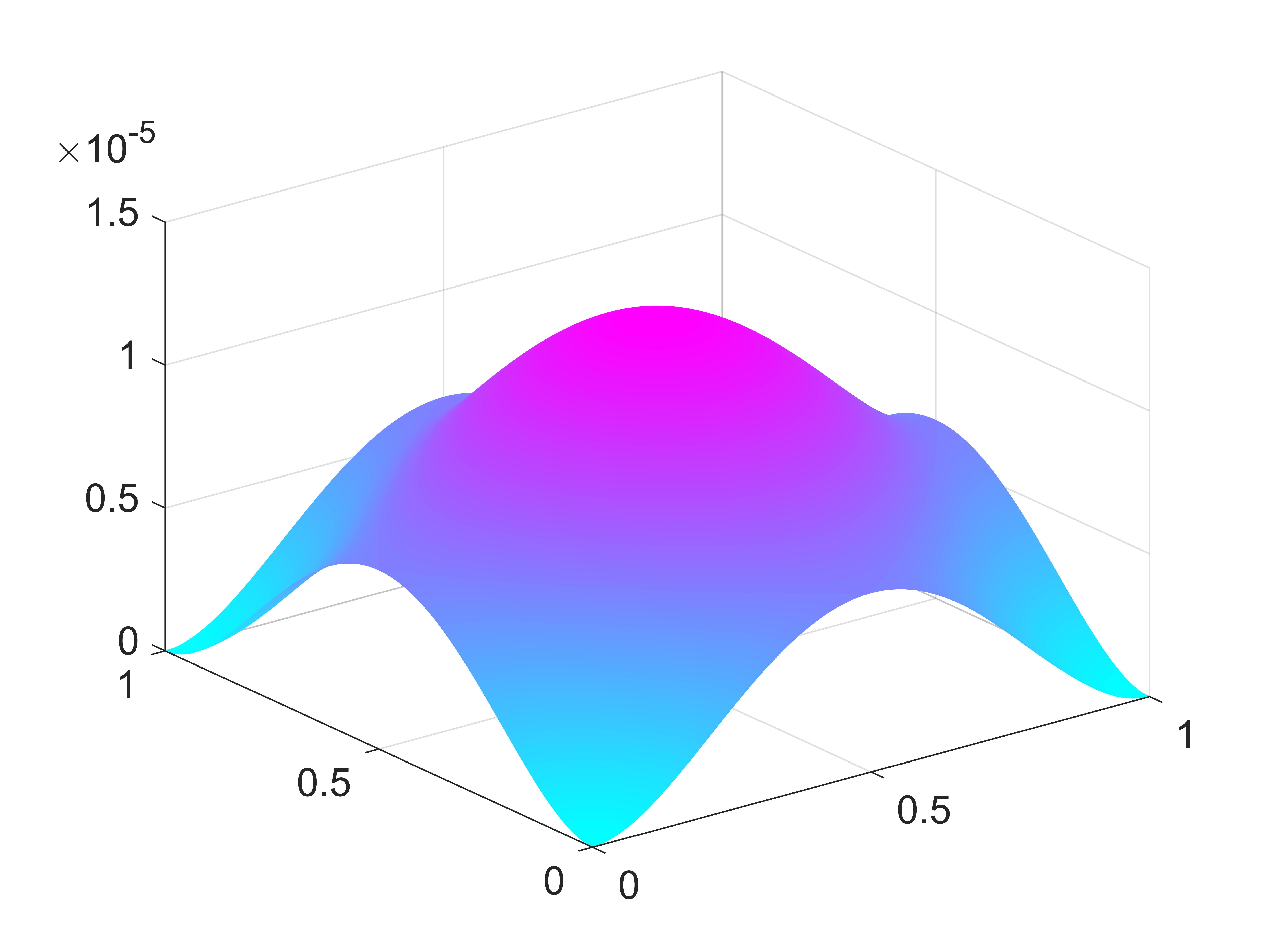}
\end{minipage}
\hfill
\begin{minipage}[t]{0.485\textwidth}
\centering
\includegraphics[width=\textwidth]{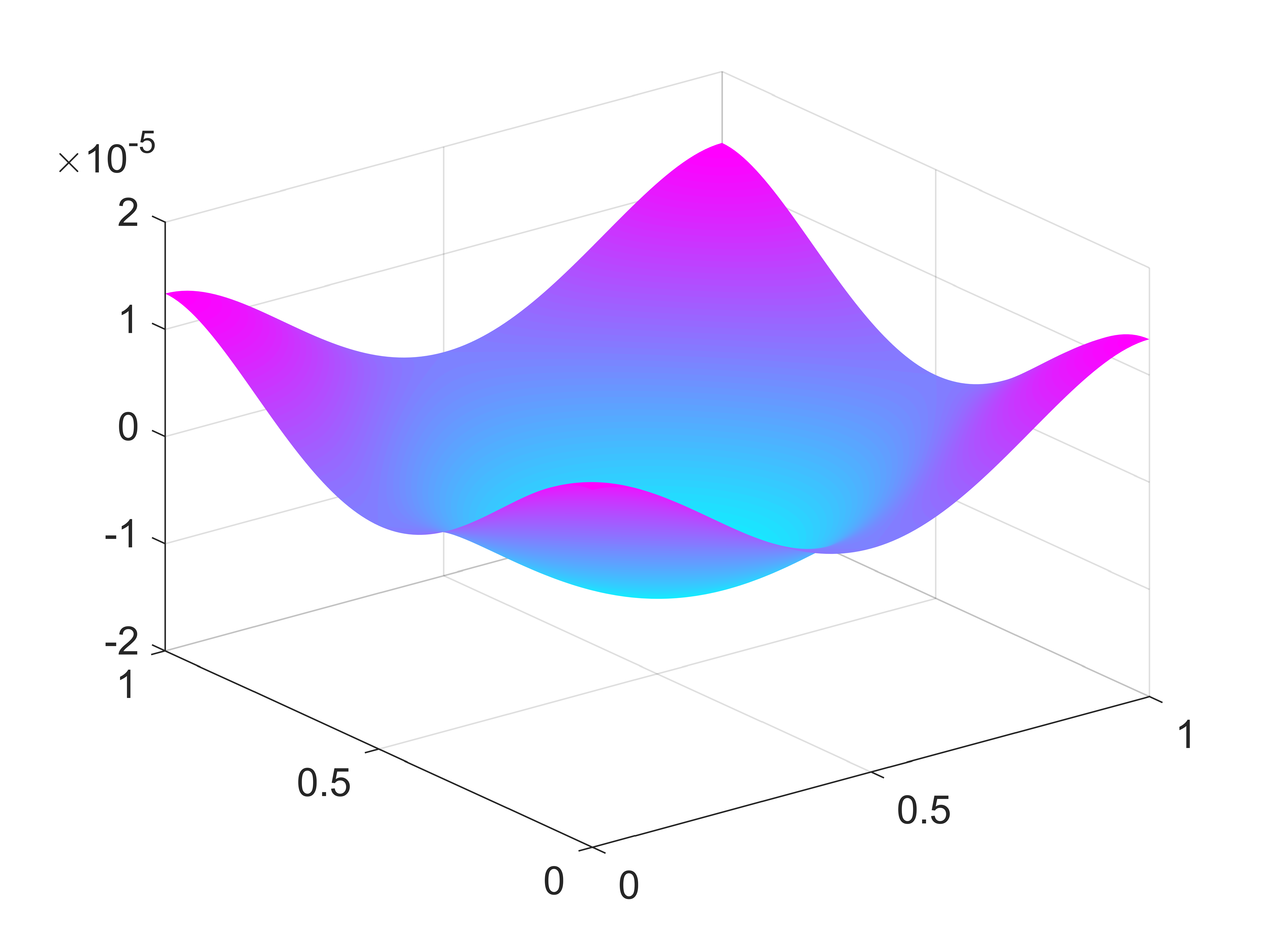}
\end{minipage}
\caption{Exact concentration $c$ (left) and pressure $p$ (right) of example 1, given by~\eqref{eq:Hu1_fcts}, at the final time $T=0.01$.}
\label{fig:Hu1_c_p}
\end{figure}

\begin{figure}[h]
\centering
\includegraphics[width=0.6\textwidth]{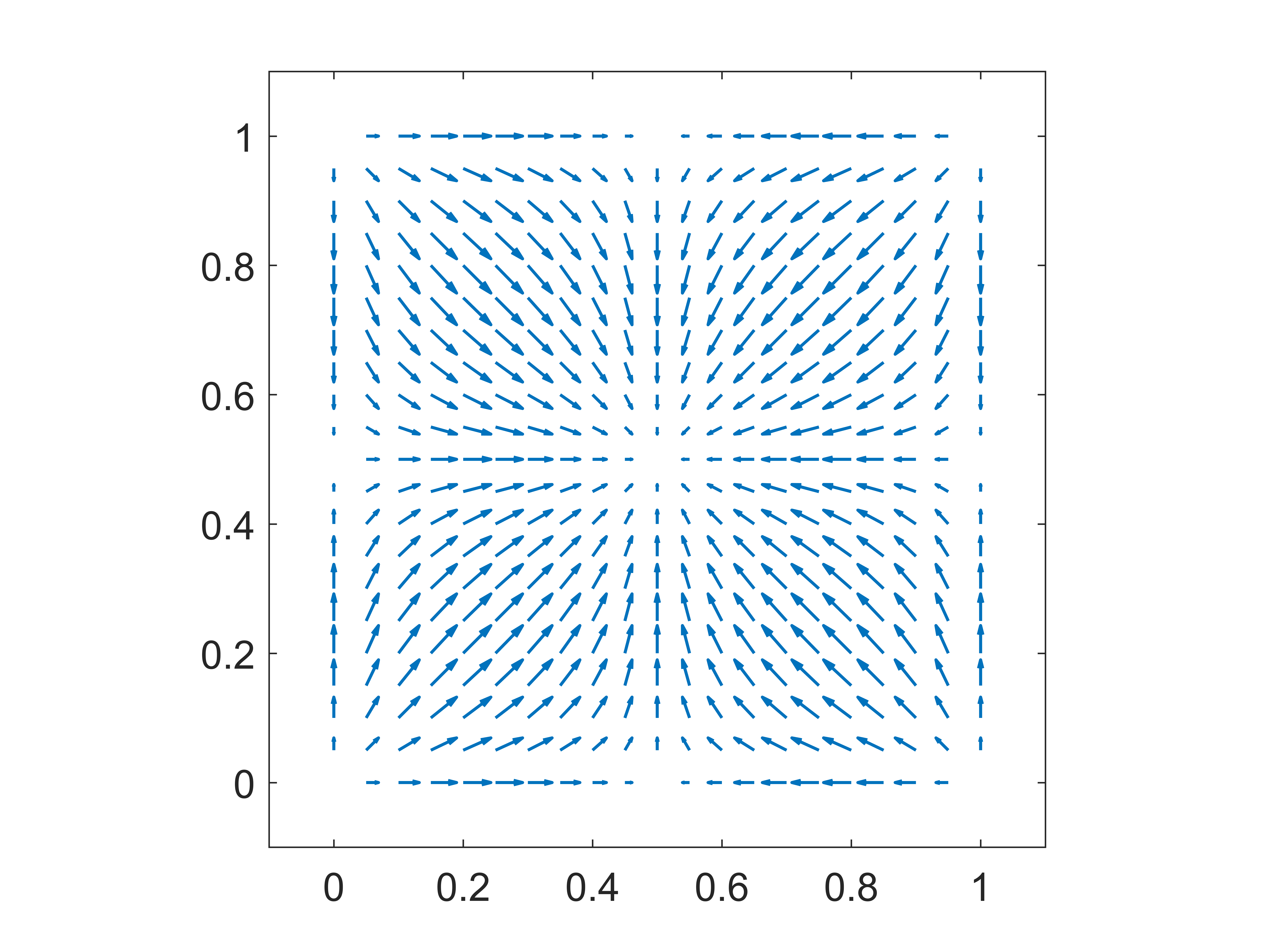}
\caption{Exact vector field $\bu$ of example 1, given by~\eqref{eq:Hu1_fcts}, at the final time $T=0.01$.}
\label{fig:Hu1_u}
\end{figure}

\vspace{0.2cm}
\textit{Example 1:} 
Here, we study a generalized version of~\eqref{eq:model_problem}, given by
\begin{equation*} \label{eq:generalized_problem}
\left\{
\begin{alignedat}{2}
\phi \, \frac{\partial c}{\partial t} + \bu \cdot \nabla c - \div (D(\bu) \nabla c) &= f \\
\div \, \bu &= g \\
\bu &= -a(c) (\nabla p - \bgamma(c)),
\end{alignedat}
\right.
\end{equation*}
endowed with the boundary and initial conditions in~\eqref{eq:bdry_cond} and~\eqref{eq:initial_cond}, respectively.
We fix $\Omega=(0,1)^2$ and pick the same choice of parameters as in~\cite{miscibledispl2018hu}, namely
$T=0.01$, $\phi=1$, $D(\bu)=|\bu|+0.02$, $d_m=0.02$, $d_{\ell}=d_t=1$, $c_0=0$, $\bgamma(c)=0$, and $a(c)=(c+2)^{-1}$, where $f$ and $g$ are taken in accordance with the analytical solutions
\begin{equation} \label{eq:Hu1_fcts}
\begin{split}
c(x,y,t)&=t^2 \left[x^2(x-1)^2 + y^2(y-1)^2\right] \\
\bu(x,y,t)&=2t^2 \begin{pmatrix} x(x-1)(2x-1) \\ y(y-1)(2y-1) \end{pmatrix} \\
p(x,y,t)&=-\frac{1}{2} c^2 -2c + \frac{17}{6300} t^4 + \frac{2}{15}t^2.
\end{split}
\end{equation}
Plots of the exact solution at the final time $T$ are shown in Figures~\ref{fig:Hu1_c_p} and~\ref{fig:Hu1_u}.

\begin{figure}[h]
\begin{minipage}[t]{0.485\textwidth}
\centering
\includegraphics[width=0.75\textwidth]{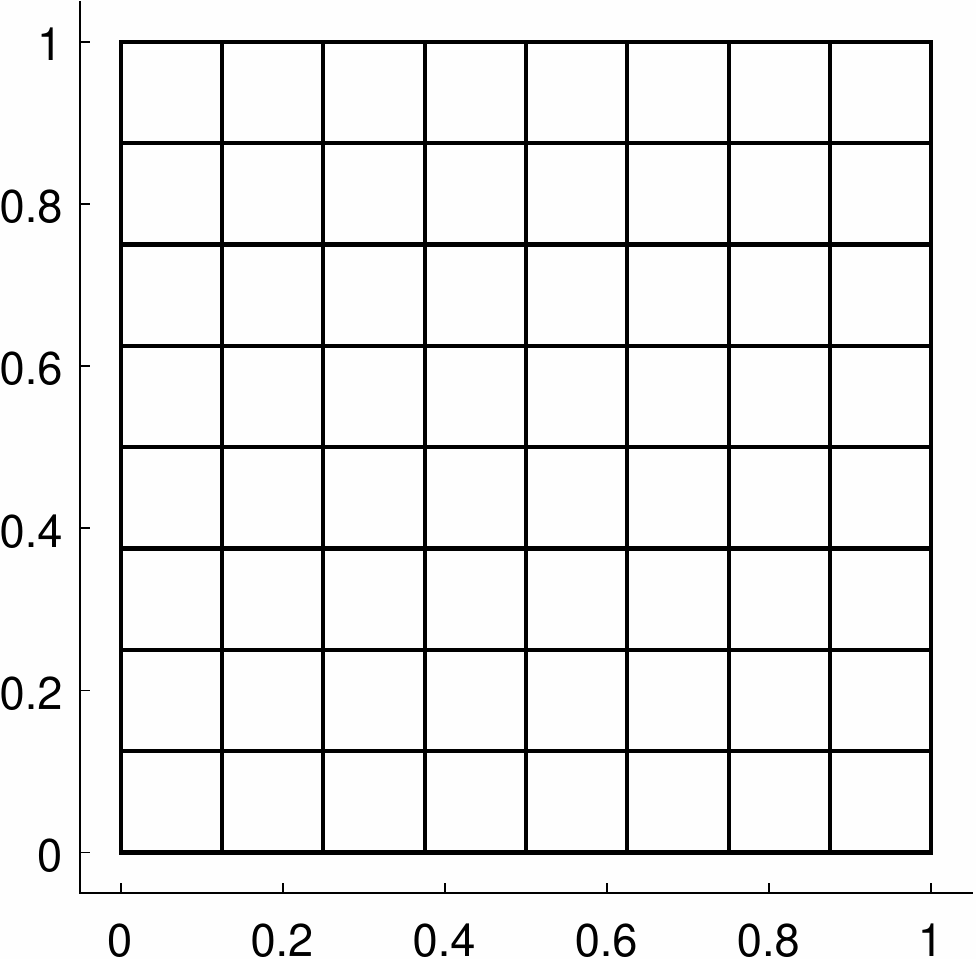}
\end{minipage}
\begin{minipage}[t]{0.485\textwidth}
\centering
\includegraphics[width=0.75\textwidth]{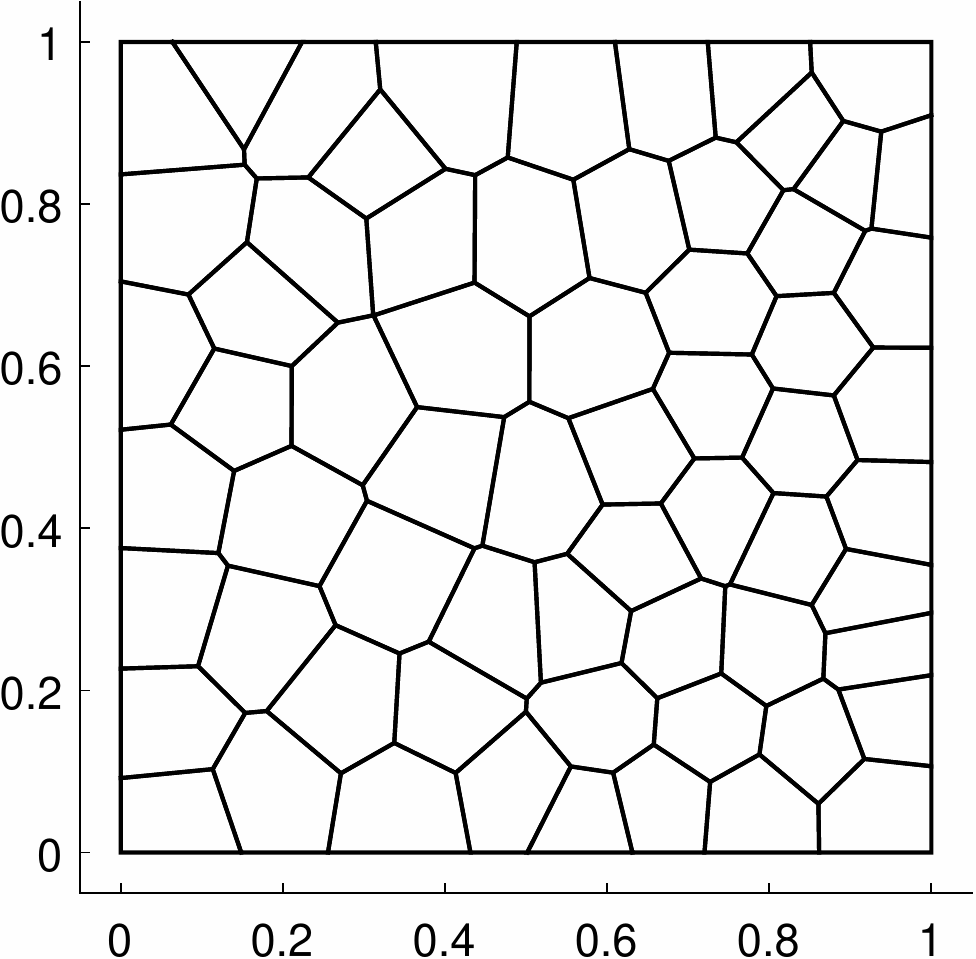}
\end{minipage}
\caption{Meshes: regular 8x8 Cartesian mesh (left); Voronoi mesh with 64 elements (right).}
\label{fig:meshes}
\end{figure}

We employ a sequence of regular Cartesian meshes and Voronoi meshes, as portrayed in Figure~\ref{fig:meshes}. In addition to the current version, we also test the method when replacing the stabilization terms in~\eqref{eq:def_Mh_loc},~\eqref{eq:def_Dcalh_loc}, and~\eqref{eq:def_Acalh_loc} by alternative ones:
\begin{equation*}
\begin{split}
\nuME(\phi) \SEM\left((I-\PizkoE) \ch, (I-\PizkoE) \zh \right) \quad &\rightsquigarrow \quad \widetilde{\SEM}\left((I-\PizkoE) \ch, (I-\PizkoE) \zh \right) \\
\nuDE(\buh) \, \SED\left((I-\PinablaE) \ch,(I-\PinablaE) \zh) \right) \quad &\rightsquigarrow \quad \widetilde{\SED}\left(\buh;(I-\PinablaE) \ch,(I-\PinablaE) \zh) \right)\\
\nuAE(\ch) \, \SEA((I-\PizbkE)\buh,(I-\PizbkE)\bvh) \quad &\rightsquigarrow \quad \widetilde{\SEA}(\ch;(I-\PizbkE)\buh,(I-\PizbkE)\bvh).
\end{split}
\end{equation*}
%Let $\{\philE\}_{\ell=1}^{\textrm{dim}\ZhE}$ and $\{\psilE\}_{\ell=1}^{\NVhE}$ be the local canonical basis functions for $\ZhE$ and $\VhE$ defined by duality, i.e.,
%\begin{equation*}
%% \dofZhE_j(\philE) = \delta_{j,\ell}, \qquad
%{\dof_j^{\ZhE}}(\philE) = \delta_{j,\ell}, \qquad
%\dofVhE_j(\psilE) = \delta_{j,\ell},
%\end{equation*}
%where $\delta_{j,\ell}$ here denotes the Kronecker delta. 
The alternative (diagonal) stabilizations are given by
\begin{equation} \label{eq:stabs_new}
\begin{split}
\widetilde{\SEM}(\ch,\zh)&=|\E| \sum_{j=1}^{\NZhE} d^{\mathcal{M}}_j \, \dofZhE_j(\ch) \, \dofZhE_j(\zh) \\
\widetilde{\SED}(\ch,\zh)&=\sum_{j=1}^{\NZhE} d^{D}_j \, \dofZhE_j(\ch) \, \dofZhE_j(\zh) \\
\widetilde{\SEA}(\buh,\bvh)&=|\E|\sum_{j=1}^{\NVhE} d^{\mathcal{A}}_j \, \dofVhE_j(\buh) \, \dofVhE_j(\bvh),
\end{split}
\end{equation}
with 
\begin{equation} \label{eq:stabs_new_coeff}
\begin{split}
d^{\mathcal{M}}_j \,&:=\max\left\{ \frac{1}{ |\E|} \,\int_\E \phi \, (\PizkoE \philE)^2 \, \dx,\, \sigma \nuME(\phi) \right\} \\
d^{D}_j \,&:=\max\left\{\int_\E D(\PizbkE \buh) \, |\PizbkE(\nabla \philE)|^2 \, \dx, \, \sigma \nuDE(\buh) \right\} \\
d^{\mathcal{A}}_j \,&:=\max\left\{\frac{1}{ |\E|} \, \int_\E A(\PizkoE \ch) \, |\PizbkE \psilE|^2 \, \dx, \, \sigma  \nuAE(\ch) \right\},
\end{split}
\end{equation} 
where $\{\philE\}_{\ell=1}^{\textrm{dim}\ZhE}$ and $\{\psilE\}_{\ell=1}^{\NVhE}$ denote the local canonical basis functions for $\ZhE$ and $\VhE$,
and $\sigma>0$ is a safety parameter. In the forthcoming experiments, we set $\sigma=1e-3$. We highlight that these stabilizations are in fact modifications of the so-called \textit{D-recipe}, which was introduced in~\cite{VEM3Dbasic} and has already been successfully applied in some variants to other model problems, such as the Helmholtz problem~\cite{TVEM_Helmholtz_num}. The first entry inside the max is simply the ``diagonal part'' of the consistency term of the local approximate forms in~\eqref{eq:def_Mh_loc},~\eqref{eq:def_Dcalh_loc}, and~\eqref{eq:def_Acalh_loc}, respectively, whereas the second terms correspond to the original stabilizations associated to the degrees of freedom in~\eqref{eq:stabs} multiplied by~$\sigma$, which acts as a positivity safeguard. Importantly, it is easy to check that the error analysis can be easily extended to the new choice of stabilizations.

Due to the virtuality of the basis functions, we measure the following relative $L^2$ errors:
\begin{equation*}
\frac{\lVert c-\Pi^0_1 \Cn \rVert_{0,\Omega}}{\lVert c \rVert_{0,\Omega}}, \quad
\frac{\lVert \bu-\boldsymbol{\Pi}^0_0 \bUn \rVert_{0,\Omega}}{\lVert \bUn \rVert_{0,\Omega}}, \quad
\frac{\lVert p-\Pi^0_0 \Pn \rVert_{0,\Omega}}{\lVert p \rVert_{0,\Omega}},
\end{equation*}
where $\Cn$, $\bUn$, and $\Pn$ are the numerical solutions at the final time $T$.

\begin{figure}[h]
\begin{minipage}[t]{0.485\textwidth}
\centering
\includegraphics[width=\textwidth]{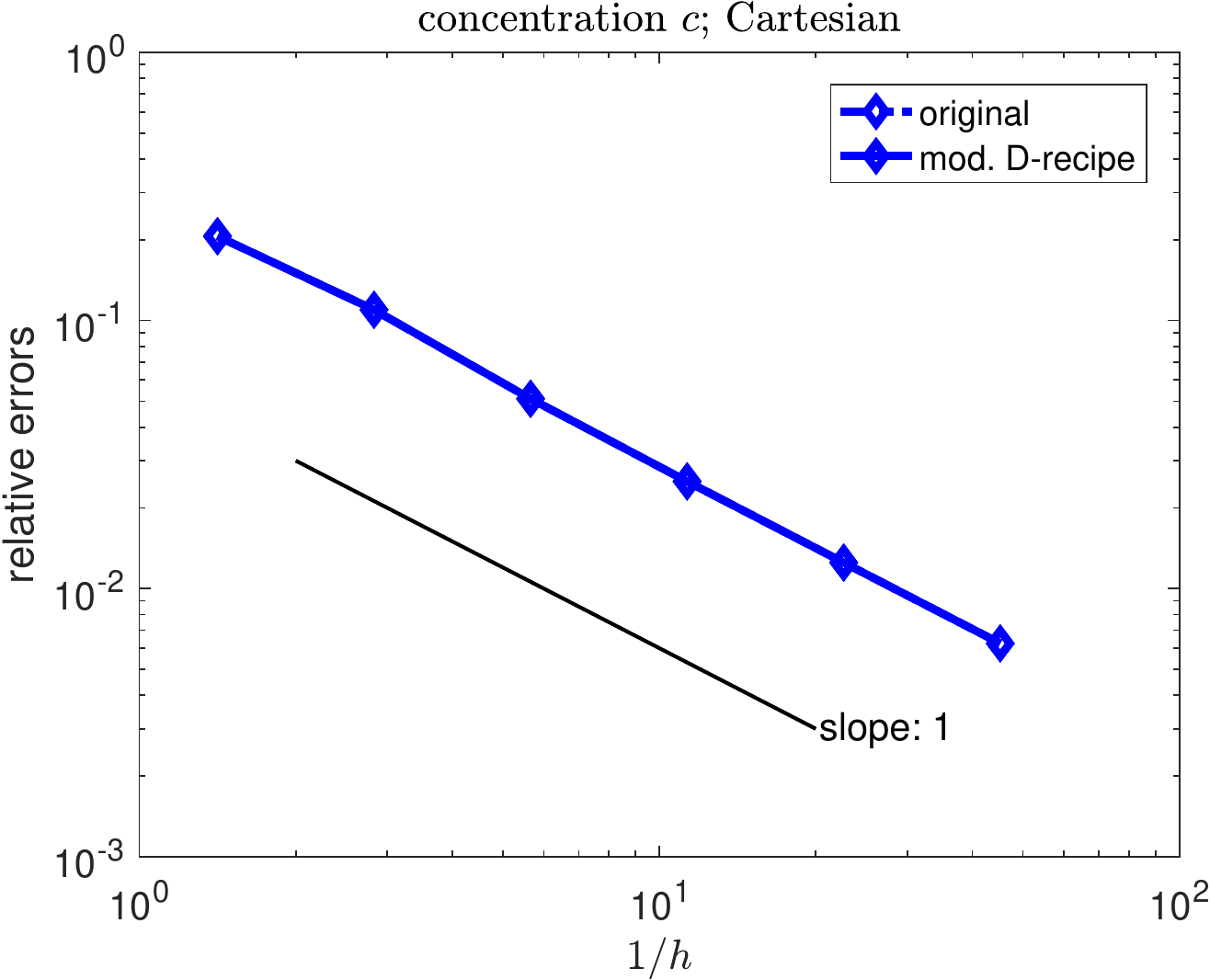}
\end{minipage}
\hfill
\begin{minipage}[t]{0.485\textwidth}
\centering
\includegraphics[width=\textwidth]{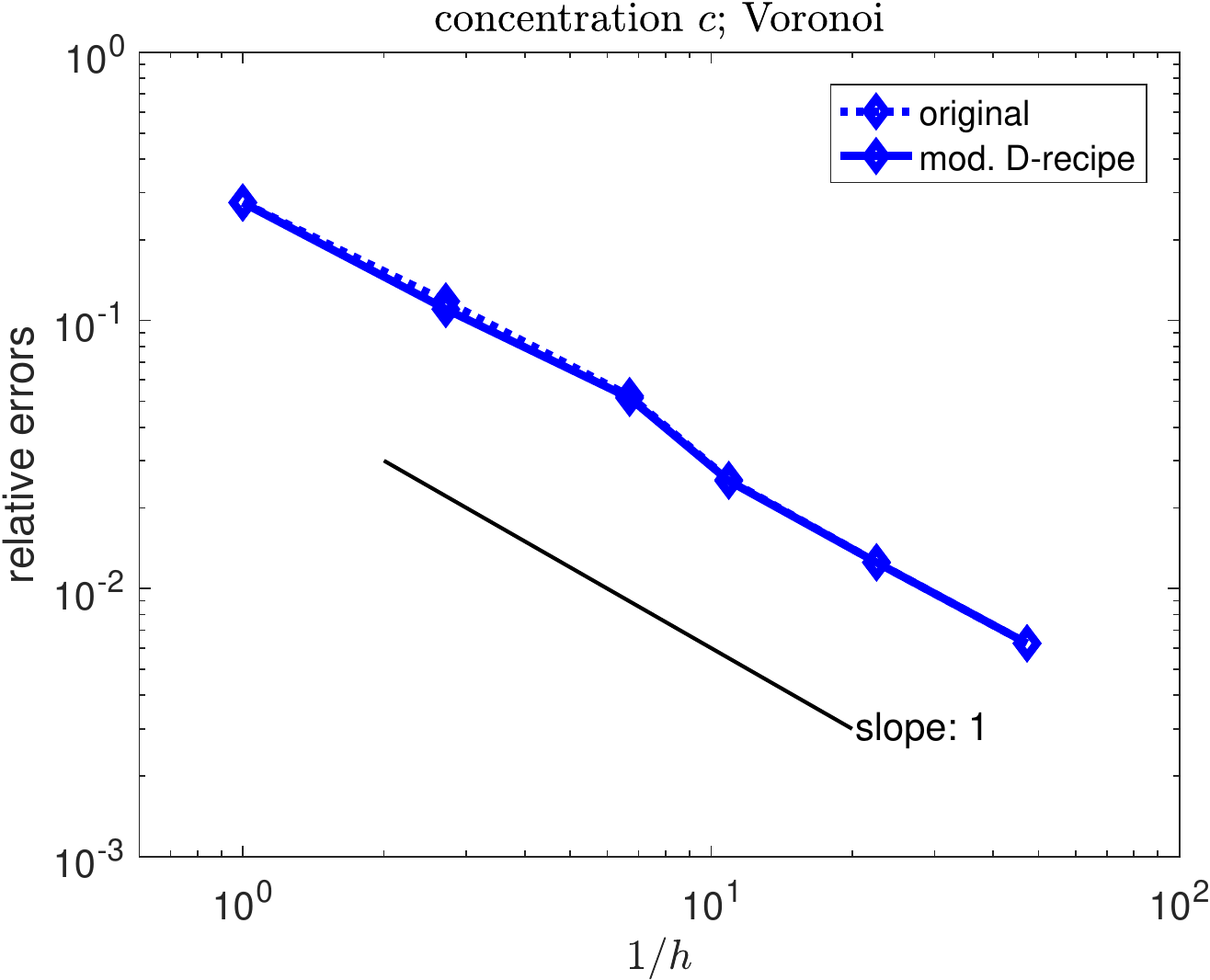}
\end{minipage}
\caption{Relative $L^2$ errors for the concentration in example~1 at the final time $T$ on regular Cartesian meshes (left) and Voronoi meshes (right). The original stabilization~\eqref{eq:stabs} and the D-recipe stabilization~\eqref{eq:stabs_new} are employed.}
\label{fig:Ex1:c}
\end{figure}

\begin{figure}[h]
\begin{minipage}[t]{0.485\textwidth}
\centering
\includegraphics[width=\textwidth]{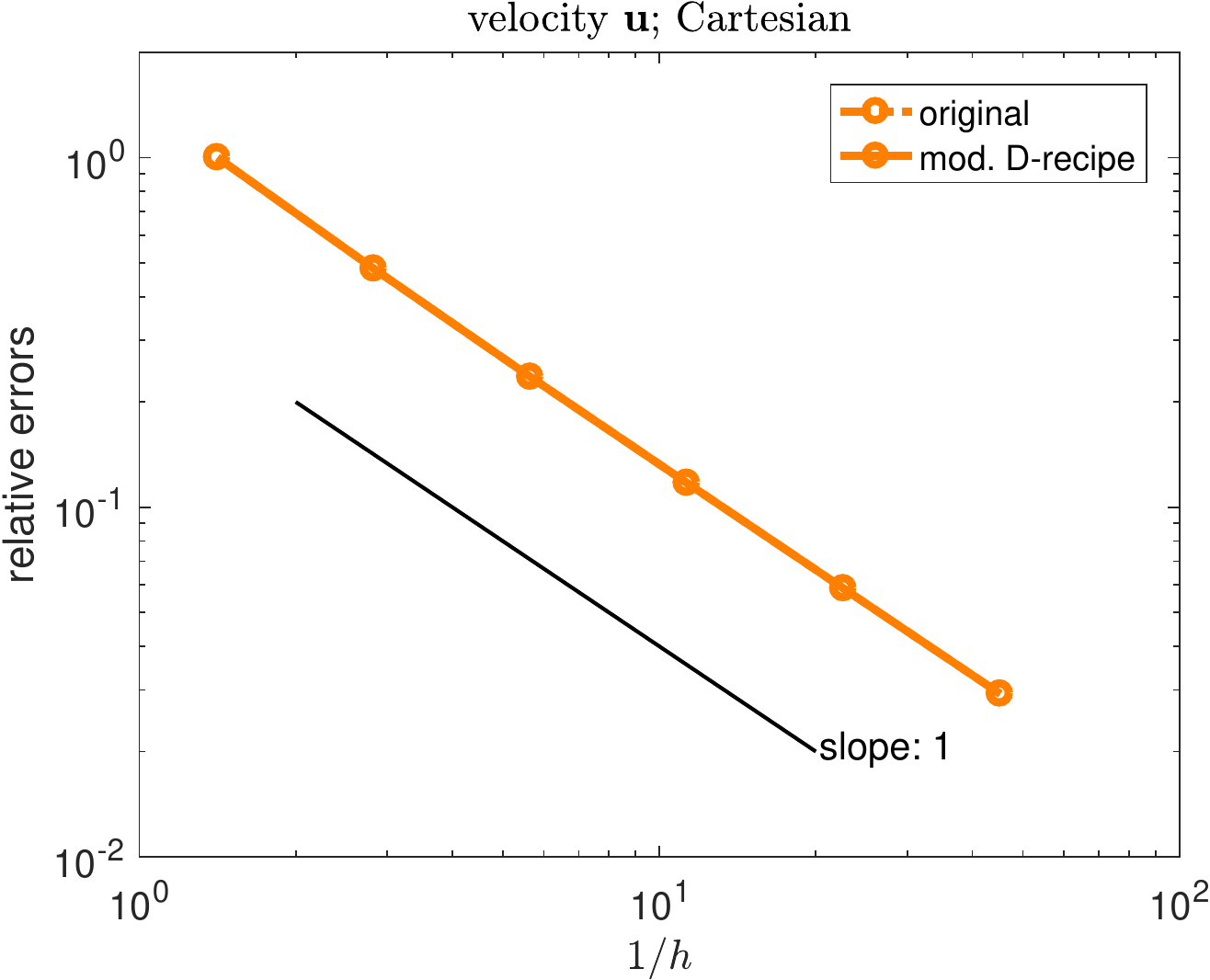}
\end{minipage}
\hfill
\begin{minipage}[t]{0.485\textwidth}
\centering
\includegraphics[width=\textwidth]{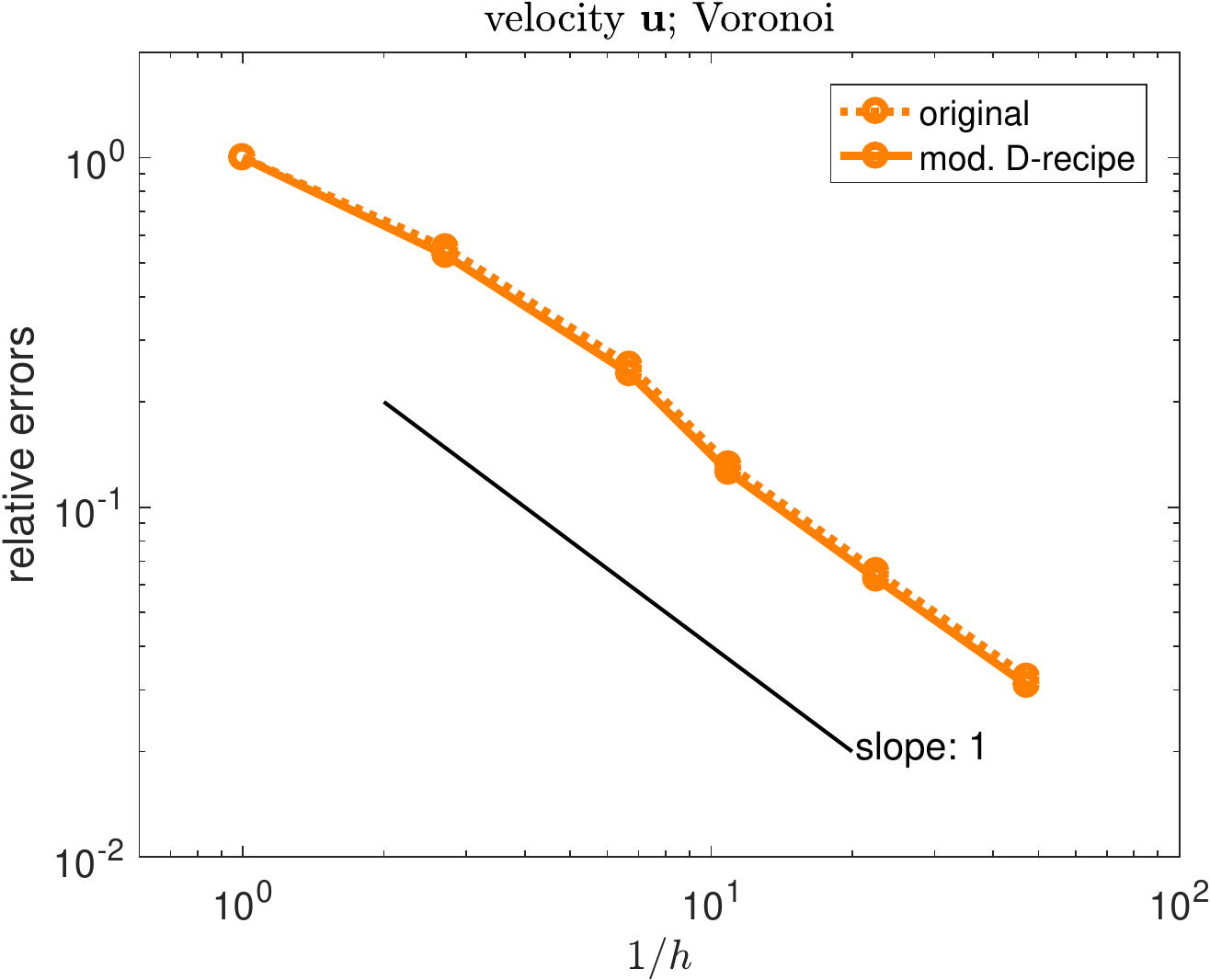}
\end{minipage}
\caption{Relative $L^2$ errors for the velocity field in example~1 at the final time $T$ on regular Cartesian meshes (left) and Voronoi meshes (right). The original stabilization~\eqref{eq:stabs} and the D-recipe stabilization~\eqref{eq:stabs_new} are employed.}
\label{fig:Ex1:u}
\end{figure}

\begin{figure}[h]
\begin{minipage}[t]{0.485\textwidth}
\centering
\includegraphics[width=\textwidth]{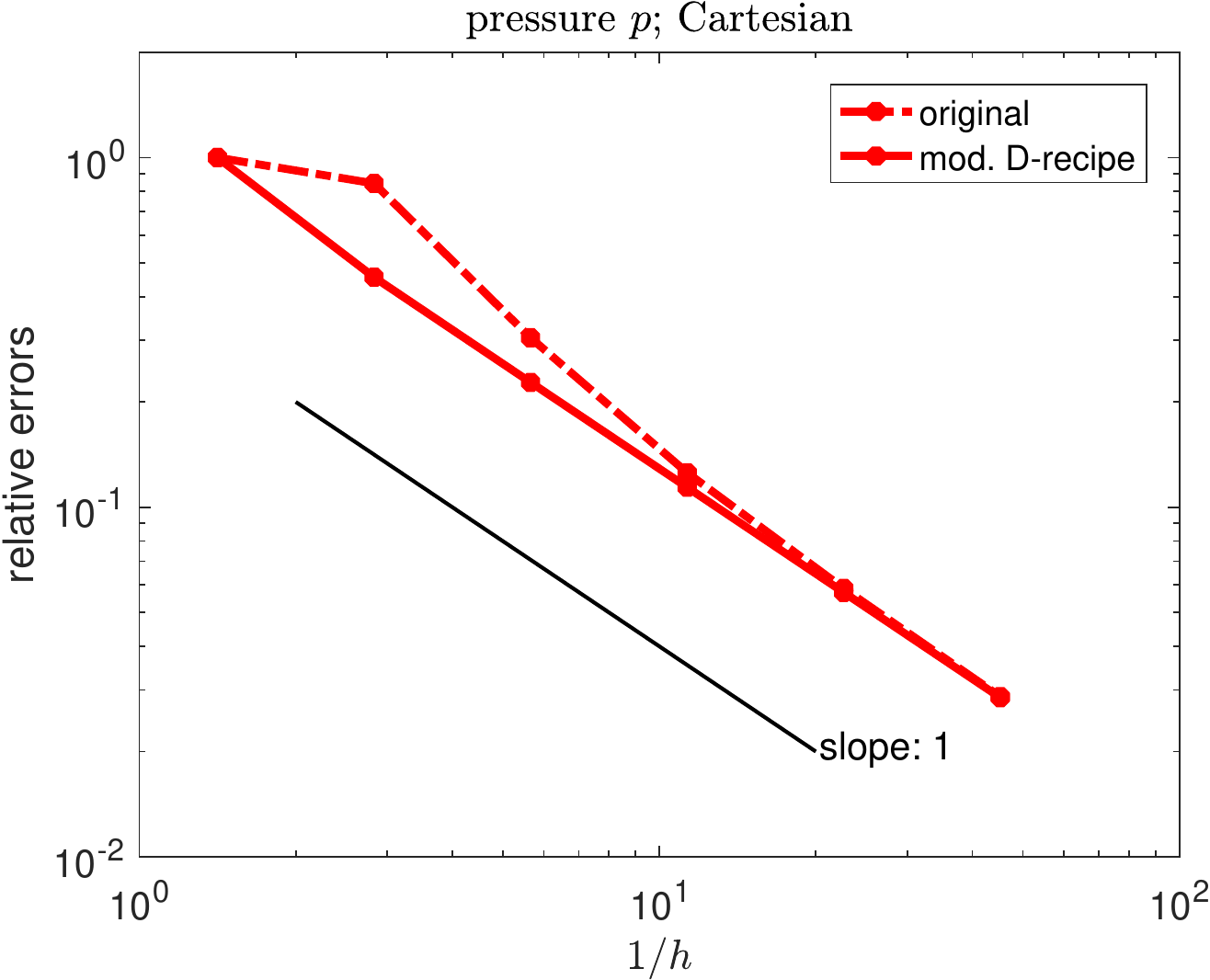}
\end{minipage}
\hfill
\begin{minipage}[t]{0.485\textwidth}
\centering
\includegraphics[width=\textwidth]{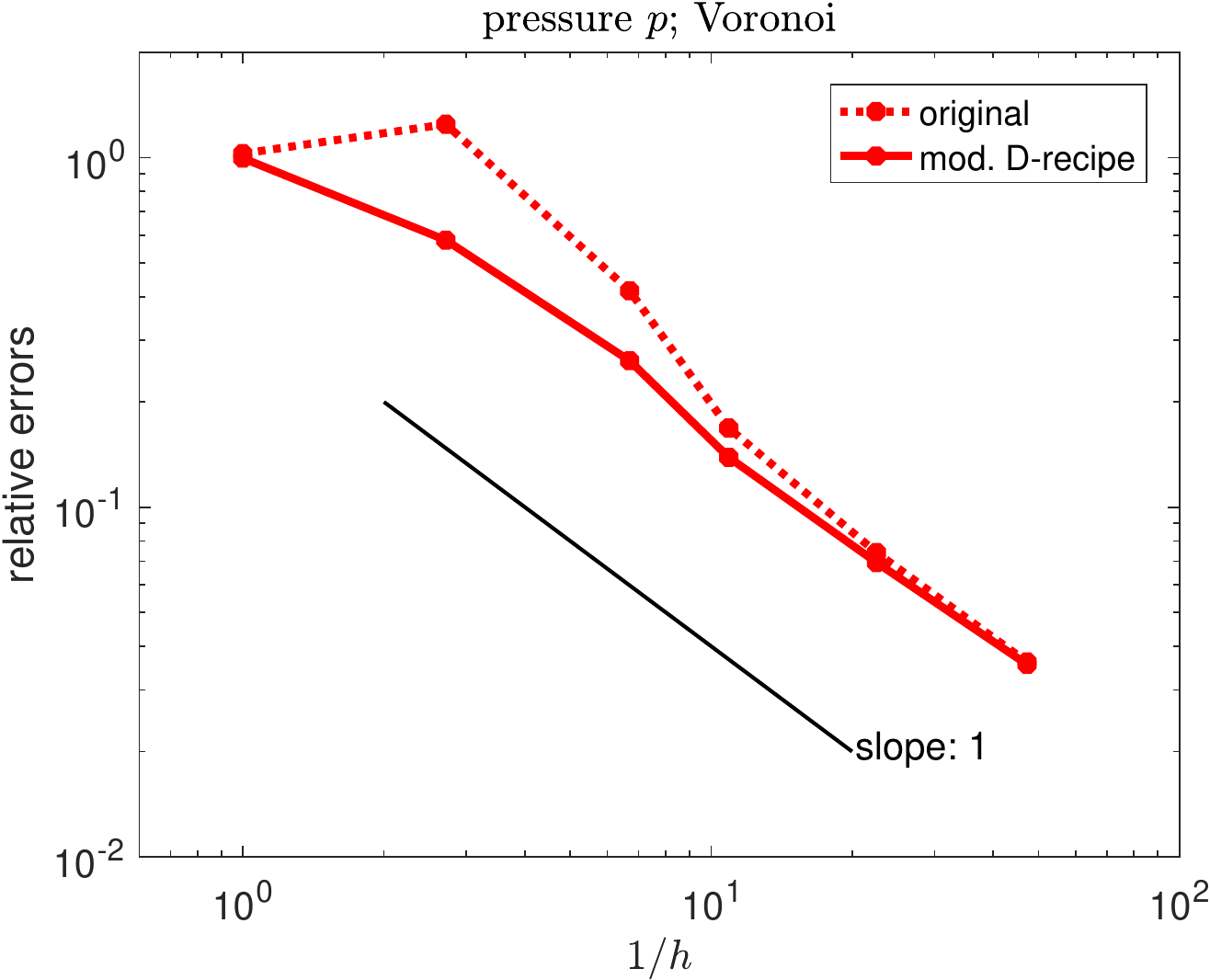}
\end{minipage}
\caption{Relative $L^2$ errors for the pressure in example~1 at the final time $T$ on regular Cartesian meshes (left) and Voronoi meshes (right). The original stabilization~\eqref{eq:stabs} and the D-recipe stabilization~\eqref{eq:stabs_new} are employed.}
\label{fig:Ex1:p}
\end{figure}

The relative $L^2$ discretization errors for the concentration are plotted in Figure \ref{fig:Ex1:c} in terms of the mesh size $h$ for both families of meshes and both variants of stabilizations. In order to better underline the expected linear convergence of the method both in $h$ and $\tau$ (see Theorem~\ref{thm:Cn_cn}, recalling that $k=0$), the time step~$\tau$ is chosen proportional to $h$. In other words, starting with the coarsest mesh and $\tau=T/5$, each subsequent case is obtained by dividing both $h$ (adopting a finer mesh) and $\tau$ by a factor of 2. Analogous plots are shown for the velocity and pressure variable errors in Figures~\ref{fig:Ex1:u} and~\ref{fig:Ex1:p}. In all cases, the linear convergence rates are in accordance with Theorem \ref{thm:Un_un_Pn_pn} and Theorem~\ref{thm:Cn_cn}. For the pressure discretization error, since the initial meshes are very coarse, we observe some pre-asymptotic regime when employing the original stabilizations in~\eqref{eq:stabs}. This effect, however, is not present for the alternative stabilizations in~\eqref{eq:stabs_new}. Both variants lead to similar results for the concentration and velocity errors.


Since the concentration often evolves more rapidly than the velocity and pressure, it could be worth to consider a cheaper variant of the discrete scheme~\eqref{eq:fully_discr_model_0}-\eqref{eq:fully_discr_model_2} where the discrete velocity-pressure pair is updated only every R time steps (with $R \in {\mathbb N}$). This leads to a smaller number of linear system resolutions (possibly with a small reduction in accuracy) since only system \eqref{eq:fully_discr_model_1} is solved at every time step, while \eqref{eq:fully_discr_model_2} is solved only every R steps.
In order to test this, we tried to run the same test above and compare the original version with the cheaper version with $R=5$. The difference in error was only at the fourth meaningful digit; we do not plot the graphs since these would completely overlap the ones of the original method.

% -----------------------------------------------------------------------------------------
\vspace{0.2cm}
\textit{Example 2:} Next, we investigate the behavior of the method for Test~1 and Test~2 in~\cite{wang2000approximation,chainais2007convergence}.

The problem is given in the form~\eqref{eq:model_problem} with boundary conditions~\eqref{eq:bdry_cond} and initial condition~\eqref{eq:initial_cond} over the spatial domain $\Omega=(0,1000)^2$ ft$^2$. Moreover, $T=3600$ days and $\tau=36$ days. At the upper right corner, i.e. at $[1000,1000]$, fluid with concentration $\chat=1.0$ is injected with rate $\qplus=30$ ft$^2$/day, whereas at the lower left corner, i.e. at $[0,0]$, material is absorbed with rate $\qminus=30$ ft$^2$/day. Both wells are henceforth treated as Dirac masses, which is admissible at the discrete level since the discrete functions are piecewise regular (which can be interpreted as an approximation of the Dirac delta by a localized function with support within the corner element and unitary integral). Furthermore, the following choices for the parameters are picked: $\phi=0.1$, $d_{\ell}=50$, $d_t=5$, $c_0=0$, $\bgamma(c)=0$, and $a(c)=80(1+(M^{\frac{1}{4}}-1)c)^4$, where 
\begin{equation*}
\text{Test A}: d_m=10, M=1; \qquad\qquad
\text{Test B}: d_m=0, M=41.
\end{equation*}

Whereas $a(c)$ is constant for Test A, it changes rapidly across the fluid interface for Test B (which is in fact not covered by the theoretical analysis since $d_m=0$, but is interesting to study numerically) resulting in a much faster propagation of the fluid concentration front along the diagonal direction ($d_\ell \gg d_t$). This effect is known as \textit{macroscopic fingering phenomenon}\cite{ewing1983mathematics}.

For this example, we used a regular 25x25 Cartesian mesh and we employed the more sophisticated stabilization in~\eqref{eq:stabs_new}. Since Test~B is highly convection-dominated, pure application of our method leads to local disturbances in the form of \textit{overshoots} and \textit{undershoots} of the numerical solution for the concentration, typical in the context of convection-dominated problem. To this purpose, for this test case, we employ the flux-corrected transport (FCT) algorithm with linearization \cite{FCT1, FCT}. 
The FCT scheme with linearization for convection-dominated flow problems  operates in two steps:
(1) advance the solution in time by a low-order overly diffusive scheme  to suppress spurious oscillations, (2) correct the solution using (linear) antidiffusive fluxes.
In that way the computed solution does not show spurious oscillations and layers are not smeared.

\begin{figure}[h]
\begin{minipage}[t]{0.485\textwidth}
\centering
\includegraphics[width=\textwidth]{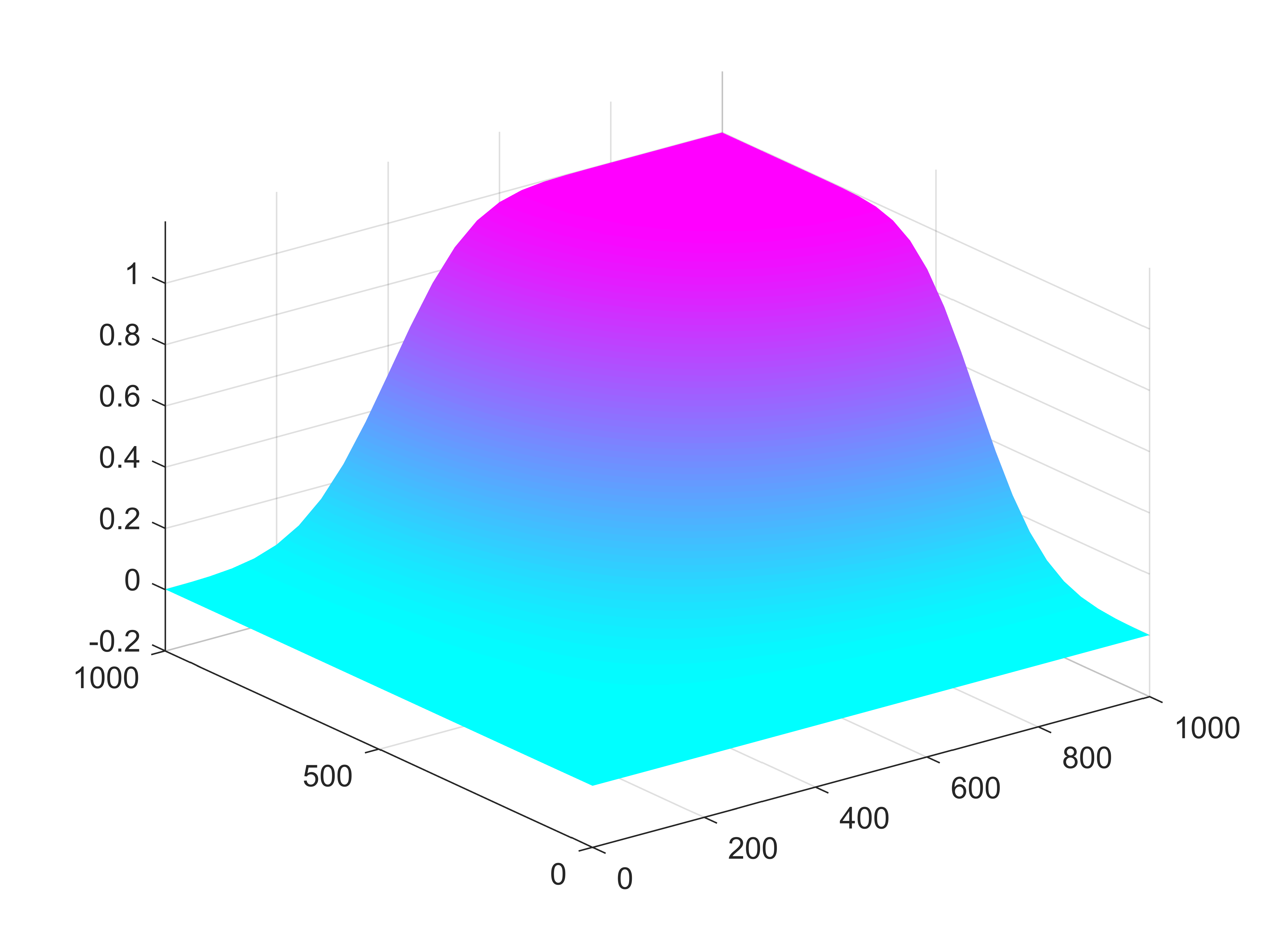}
\end{minipage}
\hfill
\begin{minipage}[t]{0.485\textwidth}
\centering
\includegraphics[width=\textwidth]{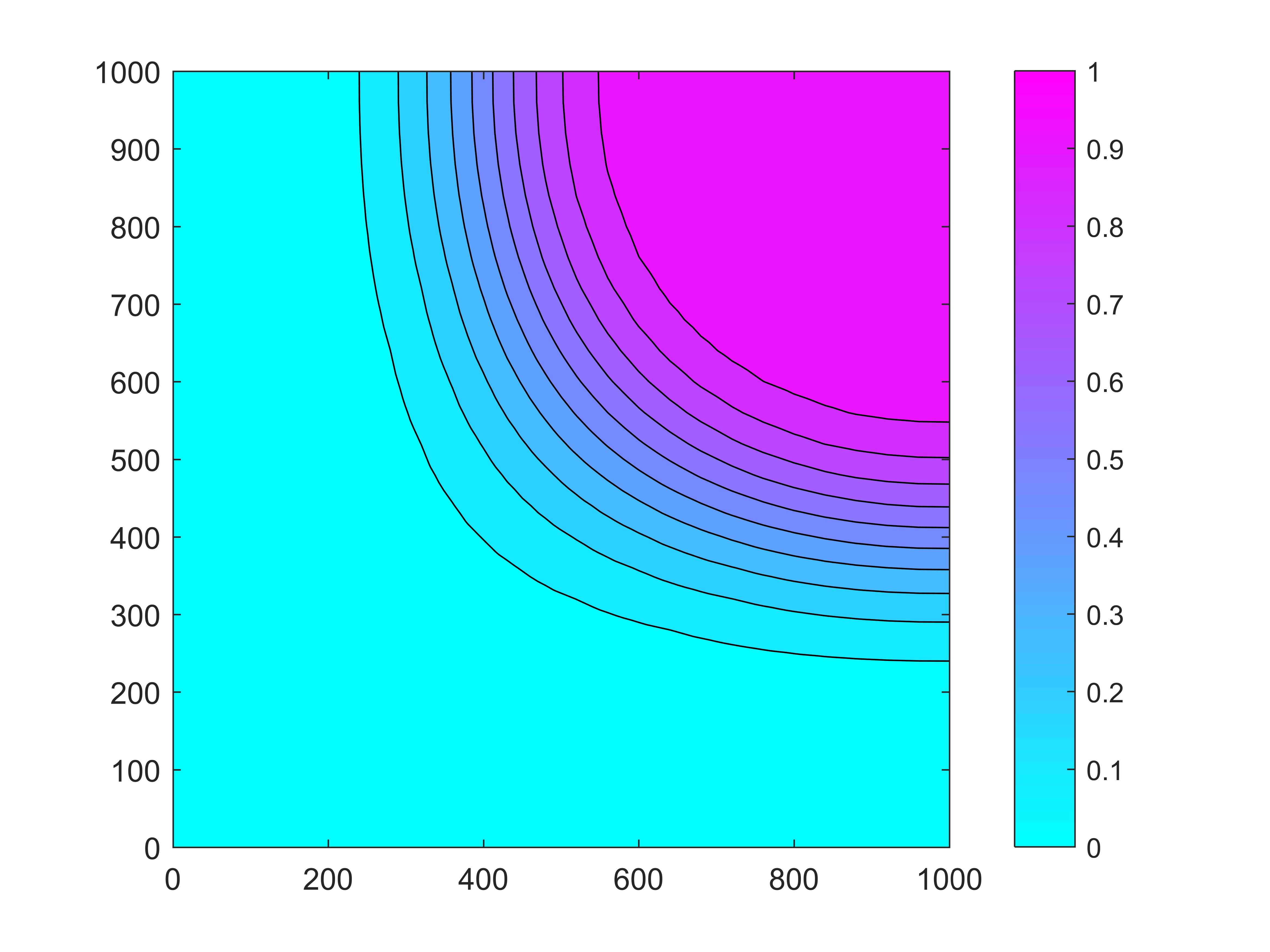}
\end{minipage}
\caption{Numerical solution for the concentration (left) and contour plot (right) after 3 years in Test~A.}
\label{fig:Ex2_testA_3}
\end{figure}

\begin{figure}[h]
\begin{minipage}[t]{0.485\textwidth}
\centering
\includegraphics[width=\textwidth]{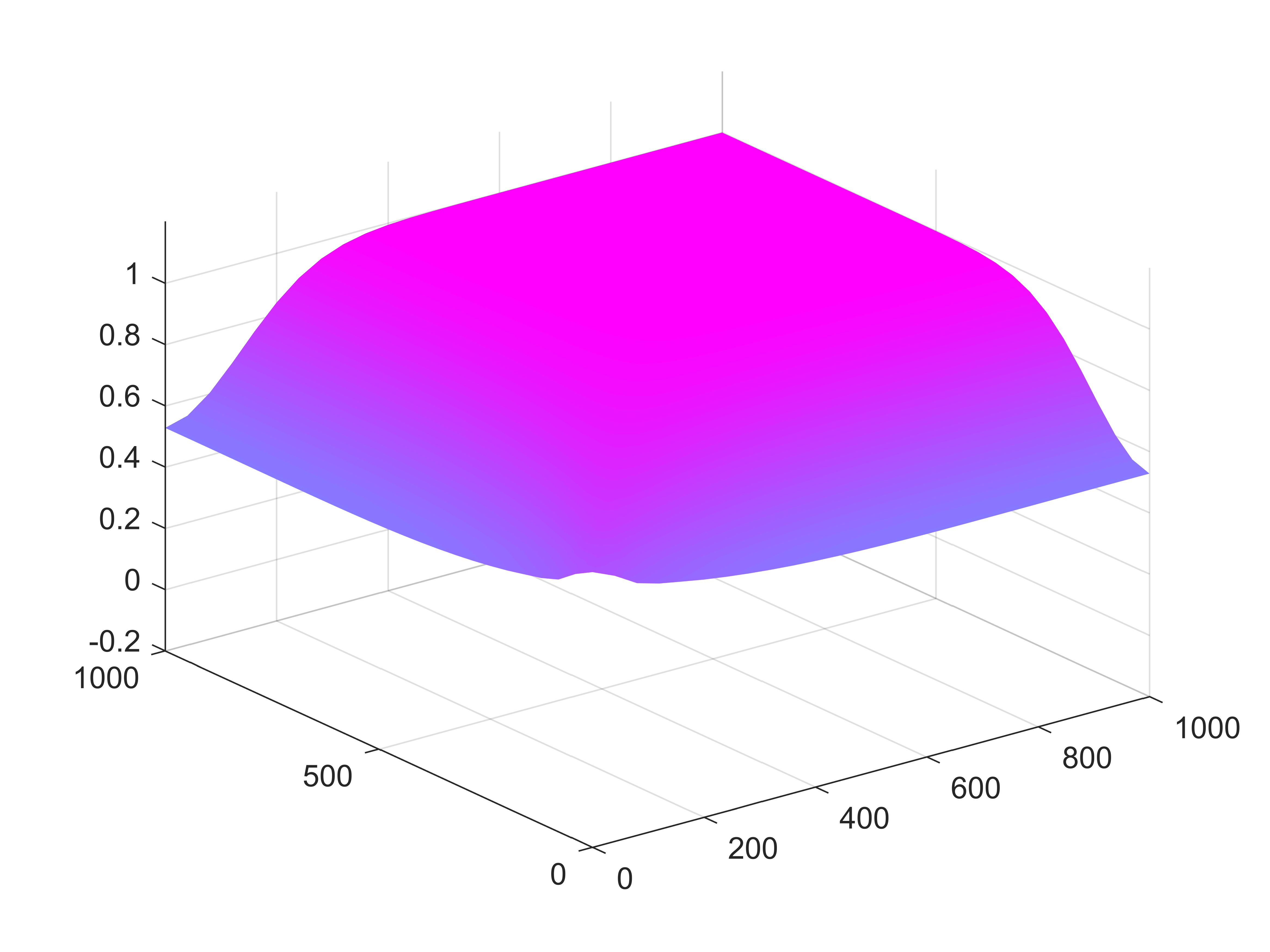}
\end{minipage}
\hfill
\begin{minipage}[t]{0.485\textwidth}
\centering
\includegraphics[width=\textwidth]{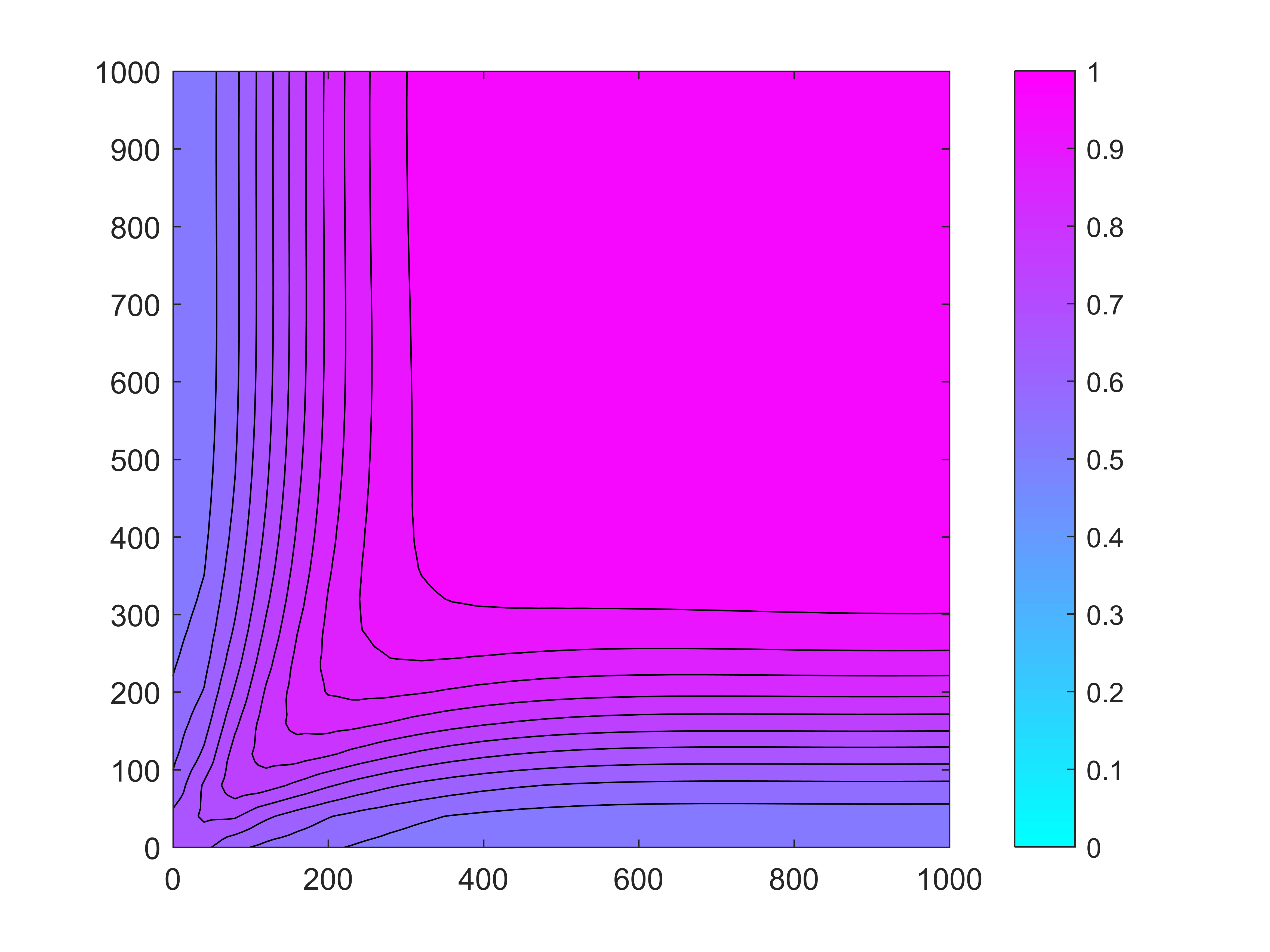}
\end{minipage}
\caption{Numerical solution for the concentration (left) and contour plot (right) after 10 years in Test~A.}
\label{fig:Ex2_testA_10}
\end{figure}

\begin{figure}[h]
\begin{minipage}[t]{0.485\textwidth}
\centering
\includegraphics[width=\textwidth]{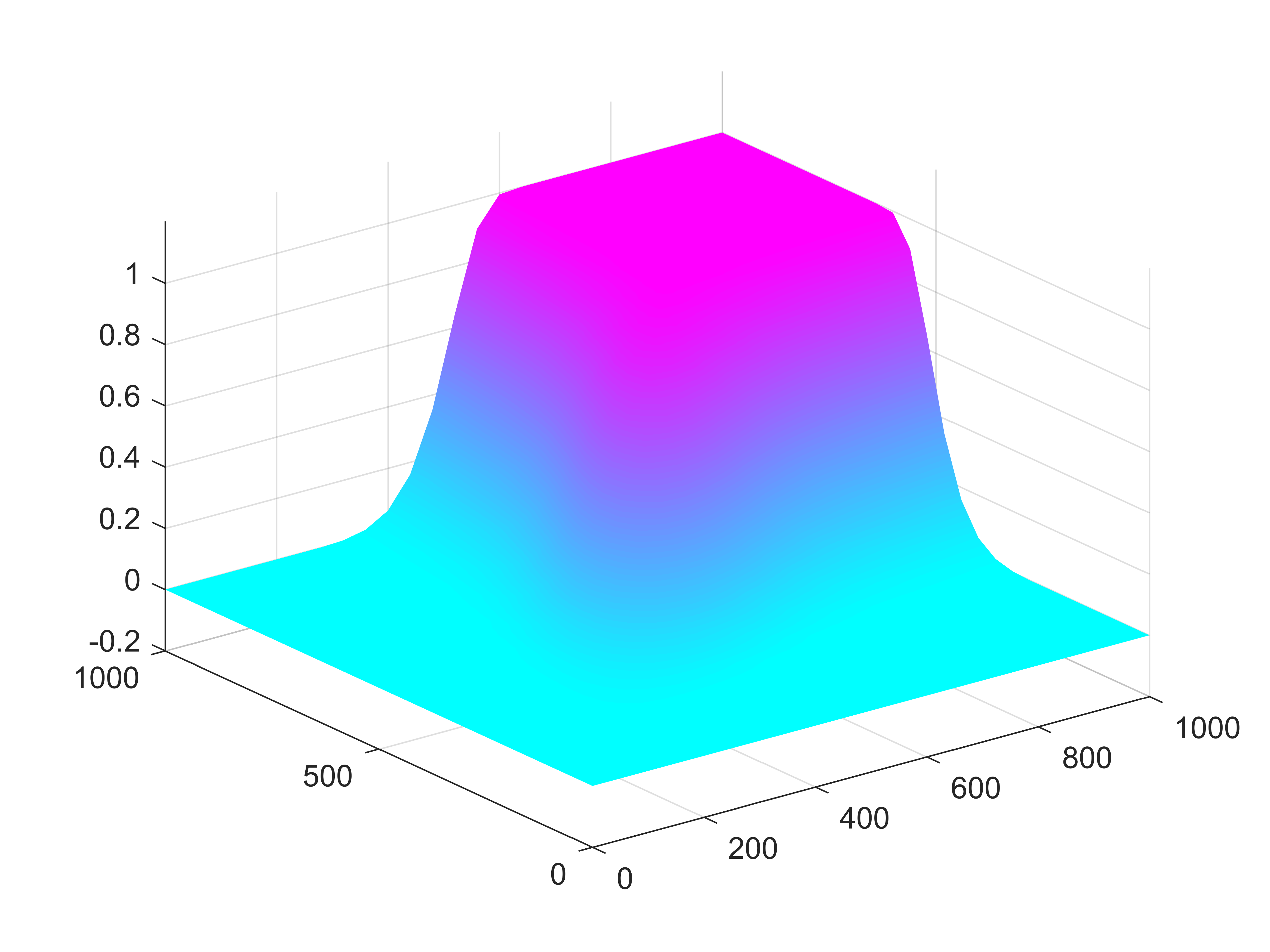}
\end{minipage}
\hfill
\begin{minipage}[t]{0.485\textwidth}
\centering
\includegraphics[width=\textwidth]{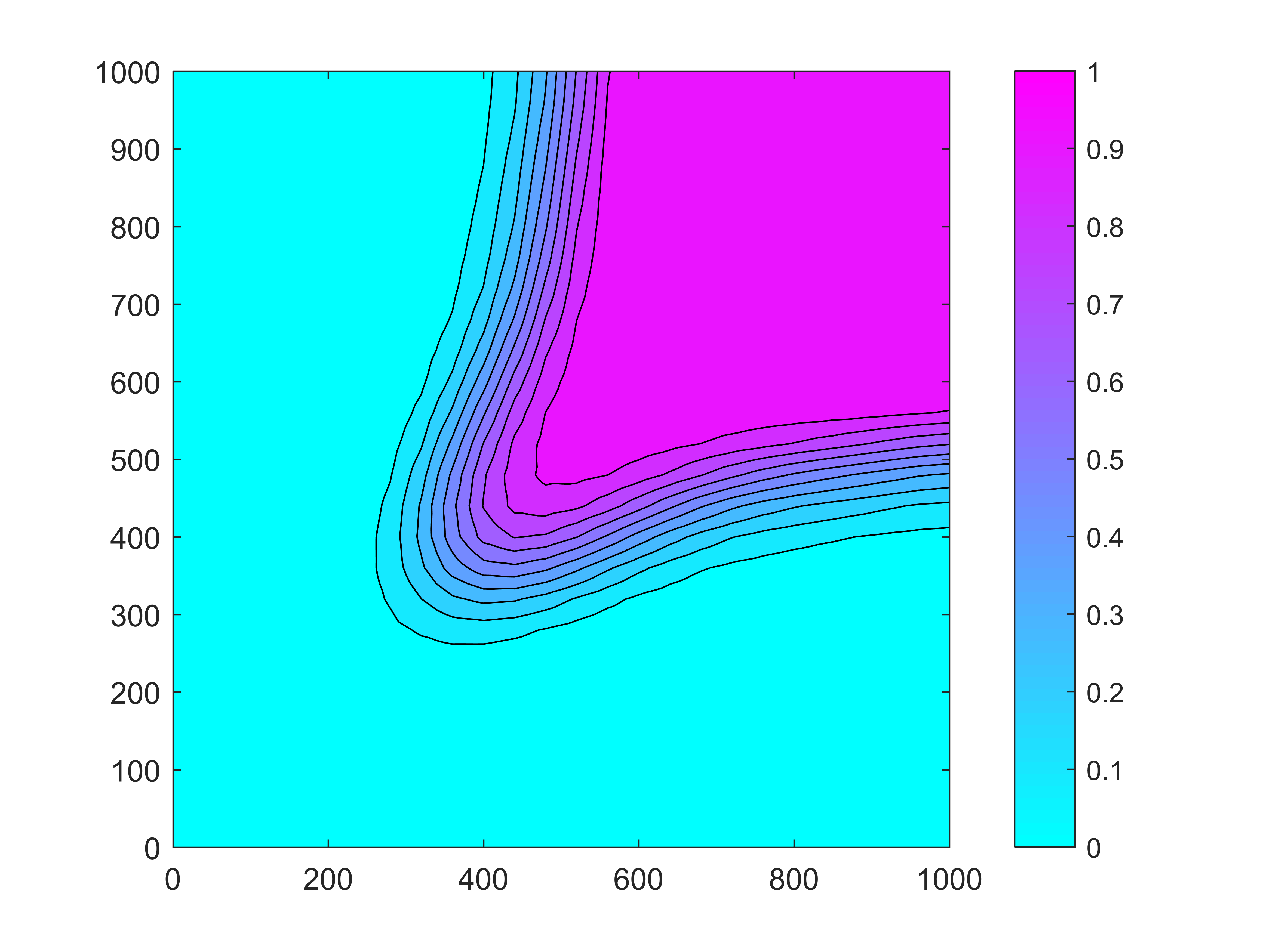}
\end{minipage}
\caption{Numerical solution for the concentration (left) and contour plot (right) after 3 years in Test~B.}
\label{fig:Ex2_testB_3}
\end{figure}

\begin{figure}[h]
\begin{minipage}[t]{0.485\textwidth}
\centering
\includegraphics[width=\textwidth]{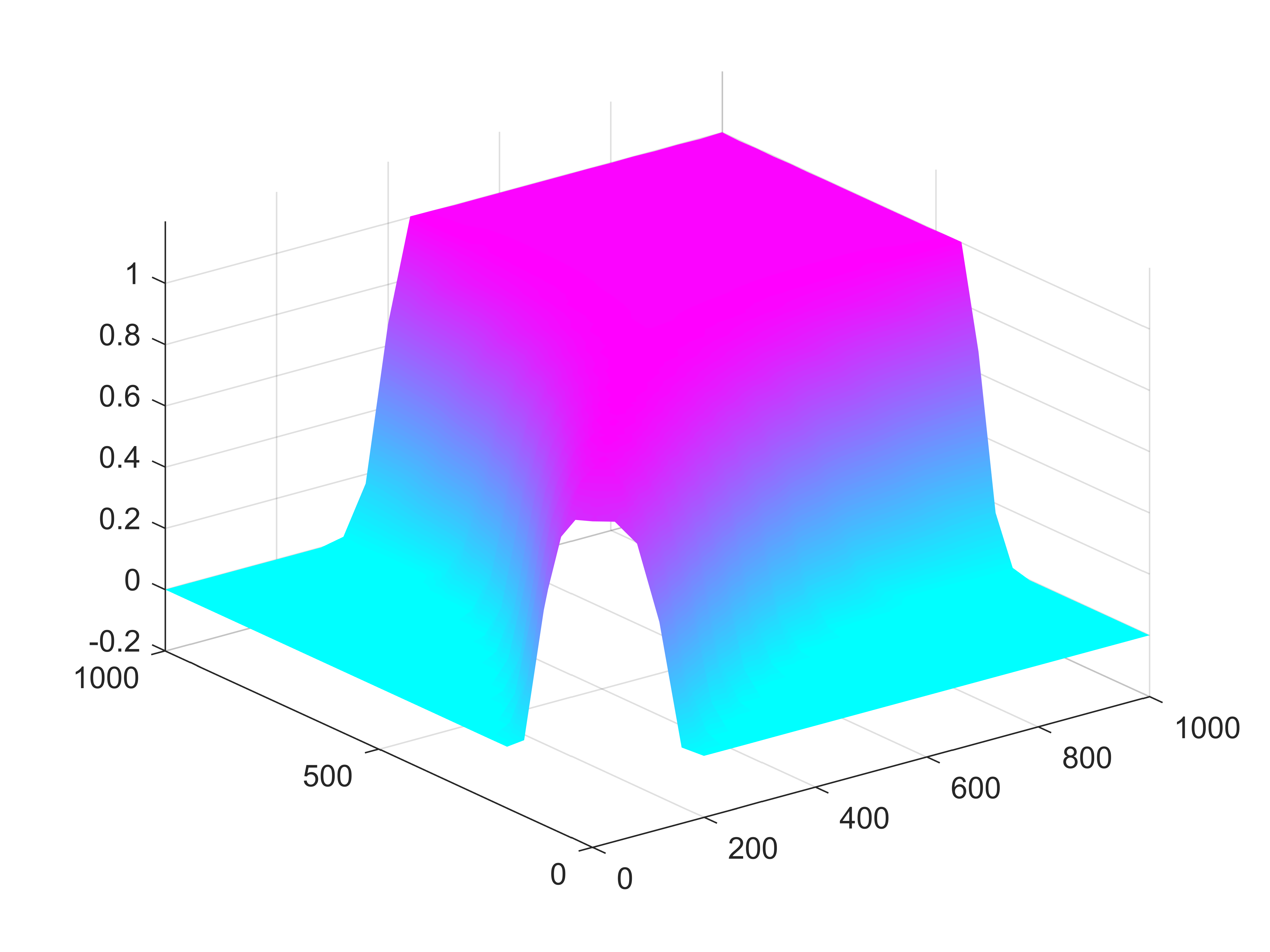}
\end{minipage}
\hfill
\begin{minipage}[t]{0.485\textwidth}
\centering
\includegraphics[width=\textwidth]{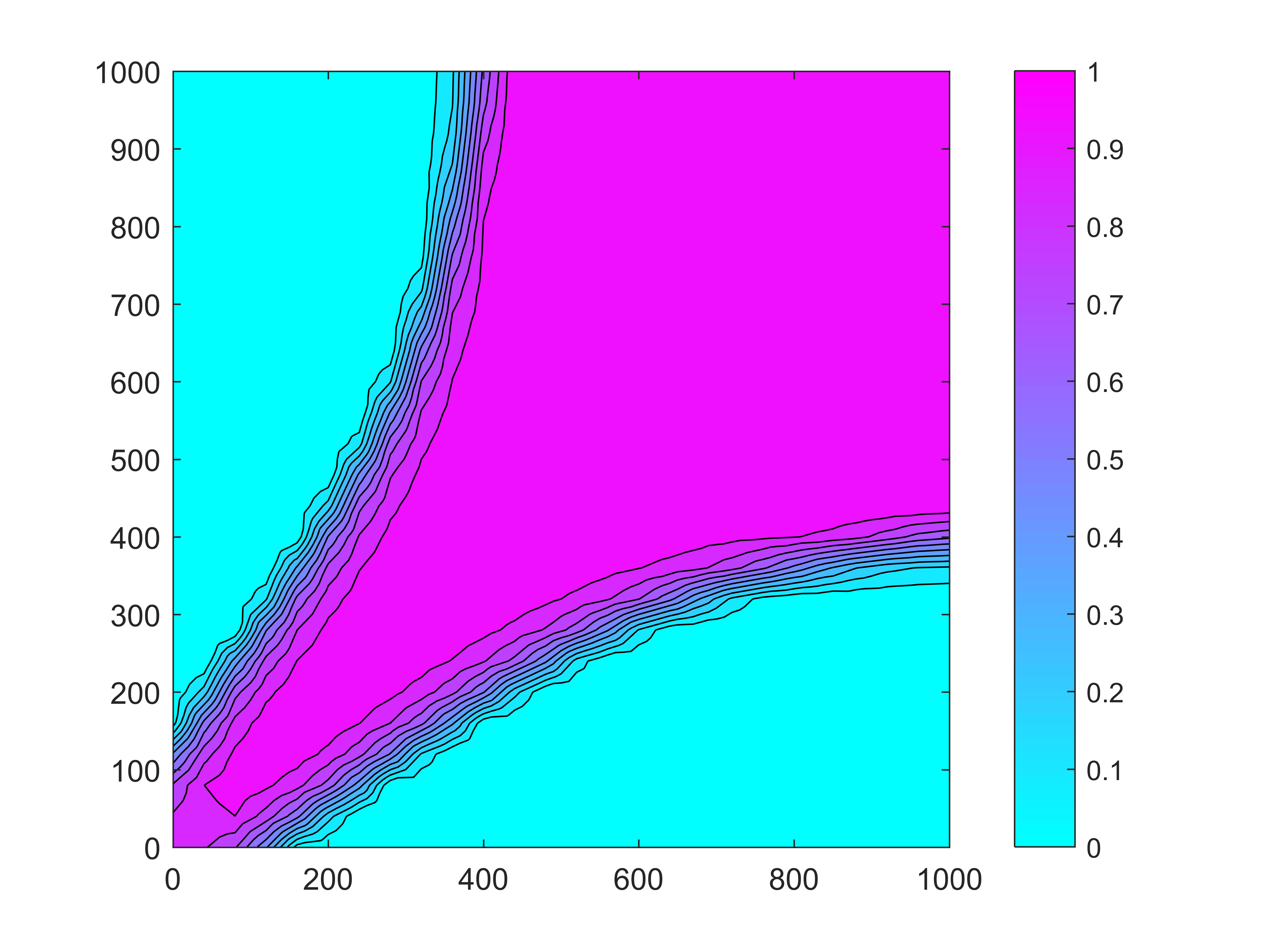}
\end{minipage}
\caption{Numerical solution for the concentration (left) and contour plot (right) after 10 years in Test~B.}
\label{fig:Ex2_testB_10}
\end{figure}

Due to the fact that no analytical solutions are available for Test~A and Test~B, we plot the numerical solutions (and the corresponding contour plots) for the concentration after 3 and 10 years. These times correspond to $n=30$ and $n=100$, respectively. For visualization of the results, since the numerical solution is virtual but the nodal values are known, we simply add, inside each square, the barycenter with associated mean value of the nodal values, then create a triangulation based upon these points, and finally interpolate the function values linearly inside each triangle. In Figures~\ref{fig:Ex2_testA_3} and~\ref{fig:Ex2_testA_10}, the results for Test~A are portrayed, and in Figures~\ref{fig:Ex2_testB_3} and~\ref{fig:Ex2_testB_10}, those for Test~B. The results are similar to those obtained in~\cite{wang2000approximation,chainais2007convergence}.

%\section{Conclusions}

\section*{Acknowledgements}
The first (L.B.d.V.) and last (G.V.) authors where partially supported by the European Research Council through the H2020 Consolidator Grant (grant no. 681162) CAVE, Challenges and Advancements in Virtual Elements. This support is gratefully acknowledged.
The second author (A.P.) has been funded by the Austrian Science Fund (FWF) through the project P~29197-N32, and by the Doctoral Program (DK) through the FWF Project W1245.

%%%%%%%%%%%%%%%%%%%%%%%%%%%%%%%%%%%%%%%%%%%%%%%%%%%%%%%%%%%%%%%%%%%%%%%%%%%
{\footnotesize
	\bibliography{bibliogr}
}
\bibliographystyle{plain}
%%%%%%%%%%%%%%%%%%%%%%%%%%%%%%%%%%%%%%%%%%%%%%%%%%%%%%%%%%%%%%%%%%%%%%%%%%%

%%%%%%%%%%%%%%%%%%%%%%%%%%%%%%%%%%%%%%%%%%%%%%%%%%%%%%
\end{document}